\documentclass[12pt]{article}
\usepackage{mathrsfs}
\usepackage{stmaryrd}
\usepackage{amsfonts,amsmath,amssymb,amscd}
\usepackage{shadow}
\usepackage{graphicx}
\usepackage{color}
\usepackage{longtable}
\usepackage{pstricks,multido}
\usepackage{hyperref}
\usepackage{qtree}
\usepackage{tikz}
\usetikzlibrary{arrows.meta}
\usepackage{subcaption}
\allowdisplaybreaks
\usetikzlibrary{shapes.geometric,shapes.misc,positioning}
\usetikzlibrary{arrows.meta,decorations.pathreplacing}

\newtheorem{thm}{Theorem}[section]
\newtheorem{lem}[thm]{Lemma}
\newtheorem{prop}[thm]{Proposition}
\newtheorem{cor}[thm]{Corollary}

\newtheorem{example}{Example}[section]

\def\Sym{\mathfrak{S}}
\def\N{\mathcal{N}^2}
\newcommand{\M}{\mathcal{M}}
\newcommand{\T}{\mathcal{T}}
\newcommand{\PP}{\mathcal{P}}
\newcommand{\K}{\mathcal{K}}
\newcommand{\D}{\mathcal{D}}
\newcommand{\Arc}{\operatorname{Arc}}
\newcommand{\rr}{\mathrm{r}}
\newcommand{\Lrminp}{\mathrm{Lrminp}}
\newcommand{\Lrminv}{\mathrm{Lrminv}}
\newcommand{\Rlmaxp}{\mathrm{Rlmaxp}}
\newcommand{\Rlmaxv}{\mathrm{Rlmaxv}}
\newcommand{\std}{\mathrm{std}}

\newcommand{\pos}{\mathrm{pos}}
\newcommand{\F}{\mathcal{F}}
\newcommand{\idr}{\mathrm{idr}}
\newcommand{\pre}{\mathrm{pre}}

\newcommand{\Des}{\mathrm{Des}}

\newcommand{\last}{\mathrm{last}}
\newcommand{\ro}{\mathrm{root}}
\newcommand{\RootPos}{\mathrm{RootPos}}
\newcommand{\core}{\mathrm{core}}
\newcommand{\jt}{\mathrm{jt}}
\newcommand{\pv}{\mathrm{pv}}
\newcommand{\SI}{\mathcal{SI}}
\newcommand{\RM}{\mathcal{RM}}
\newcommand{\WDes}{\mathrm{WDes}}
\newcommand{\E}{\mathcal{E}}
\newcommand{\B}{\mathcal{B}}
\newcommand{\Z}{\mathcal{Z}}
\newcommand{\R}{\mathcal{R}}
\newcommand{\C}{\mathcal{C}}
\newcommand{\Cut}{\mathrm{Cut}}
\newcommand{\Sep}{\mathrm{Sep}}

\newcommand{\KK}{\mathcal{K}}

\newcommand{\DD}{\mathcal{D}}

\def\qed{\hfill \rule{4pt}{7pt}}

\def\pf{\noindent {\it{Proof.} \hskip 2pt}}

\numberwithin{equation}{section}
\numberwithin{figure}{section}

\usepackage{tikz}

\tikzset{
	state arc/.style={
		line width=1.1pt,
		line cap=round,
		line join=round
	}
}

\newcommand{\statepic}[1]{%
	\begin{tikzpicture}[baseline=-0.55ex,x=1em,y=1em]
		#1
		\fill[black] (0,0) circle (1.7pt);
	\end{tikzpicture}%
}

\newcommand{\stateisolated}{%
	\statepic{}%
}

\newcommand{\stateredloop}{%
	\statepic{%
		\draw[state arc,red]
		(0.10,0.03)
		.. controls (0.42,0.85) and (-0.42,0.85)
		.. (-0.10,0.03);
	}%
}

\newcommand{\stateredout}{%
	\statepic{%
		\draw[state arc,red]
		(0.02,0.02)
		.. controls (0.20,0.70) and (0.62,0.70)
		.. (0.90,0.42);
	}%
}

\newcommand{\stateredin}{%
	\statepic{%
		\draw[state arc,red]
		(-0.90,0.42)
		.. controls (-0.62,0.70) and (-0.20,0.70)
		.. (-0.02,0.02);
	}%
}

\newcommand{\stateblueout}{%
	\statepic{%
		\draw[state arc,blue]
		(0.02,0.02)
		.. controls (0.20,0.70) and (0.62,0.70)
		.. (0.90,0.42);
	}%
}

\newcommand{\statebluein}{%
	\statepic{%
		\draw[state arc,blue]
		(-0.90,0.42)
		.. controls (-0.62,0.70) and (-0.20,0.70)
		.. (-0.02,0.02);
	}%
}

\newcommand{\stateredinblueout}{%
	\statepic{%
		\draw[state arc,red]
		(-0.90,0.42)
		.. controls (-0.62,0.70) and (-0.20,0.70)
		.. (-0.02,0.02);
		\draw[state arc,blue]
		(0.02,0.02)
		.. controls (0.20,0.70) and (0.62,0.70)
		.. (0.90,0.42);
	}%
}

\newcommand{\stateblueinredout}{%
	\statepic{%
		\draw[state arc,blue]
		(-0.90,0.42)
		.. controls (-0.62,0.70) and (-0.20,0.70)
		.. (-0.02,0.02);
		\draw[state arc,red]
		(0.02,0.02)
		.. controls (0.20,0.70) and (0.62,0.70)
		.. (0.90,0.42);
	}%
}

\begin{document}

\begin{center}
{\large\bf  Baxter Structures and the Odd-Even Enumeration of Red-Blue Motzkin Paths}
\end{center}

\begin{center}
Joanna N. Chen

School of Mathematical Sciences\\
Tianjin University of Technology\\
Tianjin 300384, P.R. China

chennajoanna@163.com\\
joannachen@tjut.edu.cn

\end{center}
\begin{abstract}
We study Motzkin paths satisfying the red-blue rule introduced by Wu, Eu, Ku and Shih. We answer their Question 8.5 and exhibit two new Baxter structures. Let $\M_m$ denote the set of such paths of
length $m$ and $a_m=|\M_m|$. 
We give a combinatorial proof of 
\[
a_{2n-1}=\mathrm{A005802}(n),\qquad
a_{2n}=\mathrm{A216947}(n) \qquad (n\geq 1).
\]
For the odd subsequence, we construct a bijection with
$1234$-avoiding permutations which sends the initial up and terminal down runs to the initial decreasing runs of a permutation and its inverse. The boundary 
correspondence arising from this
bijection
further identifies red-blue Motzkin paths of length $2n$ with $1234$-avoiding permutations of length $n+1$ in which $1$ occurs to the left of $2$.
We establish parallel binomial decompositions of this restricted permutation class and of $2$-colored noncrossing set partitions. The two decompositions have the same binomial coefficients, and their respective cores are explicitly in bijection with non-intersecting triples of lattice paths, a classical Baxter family. This establishes the even identity. 
Moreover, we construct a second bijection between odd-length red-blue Motzkin paths and $1234$-avoiding permutations. Under this bijection, splitting a path at its ground-level blue horizontal steps corresponds to the skew decomposition of the associated permutation.
\end{abstract}


\section{Introduction}

Motzkin paths are lattice paths with steps $U=(1,1)$,
$D=(1,-1)$, and $H=(1,0)$ that start and end on the horizontal axis and
never go below it. In this paper, 
steps are colored red or blue according to the parity of the number of horizontal
steps encountered. More precisely, the first horizontal step is red and
the colors of successive horizontal steps alternate; a non-horizontal
step is blue or red according to whether the number of preceding
horizontal steps is even or odd, respectively. Writing $U_b(w),D_b(w),H_b(w)$ (respectively $U_r(w),D_r(w),H_r(w)$)
for the set of the blue (respectively red) up, down and horizontal steps in $w$, we say a Motzkin path satisfies
the \emph{red-blue rule} if every prefix satisfies:
\begin{itemize}
  \item [1.] The number of red $D$ is at most that of blue  $U$;
  \item [2.] The number of blue $D$ is at most that of red $U$.
\end{itemize}
Let $\M_m$ be the set of red-blue Motzkin paths of length
$m$, and set $a_m=|\M_m|$. See the following as an example of the red-blue Motzkin path.

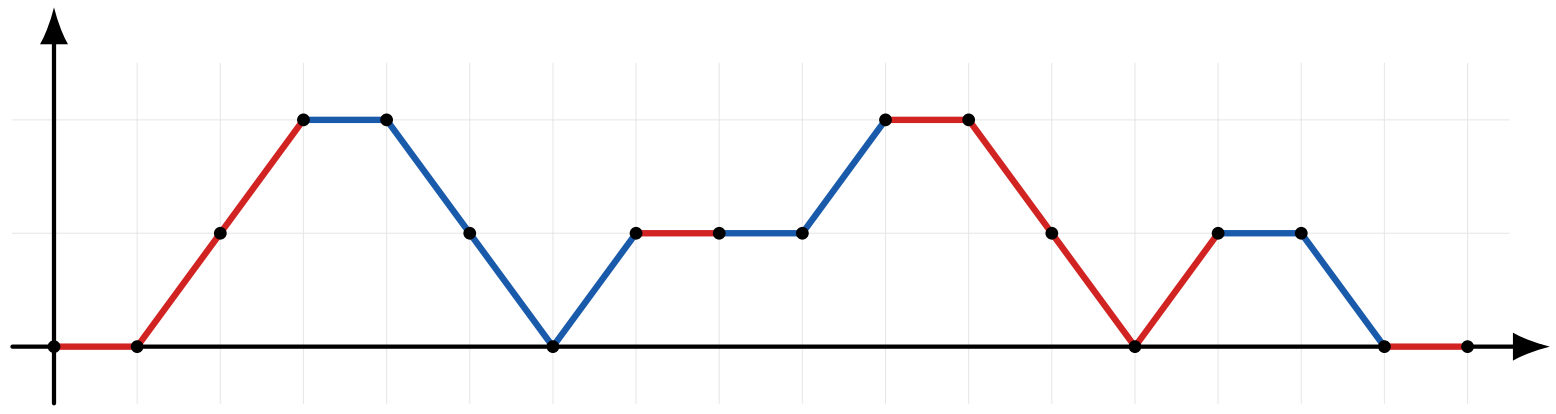
\begin{figure}[ht]
	\centering
	\begin{tikzpicture}[
		x=0.55cm,
		y=0.75cm,
		line cap=round,
		line join=round
		]
		\definecolor{pathred}{RGB}{210,35,35}
		\definecolor{pathblue}{RGB}{25,90,170}
		
		\draw[gray!18, very thin, step=1]
		(-0.5,-0.5) grid (17.5,2.5);
		
		\draw[-{Latex[length=2.5mm]}, thick]
		(-0.5,0) -- (18,0);
		\draw[-{Latex[length=2.5mm]}, thick]
		(0,-0.5) -- (0,3);
		
		\draw[pathred,  very thick] (0,0)  -- (1,0);   
		\draw[pathred,  very thick] (1,0)  -- (2,1);   
		\draw[pathred,  very thick] (2,1)  -- (3,2);   
		\draw[pathblue, very thick] (3,2)  -- (4,2);   
		\draw[pathblue, very thick] (4,2)  -- (5,1);   
		\draw[pathblue, very thick] (5,1)  -- (6,0);   
		\draw[pathblue, very thick] (6,0)  -- (7,1);   
		\draw[pathred,  very thick] (7,1)  -- (8,1);   
		\draw[pathblue, very thick] (8,1)  -- (9,1);   
		\draw[pathblue, very thick] (9,1)  -- (10,2);  
		\draw[pathred,  very thick] (10,2) -- (11,2);  
		\draw[pathred,  very thick] (11,2) -- (12,1);  
		\draw[pathred,  very thick] (12,1) -- (13,0);  
		\draw[pathred,  very thick] (13,0) -- (14,1);  
		\draw[pathblue, very thick] (14,1) -- (15,1);  
		\draw[pathblue, very thick] (15,1) -- (16,0);  
		\draw[pathred,  very thick] (16,0) -- (17,0);  
		
		\foreach \x/\y in {
			0/0,1/0,2/1,3/2,4/2,5/1,6/0,7/1,8/1,
			9/1,10/2,11/2,12/1,13/0,14/1,15/1,16/0,17/0
		}{
			\fill[black] (\x,\y) circle (1.2pt);
		}
	\end{tikzpicture}
	\caption{An example of red-blue Motzkin path with length $17$.}
	\label{fig:red-blue-path}
\end{figure}

The red-blue rule arose in the study of Cycle Petrie matrices.
Wu, Eu, Ku and Shih \cite{WuEuKuShih} showed that the associated paths satisfy this rule and
recorded the initial sequence
\[
1,1,2,3,6,11,23,47,103,225,513,1173,\ldots.
\]
They observed that its odd terms coincide with A005802 and its
even terms coincide with A216947, and posed the problem of finding a combinatorial proof of 
this fact  as given in 
Question 8.5 of \cite{WuEuKuShih}. To answer this question, we first recall some terminology.

Let $\Sym_n$ denote the set of all the permutations of $[n]=\{1,2,\ldots,n\}$. Given $\pi=\pi_1\cdots\pi_n\in\Sym_n$, 
define the \emph{initial decreasing run} of $\pi$ by the longest prefix
$\pi_1\pi_2\cdots\pi_\alpha$
whose consecutive entries are strictly decreasing.
We write $\idr(\pi)=\alpha-1$.
For $\sigma=\sigma_1\cdots\sigma_k\in\Sym_k$,
we say that $\pi$ contains a \emph{pattern} $\sigma$ if there exist indices
$1\leq i_1<i_2<\cdots<i_k\leq n$
such that $\pi_{i_a}<\pi_{i_b}$ if and only if $\sigma_a<\sigma_b$ for $1\leq a,b\leq k$.
In this case, the subsequence
$\pi_{i_1}\pi_{i_2}\cdots\pi_{i_k}$
is called an occurrence of $\sigma$ in $\pi$.
If $\pi$ contains no occurrence of $\sigma$, then $\pi$ is said to \emph{avoid}
the pattern $\sigma$. 
A \emph{vincular pattern} is a classical pattern together with a prescribed set
of adjacency requirements between consecutive entries of an occurrence.
More precisely, given a set $A\subseteq\{1,\ldots,k-1\}$, we say that
$\pi$ contains the vincular pattern $(\sigma,A)$ if $\pi_{i_1}\pi_{i_2}\cdots\pi_{i_k}$
is an instance of $\sigma$ in $\pi$
and, in addition, $i_{j+1}=i_j+1$ for every $j\in A$.
The adjacency requirements are usually indicated by underlining the corresponding adjacent letters.
 For example, $12\underline{43}$
means that the entries corresponding to $4$ and $3$ in the occurrence of $1243$ are required to
be adjacent. We write $\Sym_n(\sigma^{(1)},\ldots,\sigma^{(r)})$ as the set of all permutations in $\Sym_n$ that avoid each of the
patterns $\sigma^{(1)},\ldots,\sigma^{(r)}$ .
For further background on permutation patterns and pattern avoidance,
see Kitaev~\cite{Kitaev2011}.

A \emph{set partition} of $[n]=\{1,2,\ldots,n\}$ is a collection
$\Pi=\{B_1,\ldots,B_m\}$
of nonempty, pairwise disjoint subsets (called blocks) such that
$B_1\cup\cdots\cup B_m=[n].$
Given a block $B=\{b_1,b_2,\ldots,b_s\}$ with $b_1<b_2<\cdots<b_s$, its arcs are the pairs $(b_1,b_2),(b_2,b_3),\ldots,(b_{s-1},b_s)$, and we write $\Arc(\Pi)$ for the set of all arcs of $\Pi$. An arc $(i,i+1)$ is called an \emph{adjacent arc} of $\Pi$.
We often describe a set partition by drawing the graph whose vertex set is $[n]$ and whose edge set is $\Arc(\Pi)$. This graph is called the \emph{standard representation} of $\Pi$.
Two arcs $(i,k)$ and $(j,l)$ form a \emph{crossing} if
$i<j<k<l.$
A $2$-colored set partition of $[n]$, defined in \cite{Marberg2013}, is a pair
$(\Pi,\chi),$ where $\Pi$ is a set partition of $[n]$ and
$\chi:\operatorname{Arc}(\Pi)\longrightarrow\{\mathsf{r},\mathsf{b}\}$
is a coloring of its arcs with two colors, red and blue.
The pair $(\Pi,\chi)$ is called \emph{two-colored noncrossing} if no two
crossing arcs have the same color. 
We denote by $\N_n$ the set of all two-colored noncrossing set partitions of $[n]$.

Let $N=(0,1)$ and recall that $H=(1,0)$. 
A lattice path of length $n$ is a sequence of $n$ steps, each step being
either $N$ or $H$. Given a triple
$T=(P_{\mathrm{b}},P_{\mathrm{m}},P_{\mathrm{t}})$
of lattice paths of length $n$, we encode the step sequences of the three paths by
\[
b_1b_2\cdots b_n,\qquad
m_1m_2\cdots m_n,\qquad
t_1t_2\cdots t_n,
\]
where
\[
b_i,m_i,t_i=
\begin{cases}
1,&\text{if the corresponding step is }H,\\
0,&\text{if the corresponding step is }N.
\end{cases}
\]
We refer to $P_{\mathrm{b}}$ (resp. $P_{\mathrm{m}}$, $P_{\mathrm{t}}$) as the bottom (resp. middle, top) path of $T$.
For $n\geq0$, let $\PP_n$ be the set of such triples
of lattice paths of length $n$, starting respectively at
$(2,0),(1,1),(0,2)$ and satisfying that
\begin{equation}\label{eq:three-path-prefix}
\sum_{i=1}^{j}b_i
\geq
\sum_{i=1}^{j}m_i
\geq
\sum_{i=1}^{j}t_i,
\end{equation}
for each $0\leq j\leq n$ and
\begin{equation}\label{eq:three-path-endpoint}
\sum_{i=1}^{n}b_i
=
\sum_{i=1}^{n}m_i
=
\sum_{i=1}^{n}t_i.
\end{equation}
The inequalities in \eqref{eq:three-path-prefix} ensure that the three
paths remain strictly ordered from bottom to top, while
\eqref{eq:three-path-endpoint} ensures that their endpoints have the same
relative horizontal displacement as their starting points. 
An element of $\mathcal{P}_n$ will be called a non-intersecting triple of lattice paths.


Now we are ready to give the outline of the solution.
The two sequences have standard combinatorial interpretations. A005802 counts
$1234$-avoiding permutations \cite{GesselWeinsteinWilf1998} and 
A216947 counts two-colored noncrossing set partitions \cite{Marberg2013}. Thus, it suffices to  give 
combinatorial explanations of 
\[a_{2n-1}=|\Sym_n(1234)| \,\,\,\text{and}\,\,\, a_{2n}=|\N_n|.\]
We achieve this through the following theorems.

\begin{thm}\label{thm:bijodd}
For $n \geq 1$, we have
\[
\sum_{w\in\M_{2n-1}}
 x^{u_0(w)} y^{d_0(w)}
=\sum_{\pi\in\Sym_n(1234)}
 x^{\idr(\pi)}y^{\idr(\pi^{-1})},
\]
where $u_0(w)$ and $d_0(w)$ are the numbers of consecutive initial up steps and terminal down steps of $w$. Moreover, if $\pi=\phi(w)$, then we have 
$J_k(w)=P_k(\pi)$ and $T_k(w)=V_k(\pi)$ for $k=1,2,3$.
The sets $J_k(w)$, $T_k(w)$, $P_k(\pi)$, and $V_k(\pi)$ will be defined in Section \ref{sec:rb-bijection}.
\end{thm}

\begin{thm}\label{thm:bijeven} 
For $n \geq 1$,
we have
\[
\sum_{w\in\M_{2n}}
 x^{u_0(w)}
=\sum_{\pi\in\T_{n+1}}
 x^{\idr(\pi)},
\]
where $\T_{n+1}=\{\pi\in\Sym_{n+1}(1234): \pi^{-1}(1)<\pi^{-1}(2)\}$.
\end{thm}

The following equidistribution follows formally from Theorem \ref{thm:bijodd} and the natural decompositions of red-blue Motzkin paths and 1234-avoiding permutations
\[
\sum_{w\in\M_{2n-1}}
x^{\Sep(w)} y^{J_3(w)}z^{T_3(w)}
=\sum_{\pi\in\Sym_n(1234)}
x^{\Cut(\pi)}y^{P_3(\pi)}z^{V_3(\pi)},
\]
where  $\Sep(w)$ and $\Cut(\pi)$ will be defined in Section \ref{sec:bijodd-alternative}. However, the bijection $\phi$ in Theorem \ref{thm:bijodd} does not respect these  decompositions. The following theorem provides a direct bijection to explain this equidistribution.

\begin{thm}\label{thm:bijodd-alternative}
	For $n \geq 1$, there is an explicit bijection
	$F$ between $1234$-avoiding permutations and red-blue Motzkin path of odd length, with $|F(\pi)|=2|\pi|-1$,
    such that
	\[
	F(\sigma\ominus\tau)=F(\sigma)\,H_b\,F(\tau)
	\]
	for all nonempty $1234$-avoiding permutations $\sigma$ and $\tau$. Moreover, for every
	$\pi\in\Sym_n(1234)$,
	\[
	\Sep(F(\pi))=\Cut(\pi),\qquad
	J_3(F(\pi))=P_3(\pi),\qquad
	T_3(F(\pi))=V_3(\pi).
	\]	
	The skew sum $\ominus$ will be defined in Section \ref{sec:bijodd-alternative}.
\end{thm}

The numerical binomial-sum identity given in
Corollary \ref{cor:baxter-transform} appears in the conference poster by Gil, Kenepp and Weiner
\cite{GilKeneppWeiner2020}. See also the results of Bostan, Tirrell, Westbury and
Zhang~\cite{BostanTirrellWestburyZhang} on sequences
associated with invariant theory.
In the present work, we 
develop explicit combinatorial decompositions
for the corresponding classes of permutations
and set partitions, and identify two classes of cores enumerated by Baxter numbers.

\begin{thm}\label{thm:perdecom}
For $n \geq 1$,
we have
\[|\T_{n+1}|=\sum_{k=0}^{n-1}\binom{n-1}{k} |\K_k|,\]
where $\K_k=\{\pi\in\Sym_{k+2}(1234,12\underline{43}):\pi_{k+2}=2\}$.
\end{thm}

\begin{thm}\label{thm:percore-baxter}
For $n \geq 1$,
there is a bijection between 
 $\K_n$ and $\PP_{n}$.
\end{thm}

\begin{thm}\label{thm:setpardecom}
For $n \geq 1$,
we have
\[|\N_{n}|=\sum_{k=0}^{n-1}\binom{n-1}{k} |\D_k|,\]
where
$\D_k=\{\delta\in\N_{k+1}:\delta\text{ has no blue adjacent arcs}\}$.
\end{thm}

\begin{thm}\label{thm:setcore-baxter}
For $n \geq 1$,
there is a bijection between 
 $\D_n$ and $\PP_n$.
\end{thm}

 \begin{cor}\label{cor:baxter-transform}
	For $n\geq 1$, we have
	\[
	a_{2n}=|\T_{n+1}|=|\N_n|=\sum_{k=0}^{n-1}\binom{n-1}{k}B_{k+1},
	\]
	where $B_m$ denotes the number of Baxter permutations of length $m$.
\end{cor}

The rest of the paper is organised as follows. Section
\ref{sec:rb-bijection} establishes Theorems \ref{thm:bijodd} and
\ref{thm:bijeven} via  bijections between red-blue Motzkin
paths and $1234$-avoiding permutations.
Section \ref{sec:bijodd-alternative} constructs a second bijection compatible with skew decomposition and proves the set-valued refinement in Theorem \ref{thm:bijodd-alternative}.
Section \ref{sec:1234-core-baxter} proves Theorems
\ref{thm:perdecom} and \ref{thm:percore-baxter} using a core
decomposition of  restricted $1234$-avoiding permutations and a bijection between
these cores and non-intersecting triples of lattice paths. Section \ref{sec:partition-decomposition-baxter} proves Theorem
\ref{thm:setpardecom} by applying a blue-chain decomposition to
two-colored noncrossing set partitions, and proves Theorem
\ref{thm:setcore-baxter} by constructing a bijection between the
resulting colored cores and non-intersecting triples of lattice paths.

\section{Proofs of Theorems \ref{thm:bijodd} and \ref{thm:bijeven}}
\label{sec:rb-bijection}

This section gives a bijection between red-blue Motzkin
paths of length $2n-1$ and $1234$-avoiding permutations of length $n$.
The idea is two compressed
words extracted from the path recover the three rank levels of
the corresponding permutation. 

For a path $w=w_1w_2\cdots w_m\in\mathcal M_m,$
define two prefix sequences $x(w)=x_0x_1x_2\cdots x_m$ and $y(w)=y_0y_1y_2\cdots y_m$ by 
letting $x_0=y_0=0$ and 
\begin{equation}\label{eq:xy-prefix}
\begin{aligned} 
 x_i&=|\{1 \leq j \leq i\colon j \in U_b(w)\}|-|\{1 \leq j \leq i\colon j \in D_r(w)\}|,\\[2pt]
 y_i&=|\{1 \leq j \leq i\colon j \in U_r(w)\}|-|\{1 \leq j \leq i\colon j \in D_b(w)\}|,
\end{aligned}
\end{equation}
for $1 \leq i \leq m$.
 By the red--blue rule, we have 
$ x_i\geq0$ and $ y_i\geq0$ for each $0\leq i\leq m$.
Since $w$ is a Motzkin path, then $x_m+y_m=0$.
Since both $x_m$ and $y_m$ are nonnegative, they must both be zero.
Consequently, we have $|U_b(w)|=|D_r(w)|$ and $|U_r(w)|=|D_b(w)|$.

Now suppose that $w$ has length $2n-1$. Since the number of horizontal steps in $w$ is odd, we write it as $2h+1$. Under our coloring rule, the horizontal steps are
colored alternately, beginning and ending with red. It follows that
$|H_b(w)|=h$ and $|H_r(w)|=h+1$. Since $2|U_b(w)|+2|U_r(w)|+2h+1=2n-1$, then 
$|U_b(w)|+|U_r(w)|+h+1=n$.

For a permutation $\pi=\pi_1\cdots\pi_n$, define the \emph{rank}
of its $i$-th entry by
\begin{equation*}\label{eq:rb-rank}
 \rr_i(\pi)=1+\max\{\rr_j(\pi):j<i,\ \pi_j<\pi_i\}
\end{equation*}
with $\max\varnothing=0$.
Namely, $\rr_i(\pi)$ is the length of a longest increasing subsequence ending at
$\pi_i$ in $\pi$. We call an entry $\pi_i$ a \emph{$j$-rank entry} 
if $\rr_i(\pi)=j$. Clearly,
a permutation $\pi$ belongs to $\Sym_n(1234)$ if and only if
$\rr_i(\pi)\in\{1,2,3\}$ for each $i$. Define its corresponding position and value set by
\begin{equation*}\label{eq:rb-ranksets}
 P_k(\pi)=\{i:\rr_i(\pi)=k\} 
\,\,\,\,\text{and}\,\,\,\,
 V_k(\pi)=\{\pi_i:\rr_i(\pi)=k\}
\end{equation*}
for $k=1,2,3$.

Given $\pi \in \Sym_n$,  a  \emph{left-to-right minimum} of $\pi$ is a letter $\pi_i$ such that either
$i = 1$ or $\pi_i<\pi_j$ for all $j<i$, while a \emph{right-to-left maximum} of $\pi$ is a letter $\pi_i$ such that either
$i = n$ or $\pi_i>\pi_j$ for all $j>i$. 
Let $\Lrminp(\pi)$ and $\Lrminv(\pi)$ denote the set of the positions and values of the 
left-to-right minima of $\pi$, respectively. $\Rlmaxp(\pi)$ and $\Rlmaxv(\pi)$ can be define similarly. We have the 
following proposition.
\begin{prop}\label{prop:rankeq}
For $\pi \in \Sym_n(1234)$, entries belonging to the same rank of $\pi$ decrease from left to right. Moreover, we have
\begin{align*}
P_1(\pi)&=\Lrminp(\pi),\,\,\,V_1(\pi)=\Lrminv(\pi),\\
P_3(\pi)&\subseteq\Rlmaxp(\pi),\,\,\,V_3(\pi)\subseteq\Rlmaxv(\pi).  
\end{align*}
\end{prop}

We now construct two maps
$\phi:\M_{2n-1}\rightarrow\Sym_n(1234)$ and 
$g:\Sym_n(1234)\rightarrow\M_{2n-1}$.
Their well-definedness will be  postponed to
Lemmas \ref{lem:rb-path} and \ref{lem:rb-rank}.

Let $w\in\M_{2n-1}$. We associate with $w$ two words
$p=p_1\cdots p_n$ and $q=q_1\cdots q_n$ over $\{1,2,3\}$, and write $\jt(w)=(p,q)$.
Their reading rules are
\[
\begin{array}{c|c|c}
 &\text{position word }p_2\cdots p_{n}&\text{value word }q_2 \cdots q_{n}\\ \hline
\text{direction}&\text{left to right}&\text{right to left}\\
\text{steps}&U_b,H_b,D_b&D_r,H_b,U_r\\
\text{letters}&U_b\mapsto1,\ H_b\mapsto2,\ D_b\mapsto3
&D_r\mapsto1,\ H_b\mapsto2,\ U_r\mapsto3
\end{array}
\]
Let $p_1=q_1=1$. Write $|s|_i$ for the number of $i$'s in the word $s$.  Clearly, we have 
\begin{equation}\label{eq:rb-content}
 |p|_1=|q|_1=|U_b(w)|+1,
 \qquad
 |p|_2=|q|_2=|H_b(w)|,
 \qquad
 |p|_3=|q|_3=|U_r(w)|.
\end{equation}
For $k\in\{1,2,3\}$, write
\[J_k(w)=\{i:p_i=k\}
=\{i_1^{(k)}<\cdots<i_{s_k}^{(k)}\},
\quad 
T_k(w)=\{j:q_j=k\}
=\{v_1^{(k)}<\cdots<v_{s_k}^{(k)}\}.
\]
Define $\pi=\phi(w)$ by $\pi_{i_t^{(k)}}=v_{s_k+1-t}^{(k)}$,
where $k\in\{1,2,3\}$ and $1\leq t\leq s_k$ and write such $\pi$ as $\pv(p,q)$.

\begin{example}
	Consider $
	w=H_rU_rU_rH_bD_bD_bU_bH_rH_bU_bH_rD_rD_rU_rH_bD_bH_r
	\in\M_{17}$ given in Figure \ref{fig:red-blue-path}.
	We have $p=123312123$, $q=123112233$ and
\[
\begin{aligned}
	J_1(w)
	&=\{1,5,7\},
	&\quad
	T_1(w)
	&=\{1,4,5\},\\
	J_2(w)
	&=\{2,6,8\},
	&
	T_2(w)
	&=\{2,6,7\},\\
	J_3(w)
	&=\{3,4,9\},
	&
	T_3(w)
	&=\{3,8,9\}.
\end{aligned}
\]
Hence, $\pi=\phi(w)=579846123$.
\end{example}

We next construct a map $g$ from $\Sym_n(1234)$ to $\M_{2n-1}$. Given $\pi \in \Sym_n(1234)$, we first define
two words $p$ and $q$ by letting $p_i=k$ if and only if $i \in P_k(\pi)$ and $q_i=k$ if and only if $i \in V_k(\pi)$ for $1 \leq i \leq n$ and $k \in \{1,2,3\}$. 
It is not hard to check that $(p,q)=\pv^{-1}(\pi)$.
Let $h$ be the number of letters $2$ in either word. Write
\[
 p=1A_0\,2A_1\,2\cdots2A_h,
 \qquad
 q=1C_0\,2C_1\,2\cdots2C_h,
\]
where each $A_i$ and $C_i$ is a possibly empty word over $\{1,3\}$.
Set $B_i$ to be the word arising from $A_i$ be replacing $1$ with $U$ and $3$ with $D$. Set  $R_i$ to be the word obtained from $C_{h-i}$ be first reversing $C_{h-i}$ and then replacing $1$ with $D$ and $3$ with $U$.
Define $w=g(\pi)$ by
\[
 w=(B_0HR_0)\,H\,(B_1HR_1)\,H\cdots H\,(B_hHR_h).
\]
The horizontal step inside each parenthesized block is red, and the
horizontal steps separating consecutive blocks are blue. The steps in
$B_i$ are blue and those in $R_i$ are red. If $h=0$, the expression is
simply $B_0HR_0$.  In this process, it is easy to verify that $w=\jt^{-1}(p,q)$.

The length of $w$ can be computed as $2(n-1-h)+(2h+1)=2n-1$.
Equality of the contents of $p$ and $q$ shows that
$|U_b(w)|=|D_r(w)|$ and $|U_r(w)|=|D_b(w)|$. 

\begin{example}
	Given $\pi=579846123$, we have 
		\[
	\begin{aligned}
		P_1(\pi)
		&=\{1,5,7\},
		&\quad
		V_1(\pi)
		&=\{1,4,5\},\\
		P_2(\pi)
		&=\{2,6,8\},
		&
		V_2(\pi)
		&=\{2,6,7\},\\
		P_3(\pi)
		&=\{3,4,9\},
		&
		V_3(\pi)
		&=\{3,8,9\}.
	\end{aligned}
	\]
	It follows that $p=123312123$ and $q=123112233$. Write
	$p=1A_0\,2 A_1\, 2 A_2\, 2 A_3$ and $q=1C_0\,2 C_1\, 2 C_2 \,2 C_3$,
	where $A_0=\emptyset$, $A_1=331$, $A_2=1$, $A_3=3$, $C_0=\emptyset$, $C_1=311$, $C_2=\emptyset$ and $C_3=33$. By definition, $B_0=\emptyset$, $B_1=DDU$, $B_2=U$, $B_3=D$, $R_0=UU$, $R_1=\emptyset$, $R_2=DDU$ and $R_3=\emptyset$.
	Hence, 
	\[w=g(\pi)=H_rU_rU_rH_bD_bD_bU_bH_rH_bU_bH_rD_rD_rU_rH_bD_bH_r \in\M_{17},\] 
	which is given in Figure \ref{fig:red-blue-path}.
\end{example}

Consider two words $p,q$ that occur in the definitions of both $\phi$ and $g$.  Let the positions of $2$ be
$a_1<\cdots<a_h$ in $p$ and $b_1>\cdots>b_h$ in $q$.
Put $a_0=0$ and $a_{h+1}=n+1$. For $1\leq t\leq h$, define
\[
 \alpha_t=|\{i<a_t:p_i=1\}|
 \quad \text{and} \quad
 \beta_t=|\{j>b_t:q_j=1\}|.
\]
For $0\leq t\leq h$, define
\[
 \gamma_t=|\{i<a_{t+1}:p_i=3\}|.
\]
Set $\delta_0=0$ and for $1\leq t\leq h$ define
\[
 \delta_t=|\{j>b_t:q_j=3\}|.
\]
We say that a pair of words $(p,q)$ is \emph{compatible} if it satisfies the following two conditions: 
\begin{align}
 \alpha_t&>\beta_t  \qquad\qquad (1\leq t\leq h),\label{eq:c1}\\
 \gamma_t&\leq\delta_t \qquad\qquad (0\leq t\leq h).\label{eq:c2}
\end{align}

\begin{lem}\label{lem:rb-path}
	Suppose $p$ and $q$ are two words containing equal number of $1$, $2$ and $3$ with length $n$. 
	Then $(p,q)=\jt(w)$ for some $w \in \M_{2n-1}$
    if and only if $(p,q)$ is compatible.
\end{lem}

\pf Given a Motzkin path $w$ of odd length which is colored following the convention given in the paper and write 
\[w=(B_0HR_0)\,H\,(B_1HR_1)\,H\cdots H\,(B_hHR_h), \] 
where $B_i$ is the word on $U_b$ and $D_b$, while $R_i$ is the word on $U_r$ and $D_r$.
It is readily checked that $w \in \M_{2n-1}$ if and only if
the following two conditions hold: 
 \begin{itemize}
 	\item[1.] In each prefix $(B_0HR_0)H(B_1HR_1)H\cdots H (B_{t-1}HR_{t-1})$ with $1 \leq t \leq h$, the number of $D_r$ is at most that of $U_b$.
 	\item[2.] In each prefix $(B_0HR_0)H(B_1HR_1)H\cdots H B_t$ with $0 \leq t \leq h$, the number of $D_b$ is at most that of $U_r$.
 \end{itemize}
By the construction of $\phi$ and $g$, 
Condition $1$ equals to that $\alpha_{t}-1\geq \beta_t$ and 
Condition $2$ equals to that $\gamma_{t} \leq \delta_t$, as desired. 
\qed

\begin{lem}\label{lem:rb-rank}
	Suppose $p$ and $q$ are two words containing equal number of $1$, $2$ and $3$ with length $n$.
	Then $(p,q)=\pv^{-1}(\pi)$ for some $\pi \in \Sym_n(1234)$
    if and only if $(p,q)$ is compatible.
\end{lem}

\pf As in the preceding text, write $p=1A_0\,2A_1\,2\cdots2A_h$ and $q=1C_0\,2C_1\,2\cdots2C_h$. It can be checked  that $(p,q)=\pv^{-1}(\pi)$ with $\pi\in \Sym_n(1234)$ if and only if  $(p,q)$ satisfies
the following three conditions: 
\begin{itemize}
	\item[1.] There is no $p_i=3$ with $i<a_1$.
	\item[2.] Each entry $i$ such that $q_i=3$ in $C_{h-t}$ lies after the position $a_{t+1}$ with $1 \leq t \leq h-1$. Moreover, there is no $3$ in $C_0$.
	\item[3.] There exists an index $i$ such that
	 $q_i=1$ in the prefix $1C_0\,2C_1\cdots 2C_{h-t}$ lies before the position $a_{t}$ with $1 \leq t \leq h$.
\end{itemize}
Conditions $1$ and $2$ come from the fact that any $3$-rank entry has at least one smaller $2$-rank entry to its left. Moreover,  condition $1$ implies that $\gamma_0=0$, while
condition $2$ means that the number of $3$'s in $C_hC_{h-1}\cdots C_{h-t+1}$ is not smaller than the number of $3$'s in $A_0A_1\cdots A_{t}$ with $1 \leq t \leq h$. Equally, $\gamma_t\leq\delta_t$ for $1 \leq t \leq h$.
Condition $3$ comes from the fact that any $2$-rank entry
has at least one smaller $1$-rank entry to its left
Moreover, condition $3$ means $\beta_t \leq \alpha_t-1$ for  $1 \leq t \leq h$, as desired.
\qed

\begin{lem}\label{lem:rb-bijection}
For every $n\geq1$, the map $\phi$ is a bijection 
from $\M_{2n-1}$ to $\Sym_n(1234)$ with the inverse map $g$. Moreover, if $\pi=\phi(w)$, then $P_i(\pi)=J_i(w)$ and $V_i(\pi)=T_i(w)$ for each $i=1,2,3$.
\end{lem}

\pf 
By Lemmas \ref{lem:rb-path} and \ref{lem:rb-rank}, we see that the maps $\phi$ and $g$ are well-defined. The compatible pair  $(p,q)$ is the bridge of the two maps $\phi$ and $g$.
From above, we have already given the inverses of  $\jt$ and $\pv$. It follows that $\phi$ and $g$ are inverses to each other. Since $P_i(\pi),V_i(\pi),J_i(w)$ and $T_i(w)$ are totally determined by $(p,q)$, then $P_i(\pi)=J_i(w)$ and $V_i(\pi)=T_i(w)$ for each $i=1,2,3$ follows.
\qed

For a word $s=s_1\cdots s_n$ with repeated letters, define its weak
descent set by
\[
 \WDes(s)=\{i\in[n-1]:s_i\geq s_{i+1}\}.
\]

\begin{lem}\label{lem:rb-desc}
If $\pi=\phi(w)$ and $(p,q)=\jt(w)$, then we have 
\begin{equation*} 
 \Des(\pi)=\WDes(p),
 \quad
 \Des(\pi^{-1})=\WDes(q).
\end{equation*}
\end{lem}

\pf Write $r=r_1 r_2 \cdots r_n$ with $r_i=\rr_i(\pi)$, and call it the rank word of $\pi$.  One has $\pi_i<\pi_{i+1}$ if and only if $r_i<r_{i+1}$.
 The forward direction follows directly from the definition of rank. If $\pi_i>\pi_{i+1}$, we can change the longest
 increasing subsequence ending at $\pi_{i+1}$ into that ending at $\pi_i$ by replacing $\pi_{i+1}$ with $\pi_i$.
This shows $r_i\geq r_{i+1}$. Hence
$i\in\Des(\pi)$ exactly when $p_i\geq p_{i+1}$.

Since the permutation diagram of $\pi$ and $\pi^{-1}$ are reflections of each other across $y=x$, the positions and values of the entries in each increasing subsequence are interchanged, while its length remains unchanged.
This implies that $q$ is the rank word of $\pi^{-1}$. Applying the same argument above to $q$ proves the
second equality. \qed

\noindent { \textbf{Proof of Theorem \ref{thm:bijodd}}}

\pf
Following from Lemma \ref{lem:rb-bijection}, it suffices to show that $u_0(w)=\idr(\pi)$ and $d_0(w)=\idr(\pi^{-1})$.
Notice that the last block of $w$ is red. By the  red-blue rule, the first blue block of $w$ contains no down step, and its last red block contains no up
step. Therefore the initial runs of $1$'s in $p$ and $q$ have lengths $u_0(w)+1$ and $d_0(w)+1$, respectively. Together with Lemma \ref{lem:rb-desc}, this yields
$u_0(w)=\idr(\pi)$ and $d_0(w)=\idr(\pi^{-1})$. This completes the proof.
\qed

\noindent { \textbf{Proof of Theorem \ref{thm:bijeven}}}

\pf
For $n\geq1$, appending a final horizontal step gives a bijection from $\M_{2n}$ to the set $\{w\in\M_{2n+1}:w\text{ ends with }H\}$.
By Lemma \ref{lem:rb-desc} and setting $d_0(w)=\idr(\pi^{-1})=0$, the map $\phi$ restricts a bijection from $\{w\in\M_{2n+1}:w\text{ ends with }H\}$ to 
$\{\pi\in\Sym_{n+1}(1234):
\pi^{-1}(1)<\pi^{-1}(2)\}$, as desired.
\qed

\section{A proof of Theorem \ref{thm:bijodd-alternative}}
\label{sec:bijodd-alternative}

For $\sigma \in \Sym_a$ and $\tau \in \Sym_b$, the \emph{skew sum} of $\sigma$ and $\tau$, denoted   $\sigma\ominus\tau$, is defined by adding $b$ to each
entry of $\sigma$ and then appending $\tau$.
Define the cut set of $\pi$ by
\begin{equation*}
	\Cut(\pi)=\{c:1\le c<n,\ \{\pi_1,\ldots,\pi_c\}
	=\{n-c+1,\ldots,n\}\}.
\end{equation*}
A permutation is \emph{skew-indecomposable} if its cut set is empty. Define 
\[\SI_n=\{\sigma \in \Sym_n(1234): \sigma \text{ is skew-indecomposable}\}\]
and 
\[
\RM_n=\{ w \in \M_n: w \,\,\text{has no} \, H_b \,\,\text{at x-axis}\}.
\]

\begin{thm}\label{thm:odd-indecomposable}
	For $n \geq 1$, there is a bijection from  $\SI_n$ to $\RM_{2n-1}$.
\end{thm}

Given $\sigma \in \Sym_n$,  $(i,\sigma_i)$ is said to be a \emph{middle point} if $\sigma_i$ is the middle entry of some $123$-pattern. Set  
\begin{equation*}
	\E(\sigma)=\{(i,\sigma_i):\text{there exist }a<i<b
	\text{ with }\sigma_a<\sigma_i<\sigma_b\}.
\end{equation*}
If $\sigma \in \Sym_n(1234)$ with $n \geq 2$ and $(i,\sigma_i) \in \E(\sigma)$, then
$\rr_i(\sigma)=2$. 
Let the right-to-left maxima of $\sigma$ from left to right be $(t_1,u_1),\ldots,(t_s,u_s)$ with $t_1<\cdots<t_s=n$
and $n=u_1>\cdots>u_s$.
For each $1 \leq j <s$, define a \emph{candidate corner} $c_j=(t_j+1,u_{j+1}+1)$.
Clearly, each $c_j$ lies in $[n]\times[n]$. The corner $c_j=(x,y)$ is \emph{retained} if it satisfies that 
$\E(\sigma)\cap ([1,x]\times[1,y])=\varnothing$.
Let $\B(\sigma)$ be the set of retained corners of $\sigma$ 
and
\[\Z(\sigma)=\E(\sigma)\cup\B(\sigma).\]
Define two words $p=p_1\cdots p_n$ and $q=q_1\cdots q_n$ by
letting
\begin{equation}\label{eq:word-sets}
 \begin{aligned}
    \{i:p_i=2\}&=\{x: (x,y) \in \Z(\sigma)\},
    \qquad\{i:p_i=3\}=P_3(\sigma),\\
    \{v:q_v=2\}&=\{y: (x,y) \in \Z(\sigma)\},
    \qquad \{v:q_v=3\}=V_3(\sigma).
 \end{aligned}
\end{equation}
Other entries of $p$ and q are set to be $1$.
Finally, define $w=f(\sigma)=\jt^{-1}(p,q)$.
For $n=1$, let $w=f(1)=H_r$.

\begin{example}
	Let $\sigma=479861523$. Then, $P_3(\sigma)=\{3,4,9\}$
	and $V_3(\sigma)=\{3,8,9\}$. Its middle points and retained corners are
	$\E(\sigma)=\{(2,7),(8,2)\}$ and 
	$\B(\sigma)=\{(6,6)\}.$ The procedure of finding $\B(\sigma)$ is shown in Figure \ref{fig:479861523}.
	
\begin{figure}[htbp]
	\centering
	 \scalebox{0.8}{
	\begin{tikzpicture}[
		x=0.62cm,
		y=0.62cm,
		font=\small,
		every node/.style={anchor=center},
		gridline/.style={
			draw=gray!28,
			line width=0.35pt
		},
		middle/.style={
			circle,
			fill=orange!85!black,
			minimum size=7pt,
			inner sep=0pt,
			outer sep=0pt
		},
		rtl/.style={
			circle,
			draw=blue!65!black,
			line width=0.9pt,
			minimum size=10pt,
			inner sep=0pt,
			outer sep=0pt
		},
		rankthree/.style={
			rectangle,
			fill=blue!65!black,
			minimum size=4pt,
			inner sep=0pt,
			outer sep=0pt
		},
		ranktwo/.style={
			circle,
			fill=black,
			minimum size=4pt,
			inner sep=0pt,
			outer sep=0pt
		},
		retained/.style={
			diamond,
			fill=green!60!black,
			minimum size=8pt,
			inner sep=0pt,
			outer sep=0pt
		},
		rejected/.style={
			cross out,
			draw=red!80!black,
			line width=0.85pt,
			minimum size=6pt,
			inner sep=0pt,
			outer sep=0pt
		},
		legend text/.style={
			anchor=west,
			font=\footnotesize,
			align=left
		}
		]
		
		\fill[gray!12] (0,0) rectangle (6,6);
		
		\draw[gridline,step=1] (0,0) grid (9,9);
		
		\draw[->,line width=0.7pt] (0,0) -- (9.45,0);
		\draw[->,line width=0.7pt] (0,0) -- (0,9.45);
		
		\foreach \i in {1,...,9} {
			\node[below=2pt] at (\i,0) {\i};
			\node[left=2pt] at (0,\i) {\i};
		}
		
		\node[rejected] at (4,9) {};
		\node[rejected] at (5,7) {};
		\node[rejected] at (8,4) {};
		
		\fill (1,4) circle (2pt);
		\fill (6,1) circle (2pt);
		
		\node[middle] at (2,7) {};
		\node[middle] at (8,2) {};
		
		\node[rtl] at (3,9) {};
		\node[rankthree] at (3,9) {};
		
		\node[rtl] at (4,8) {};
		\node[rankthree] at (4,8) {};
		
		\node[rtl] at (9,3) {};
		\node[rankthree] at (9,3) {};
		
		\node[rtl] at (5,6) {};
		\node[ranktwo] at (5,6) {};
		
		\node[rtl] at (7,5) {};
		\node[ranktwo] at (7,5) {};
		
		\node[retained] at (6,6) {};
		
		\begin{scope}[shift={(11,7.4)}]
			
			\node[middle] at (0,0) {};
			\node[legend text] at (0.28,0)
			{middle point};

			\node[rtl] at (0,-1) {};
			\node[rankthree] at (0,-1) {};
			\node[legend text] at (0.28,-1)
			{right-to-left maximum \\with rank $3$};
			
			\node[rtl] at (0,-2.3) {};
			\node[ranktwo] at (0,-2.3) {};
			\node[legend text] at (0.28,-2.3)
			{right-to-left maximum \\with rank $2$};
			
			\node[rejected] at (0,-3.4) {};
			\node[legend text] at (0.28,-3.4)
			{rejected corner};
			
			\node[retained] at (0,-4.3) {};
			\node[legend text] at (0.28,-4.3)
			{retained corner};
			
		\end{scope}
	\end{tikzpicture}
		}
	\caption{The choice of $\B(\sigma)$ with $\sigma=479861523$.}
	\label{fig:479861523}
\end{figure}
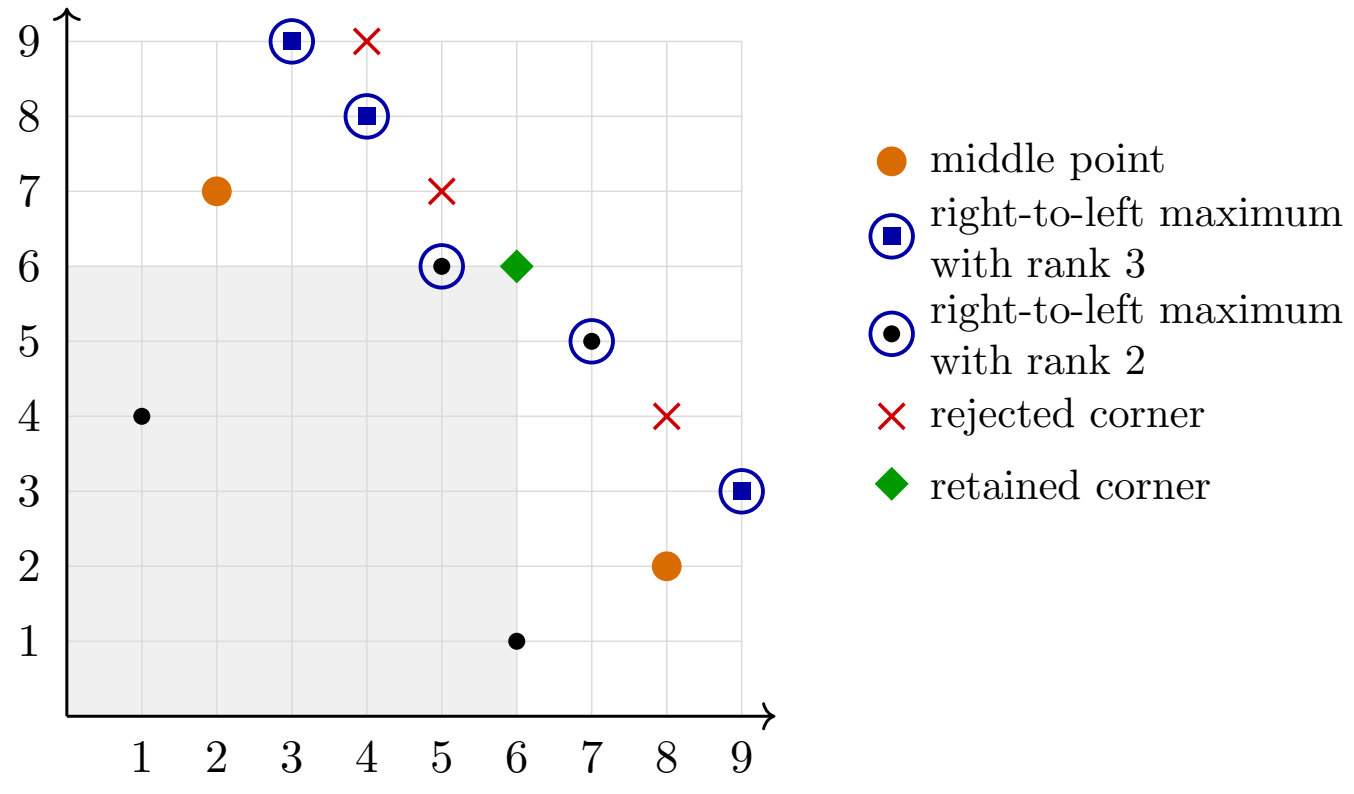
Thus, $\Z(\sigma)=\{(2,7),(6,6),(8,2)\}.$
Consequently,
	\begin{align*}
		\{i:p_i=2\}&=\{2,6,8\},
		\qquad\{i:p_i=3\}=\{3,4,9\},\\
		\{v:q_v=2\}&=\{2,6,7\},
		\qquad \{v:q_v=3\}=\{3,8,9\}.
	\end{align*}
Therefore the two encoding words of $f(\sigma)$ are
$p=123312123$ and $q=123112233$, and hence
	$f(479861523)$ is the path given in Figure \ref{fig:red-blue-path}.
\end{example}

To show the map $f$ is well-defined, we need the following lemma.

\begin{lem}\label{lem:per-pq-compatible}
	 The pair of words $(p,q)$ defined in (\ref{eq:word-sets}) is compatible.
\end{lem}

\pf First, we wish to show that 
$\{x: (x,y) \in \Z(\sigma)\} \cap P_3(\sigma)= \emptyset$.
Since $\{x: (x,y) \in \E(\sigma)\} \cap P_3(\sigma)= \emptyset$, it suffices to show  $\{x: (x,y) \in \B(\sigma)\} \cap P_3(\sigma)= \emptyset.$ 
Set $\Rlmaxp(\sigma)=\{t_1,\ldots,t_s\}$ and $\Rlmaxv(\sigma)=\{u_s,\ldots,u_1\}$.
Given $(x,y) \in \B(\sigma)$, let $(x,y)=(t_j+1,u_{j+1}+1)$
for some $1 \leq j <s$. If $x=t_{j+1}$, then $\E(\sigma) \cap [1,t_{j+1}] \times [1,u_{j+1}]=\emptyset$ by the fact that $\E(\sigma) \cap [1,x] \times [1,y]=\emptyset$. It implies that $\rr_{t_{j+1}}(\sigma)=2$, and hence $x \notin P_3(\sigma)$. If $x \neq t_{j+1}$, by Proposition \ref{prop:rankeq} we have $x \notin P_3(\sigma)$. Conversely, for $x \in P_3(\sigma)$ we may similarly show that $x$ is not the first coordinate of a point in $\B(\sigma)$. Following from the same reason, we have 
$\{y: (x,y) \in \Z(\sigma)\} \cap V_3(\sigma)= \emptyset$.

Next, we show that $(p,q)$ satisfies (\ref{eq:c1}) and (\ref{eq:c2}). Write the points of $\mathcal Z(\sigma)$ as
$(x_1,y_1),\ldots,(x_h,y_h)$ with $x_1<x_2 <\cdots <x_h$.
The construction of $\mathcal Z(\sigma)$ gives $y_1>\cdots>y_h$. If $h=0$, they certainly hold.
Fix $1\le t\le h$, and set
	\[
	A_t=\{i<x_t: i\notin P_3(\sigma)\},
	\qquad
	D_t=\{v>y_t:v\notin V_3(\sigma)\}.
	\]
The number of positions (resp. values)  labeled $2$ before $x_t$ (resp. greater than $y_t$) is $t-1$. Therefore,
$\alpha_t=|A_t|-(t-1)$ and $\beta_t=|D_t|-(t-1).$
We now show that $|A_t|>|D_t|$ through the following  two cases.
	
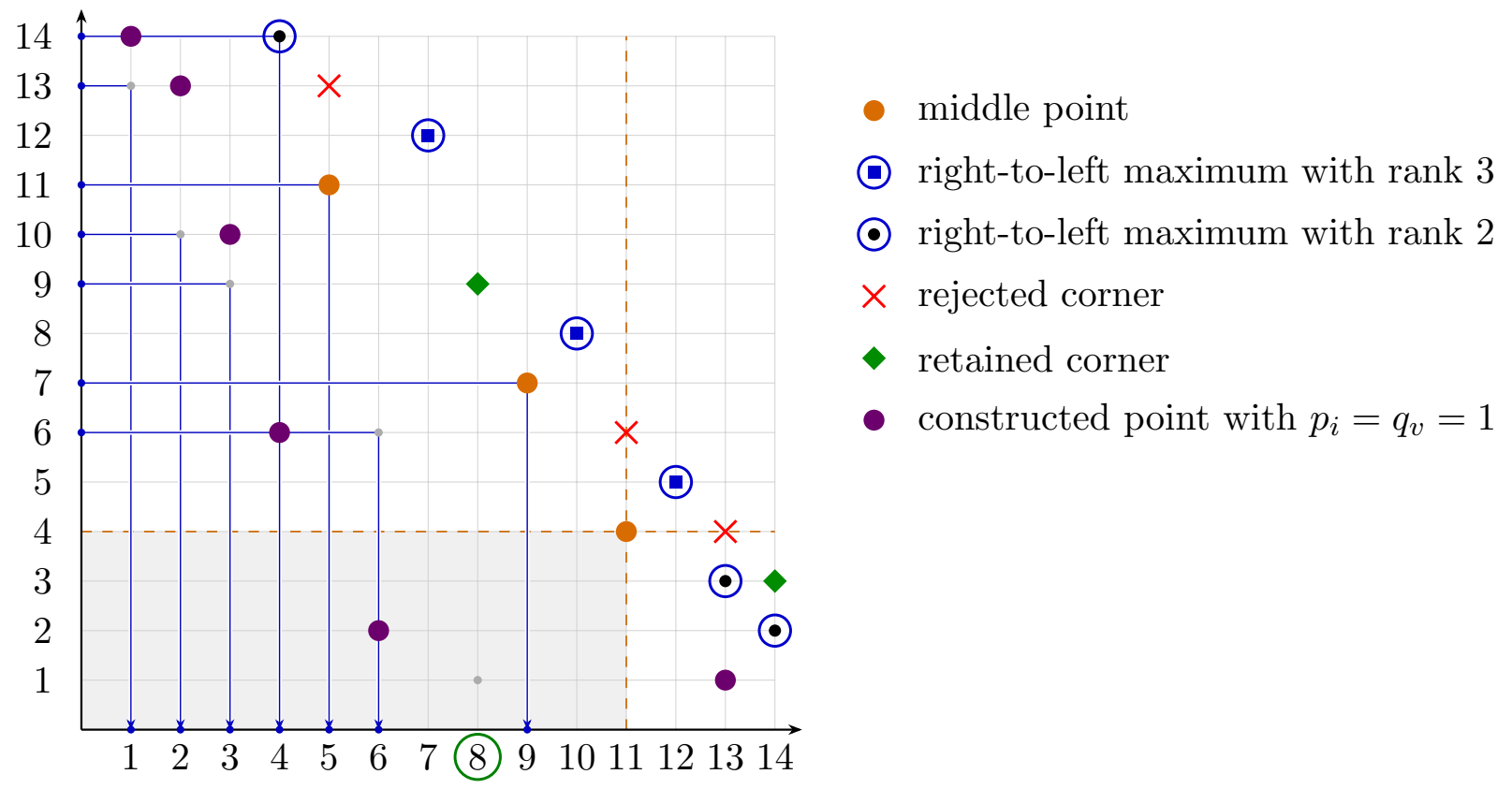
\begin{figure}[ht]
	\centering
	\scalebox{0.7}{
	\begin{tikzpicture}[x=0.50cm,y=0.50cm,
		middle/.style={circle,fill=orange!85!black,draw=none,inner sep=0pt,minimum size=6pt},
		outer/.style={circle,draw=blue!80!black,line width=.8pt,inner sep=0pt,minimum size=9pt},
		three/.style={rectangle,fill=blue!80!black,inner sep=0pt,minimum size=3.6pt},
		two/.style={circle,fill=black,inner sep=0pt,minimum size=3.5pt},
		kept/.style={diamond,fill=green!55!black,inner sep=0pt,minimum size=7pt},
		rejected/.style={cross out,draw=red,line width=.8pt,inner sep=0pt,minimum size=5.5pt},
		constructed/.style={circle,fill=violet!85!black,inner sep=0pt,minimum size=6pt},
		match/.style={blue!75!black,line width=.38pt,-{Stealth[length=3pt,width=2.5pt]},
			preaction={draw=white,line width=1.3pt}},
		every node/.style={font=\small}]
		
		\fill[gray!12] (0,0) rectangle (11,4);
		\draw[step=1,gray!35,line width=.25pt] (0,0) grid (14,14);
		\draw[orange!80!black,dashed,line width=.45pt] (11,0)--(11,14);
		\draw[orange!80!black,dashed,line width=.45pt] (0,4)--(14,4);
		\draw[-{Stealth[length=4pt]},line width=.6pt] (0,0)--(14.55,0);
		\draw[-{Stealth[length=4pt]},line width=.6pt] (0,0)--(0,14.55);
		
		\foreach \v/\i in {14/4,13/1,11/5,10/2,9/3,7/9,6/6}{
			\draw[match] (0,\v)--(\i,\v)--(\i,0);
			\fill[blue!75!black] (0,\v) circle (1.1pt);
			\fill[blue!75!black] (\i,0) circle (1.1pt);
		}
		
		\foreach \i/\v in {1/13,2/10,3/9,4/14,5/11,6/6,7/12,
			8/1,9/7,10/8,11/4,12/5,13/3,14/2}
		\fill[gray!65] (\i,\v) circle (1.2pt);
		
		\foreach \i/\v in {1/14,2/13,3/10,4/6,6/2,13/1}
		\node[constructed] at (\i,\v) {};
		\foreach \i/\v in {5/11,9/7,11/4}
		\node[middle] at (\i,\v) {};
		\foreach \i/\v in {7/12,10/8,12/5}{
			\node[outer] at (\i,\v) {};
			\node[three] at (\i,\v) {};
		}
		\foreach \i/\v in {4/14,13/3,14/2}{
			\node[outer] at (\i,\v) {};
			\node[two] at (\i,\v) {};
		}
		\foreach \i/\v in {5/13,11/6,13/4}
		\node[rejected] at (\i,\v) {};
		\foreach \i/\v in {8/9,14/3}
		\node[kept] at (\i,\v) {};
		
		\foreach \k in {1,...,14}{
			\node at (\k,-.55) {\k};
			\node[anchor=east] at (-.3,\k) {\k};
		}
		
		\node[circle,draw=green!50!black,line width=.8pt,minimum size=13pt,
		inner sep=0pt] at (8,-.55) {};
		
		\node[middle] at (16,12.5) {};
		\node[anchor=west] at (16.6,12.5) {middle point};
		\node[outer] at (16,11.25) {};
		\node[three] at (16,11.25) {};
		\node[anchor=west] at (16.6,11.25) {right-to-left maximum with rank $3$};
		\node[outer] at (16,10) {};
		\node[two] at (16,10) {};
		\node[anchor=west] at (16.6,10) {right-to-left maximum with rank $2$};
		\node[rejected] at (16,8.75) {};
		\node[anchor=west] at (16.6,8.75) {rejected corner};
		\node[kept] at (16,7.5) {};
		\node[anchor=west] at (16.6,7.5) {retained corner};
		\node[constructed] at (16,6.25) {};
		\node[anchor=west] at (16.6,6.25)
		{constructed point with $p_i=q_v=1$};
		
	\end{tikzpicture}
}
	\caption{The injection from
		$D_4$ to $A_4$
		for $\sigma=(13,10,9,14,11,6,12,1,7,8,4,5,3,2)$
		with $(x_4,y_4)=(11,4)$ and $8$ is an extra position.}
	\label{fig:injection-case1}
\end{figure}

\noindent\textbf{Case 1: $(x_t,y_t)\in\mathcal E(\sigma)$.}
	Here $\sigma_{x_t}=y_t$. Suppose $a<x_t$ with
	$w=\sigma_a=\min\{\sigma_1,\ldots,\sigma_{x_t-1}\}.$
	Since $(x_t,y_t)$ is a middle point, we have $w<y_t$.
	Moreover, $\sigma_a$ is a $1$-rank entry in $\sigma$, and hence $a\in A_t$ and $\sigma_a \notin D_t$.
	For each $v\in D_t$, set $b=\sigma^{-1}(v)$.
	If $b>x_t$, then $a<x_t<b$ and $\sigma_a<\sigma_{x_t}<\sigma_b$. Combining the fact that $\sigma$ avoids $1234$, we have $\rr_b(\sigma)=3$. Hence, $v\in V_3(\sigma)$, a contradiction.
	Since $b\ne x_t$, we have $b<x_t$ and $b\in A_t$.
	It follows that the map $v\mapsto\sigma^{-1}(v)$ is an injection from $D_t$ into $A_t$.
	Since  $a\in A_t$ and $\sigma_a \notin D_t$,
    the map $v\mapsto\sigma^{-1}(v)$
	is an injection from $D_t$ into $A_t\setminus\{a\}$,
	and hence $|A_t|\geq |D_t|+1$. See Figure \ref{fig:injection-case1} as a description of this case.
	
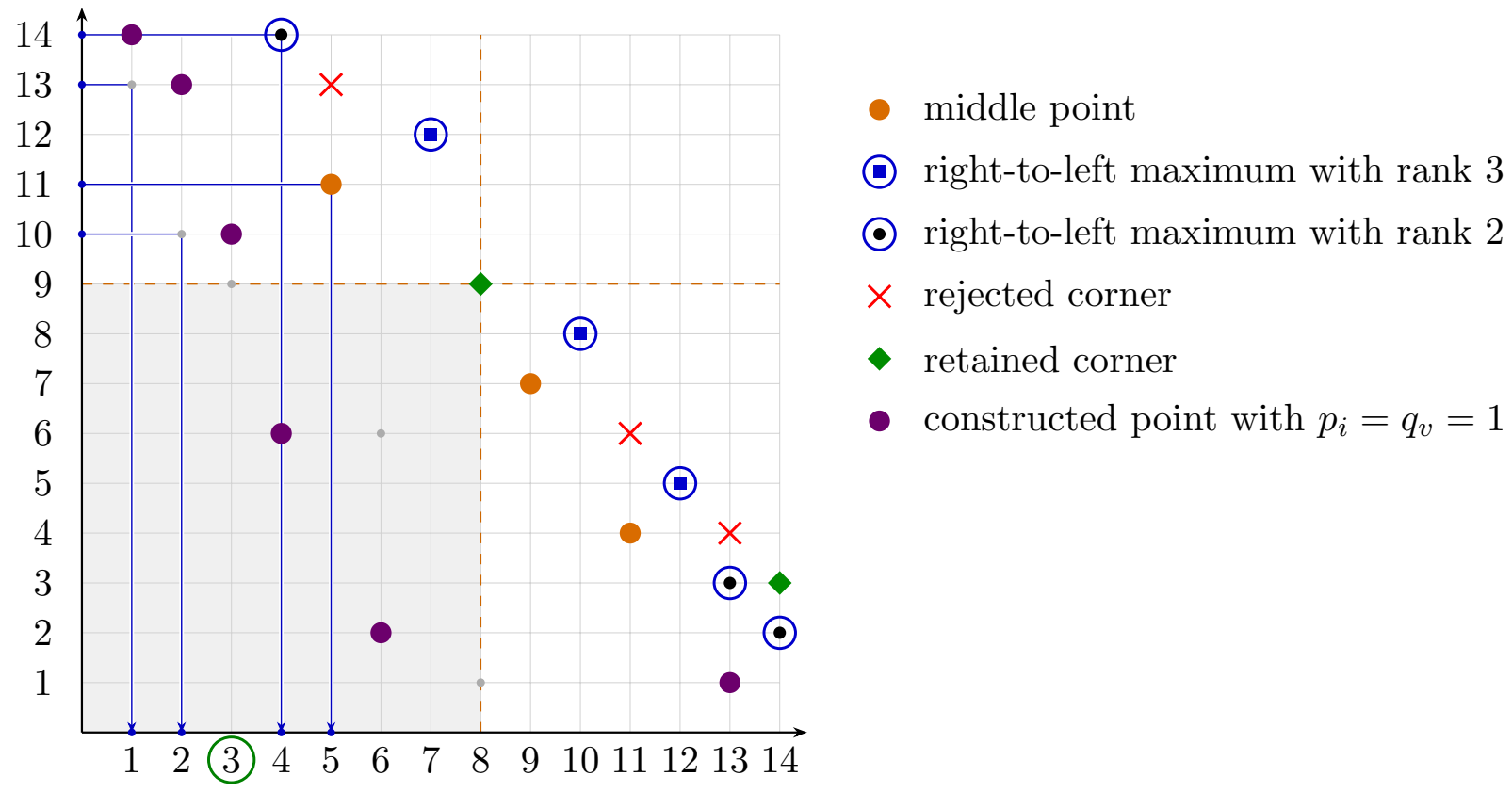
\begin{figure}[ht]
	\centering
	\scalebox{0.7}{
	\begin{tikzpicture}[x=0.50cm,y=0.50cm,
		middle/.style={circle,fill=orange!85!black,draw=none,
			inner sep=0pt,minimum size=6pt},
		outer/.style={circle,draw=blue!80!black,line width=.8pt,
			inner sep=0pt,minimum size=9pt},
		three/.style={rectangle,fill=blue!80!black,
			inner sep=0pt,minimum size=3.6pt},
		two/.style={circle,fill=black,
			inner sep=0pt,minimum size=3.5pt},
		kept/.style={diamond,fill=green!55!black,
			inner sep=0pt,minimum size=7pt},
		rejected/.style={cross out,draw=red,line width=.8pt,
			inner sep=0pt,minimum size=5.5pt},
		constructed/.style={circle,fill=violet!85!black,
			inner sep=0pt,minimum size=6pt},
		match/.style={blue!75!black,line width=.38pt,
			-{Stealth[length=3pt,width=2.5pt]},
			preaction={draw=white,line width=1.3pt}},
		every node/.style={font=\small}]
		
		\fill[gray!12] (0,0) rectangle (8,9);
		\draw[step=1,gray!35,line width=.25pt]
		(0,0) grid (14,14);
		\draw[orange!80!black,dashed,line width=.45pt]
		(8,0)--(8,14);
		\draw[orange!80!black,dashed,line width=.45pt]
		(0,9)--(14,9);
		\draw[-{Stealth[length=4pt]},line width=.6pt]
		(0,0)--(14.55,0);
		\draw[-{Stealth[length=4pt]},line width=.6pt]
		(0,0)--(0,14.55);
		
		\foreach \v/\i in {14/4,13/1,11/5,10/2}{
			\draw[match] (0,\v)--(\i,\v)--(\i,0);
			\fill[blue!75!black] (0,\v) circle (1.1pt);
			\fill[blue!75!black] (\i,0) circle (1.1pt);
		}
		
		\foreach \i/\v in {
			1/13,2/10,3/9,4/14,5/11,6/6,7/12,
			8/1,9/7,10/8,11/4,12/5,13/3,14/2}
		\fill[gray!65] (\i,\v) circle (1.2pt);
		
		\foreach \i/\v in {1/14,2/13,3/10,4/6,6/2,13/1}
		\node[constructed] at (\i,\v) {};
		
		\foreach \i/\v in {5/11,9/7,11/4}
		\node[middle] at (\i,\v) {};
		
		\foreach \i/\v in {7/12,10/8,12/5}{
			\node[outer] at (\i,\v) {};
			\node[three] at (\i,\v) {};
		}
		
		\foreach \i/\v in {4/14,13/3,14/2}{
			\node[outer] at (\i,\v) {};
			\node[two] at (\i,\v) {};
		}
		
		\foreach \i/\v in {5/13,11/6,13/4}
		\node[rejected] at (\i,\v) {};
		
		\foreach \i/\v in {8/9,14/3}
		\node[kept] at (\i,\v) {};
		
		\foreach \k in {1,...,14}{
			\node at (\k,-.55) {\k};
			\node[anchor=east] at (-.3,\k) {\k};
		}
		
		\node[circle,draw=green!50!black,line width=.8pt,
		minimum size=13pt,inner sep=0pt] at (3,-.55) {};
		
		\node[middle] at (16,12.5) {};
		\node[anchor=west] at (16.6,12.5)
		{middle point};
		
		\node[outer] at (16,11.25) {};
		\node[three] at (16,11.25) {};
		\node[anchor=west] at (16.6,11.25)
		{right-to-left maximum with rank $3$};
		
		\node[outer] at (16,10) {};
		\node[two] at (16,10) {};
		\node[anchor=west] at (16.6,10)
		{right-to-left maximum with rank $2$};
		
		\node[rejected] at (16,8.75) {};
		\node[anchor=west] at (16.6,8.75)
		{rejected corner};
		
		\node[kept] at (16,7.5) {};
		\node[anchor=west] at (16.6,7.5)
		{retained corner};
		
		\node[constructed] at (16,6.25) {};
		\node[anchor=west] at (16.6,6.25)
		{constructed point with $p_i=q_v=1$};
		
	\end{tikzpicture}
}
	\caption{The injection from
		$D_2$ to $A_2$
		for $\sigma=(13,10,9,14,11,6,12,1,7,8,4,5,3,2)$
		with $(x_2,y_2)=(8,9)$ and $3$ is an extra position.}
	\label{fig:injection-case2}
\end{figure}

\noindent \textbf{Case 2: $(x_t,y_t)\in\mathcal B(\sigma)$.}
	For some $1 \leq j <s$, we have
	$x_t=t_j+1$ and $y_t=u_{j+1}+1.$ Every value not smaller than $y_t$ occurs before $x_t$.
	Again, $v\mapsto\sigma^{-1}(v)$ defines an injection
	from $D_t$ into $A_t$. 
	We claim that $y_t\notin V_3(\sigma)$.
	Otherwise, consider $b=\sigma^{-1}(y_t)<x_t$. There would exist	$a<c<b$ such that $\sigma_a<\sigma_c<\sigma_b=y_t.$
	Then $(c,\sigma_c)$ would be a middle point in
	$[1,x_t]\times[1,y_t]$, contrary to the fact that 
	$(x_t,y_t)$ is a retained corner.
	Therefore $\sigma^{-1}(y_t)\in A_t$.
	Notice that $y_t\notin D_t$, we have $|A_t|\ge |D_t|+1$. See Figure \ref{fig:injection-case2} as a description of this case. This completes the proof of relation (\ref{eq:c1}).

We now show that $\gamma_t\leq\delta_t$ with $0\leq t\leq h$. For $t=0$, we claim that $\gamma_0=0$. Suppose, to the contrary,
 that there exists an entry $\sigma_i$ with $i<x_1$ and $\rr_i(\sigma)=3$. Let $a<b<i$ such that
 $\sigma_a<\sigma_b<\sigma_i.$
 Then $(b,\sigma_b)$ is a middle point, so
 $(b,\sigma_b)\in\mathcal E(\sigma)\subseteq\mathcal Z(\sigma)$.
 Since $\sigma_b$ has rank exactly $2$, we have $b \geq x_1$. This is impossible through $i<x_1$. 
 Hence $\gamma_0=0=\delta_0$.

 For $1\leq t\leq h$, write $C_t=\{i: i<x_{t+1}, i \in P_3(\sigma)\}$ and $E_t=\{v: v>y_{t}, v \in V_3(\sigma)\}$ with $x_{h+1}=n+1$.
 It suffices to construct an injection from $C_t$ to $E_t$.
For each $i<x_{t+1}$ with $i\in P_3(\sigma)$, since 
$\{x: (x,y) \in \Z(\sigma)\} \cap P_3(\sigma)= \emptyset$ we have  $i\neq x_t$. We wish to show that $\sigma_i>y_t$.
 Assume to the contrary that $\sigma_i\leq y_t$.
 Since $\rr_i(\sigma)=3$, there exist $a<b<i$ such that
 $\sigma_a<\sigma_b<\sigma_i$, and hence
 $(b,\sigma_b)\in\E(\sigma)$. In particular,
 $(b,\sigma_b)\in\mathcal Z(\sigma)$,
 $b<i<x_{t+1}$ and $\sigma_b<\sigma_i\leq y_t.$
 If $b<x_t$, then $\sigma_b>y_t$, a contradiction.
 If $b=x_t$ and $(x_t,y_t)\in\E(\sigma)$, then $\sigma_b=y_t$, a contradiction.
  For $b=x_t$ and  $(x_t,y_t)\in\B(\sigma)$, this is impossible since the retained-corner condition excludes a middle point in $[1,x_t]\times[1,y_t]$. 
 Finally, if $x_t<b<x_{t+1}$, then
 $(b,\sigma_b)\in\Z(\sigma)$ lies strictly between $(x_t,y_t)$ and $(x_{t+1},y_{t+1})$, again a contradiction.
 Thus, we have $\sigma_i>y_t$. 
 
 Now for $i \in C_t$, we have $\sigma_i \in E_t$,
 so $i\mapsto\sigma_i$ is an injection from $C_t$ to $E_t$. Therefore, $|C_t| \leq |E_t|$, as desired. This completes the proof.
\qed

\begin{lem}\label{lem:maptoRM}
	For $\sigma \in \SI_{n}$, we have $f(\sigma) \in \RM_{2n-1}$.
\end{lem}
\pf	The case $n=1$ is immediate. Suppose $n \geq 2$.
Let $w=f(\sigma)$ and write $(x_1,y_1),\ldots,(x_h,y_h)$ 
as the points in $\Z(\sigma)$ with 
$x_1<\cdots<x_h$. The $t$-th $H_b$-step of $w$ corresponds to $(x_t,y_t)$. Set $(p,q)=\jt(w)$ and define 
	\[
	X_t=\alpha_t-\beta_t-1,
	\qquad
	Y_t=\delta_t-\gamma_{t-1}.
	\]	
Since $(p,q)$ is compatible, we have $X_t \geq 0$ and $Y_t\geq 0$. In view of the fact that the height of $H_b$ is $X_t+Y_t$, it suffices to show that $X_t$ and $Y_t$ cannot both be zero. Assume to the contrary that $X_t=Y_t=0$. We claim $c=x_t-1=n-y_t+1$ is a cut of $\sigma$. That is, 
$\{\sigma_1,\ldots,\sigma_c\}=\{y_t,y_t+1,\ldots,n\}$.
The condition $X_t=Y_t=0$ and the definitions
of $\alpha_t,\beta_t,\gamma_{t-1},\delta_t$ give
	\[
	x_t-1=\alpha_t+(t-1)+\gamma_{t-1}	=\beta_t+1+(t-1)+\delta_t
	=n-y_t+1.
	\]

If $(x_t,y_t)\in\E(\sigma)$, we have $\sigma_{x_t}=y_t$ and there exist $a<x_t<b$ with $\sigma_a<y_t<\sigma_b$.
Thus $\rr_b(\sigma)=3$. By Proposition \ref{prop:rankeq},  every $3$-rank entry lying before $x_t$ is greater than $\sigma_b$, and hence greater than $y_t$.
Together with $\sigma_b$, these give at least
$\gamma_{t-1}+1$ rank-$3$ values greater than $y_t$.
Therefore, $\delta_t\geq\gamma_{t-1}+1$,  a contradiction with the fact $Y_t=0$.
This case is impossible.

If $(x_t,y_t)\in\mathcal B(\sigma)$,
let $(x_t,y_t)=(t_j+1,u_{j+1}+1),$
where $(t_j,u_j)$ and $(t_{j+1},u_{j+1})$
are consecutive right-to-left maxima of $\sigma$.
The largest entry to the right of position $t_j$
is $u_{j+1}$.
Hence every value at least $y_t=u_{j+1}+1$
occurs before $x_t$, so $\{y_t,y_t+1,\ldots,n\}
\subseteq \{\sigma_1,\ldots,\sigma_{x_t-1}\}.$
Since both sets have cardinality
$n-y_t+1=x_t-1$, they are equal.
Since $1\leq c=x_t-1<n$, so $c$ is a cut of $\sigma$.
Since $\sigma$ is skew-indecomposable,
this can not occur. Thus every $H_b$-step of $w$ has positive height, and we have $w\in\mathcal{RM}_{2n-1}$.
\qed

Given $w \in \RM_{2n-1}$ with $n \geq 2$, let $(p,q)=\jt(w)$.
For $k\in\{2,3\}$,   write
\[J_k(w)=\{i:p_i=k\}
=\{i_1^{(k)}<\cdots<i_{s_k}^{(k)}\},
\quad 
T_k(w)=\{j:q_j=k\}
=\{v_1^{(k)}>\cdots>v_{s_k}^{(k)}\}.
\]
Set
\[
\R_k(w)=\{(i^{(k)}_1,v^{(k)}_1),\ldots,(i^{(k)}_{s_k},v^{(k)}_{s_k})\}.
\]
Split $\R_2(w)$ by
\begin{align*}
	\E_w
	&=\{(x,y)\in\R_2(w):
	\text{there exists }(r,s)\in\R_3(w)
	\text{ with }r>x\text{ and }s>y\},\\
	\B_w&=\R_2(w)\setminus\mathcal E_w.
\end{align*}

\begin{lem}\label{lem:BEZ}
	Let $\sigma\in\SI_n$ and let $w=f(\sigma)$. 
	Then, we have
	\[
	\E_w=\E(\sigma),
	\quad
	\B_w=\B(\sigma),
	\quad
	\R_2(w)=\Z(\sigma).
	\]
\end{lem}
\pf
	By the definition of $f$, clearly we have   $\R_2(w)=\Z(\sigma).$
	Given $(x,\sigma_x)\in\E(\sigma)$, there exist
	$a<x<b$ such that $	\sigma_a<\sigma_x<\sigma_b.$
	Since $\sigma_x$ has rank $2$, the entry $\sigma_b$ has rank $3$. Hence $(b,\sigma_b) \in \R_3(w)$, $x <b$ and $\sigma_x <\sigma_b$. Therefore, $(x,\sigma_x)\in\E_w.$
	Conversely, suppose that $(x,y)\in\E_w$. Then
	$(x,y)\in\mathcal Z(\sigma)$ and there exists $(r,s) \in \R_3(w)$ with $r>x$ and $s>y$. The point $(x,y)$ cannot be a retained corner. Assume to the contrary that $(x,y)=(t_j+1,u_{j+1}+1)$ is a candidate corner. Then each entry to the right of position $t_j$ is at most $u_{j+1}$, which is less than $y$, against the existence of $(r,s)$ above.
	 Thus, $(x,y)\in\E(\sigma)$.  Combining the fact that
	 $\R_2(w)=\Z(\sigma)$, we have $\B_w=\B(\sigma)$.
\qed

Define a set of points $\C_w$ according to $w$ geometrically as follows:
\begin{itemize}
	\item[1.] On the $x$-axis (resp. $y$-axis), mark $n$ and $x-1$ (resp. $y-1$) with each $(x,y)\in\R_2(w)$. 
	\item[2.] Retain $\R_3(w)$ and move each  $(x,y)\in\R_2(w)$ to $\left(x-\frac12,y-\frac12\right).$
	\item[3.] Shade the closed northeast and southwest regions of each retained or shifted point in step $2$.
	\item[4.] The elements of $\C_w$ are defined to be the unshaded points with its $x$-coordinate and $y$-coordinate both marked in step $1$.
\end{itemize}

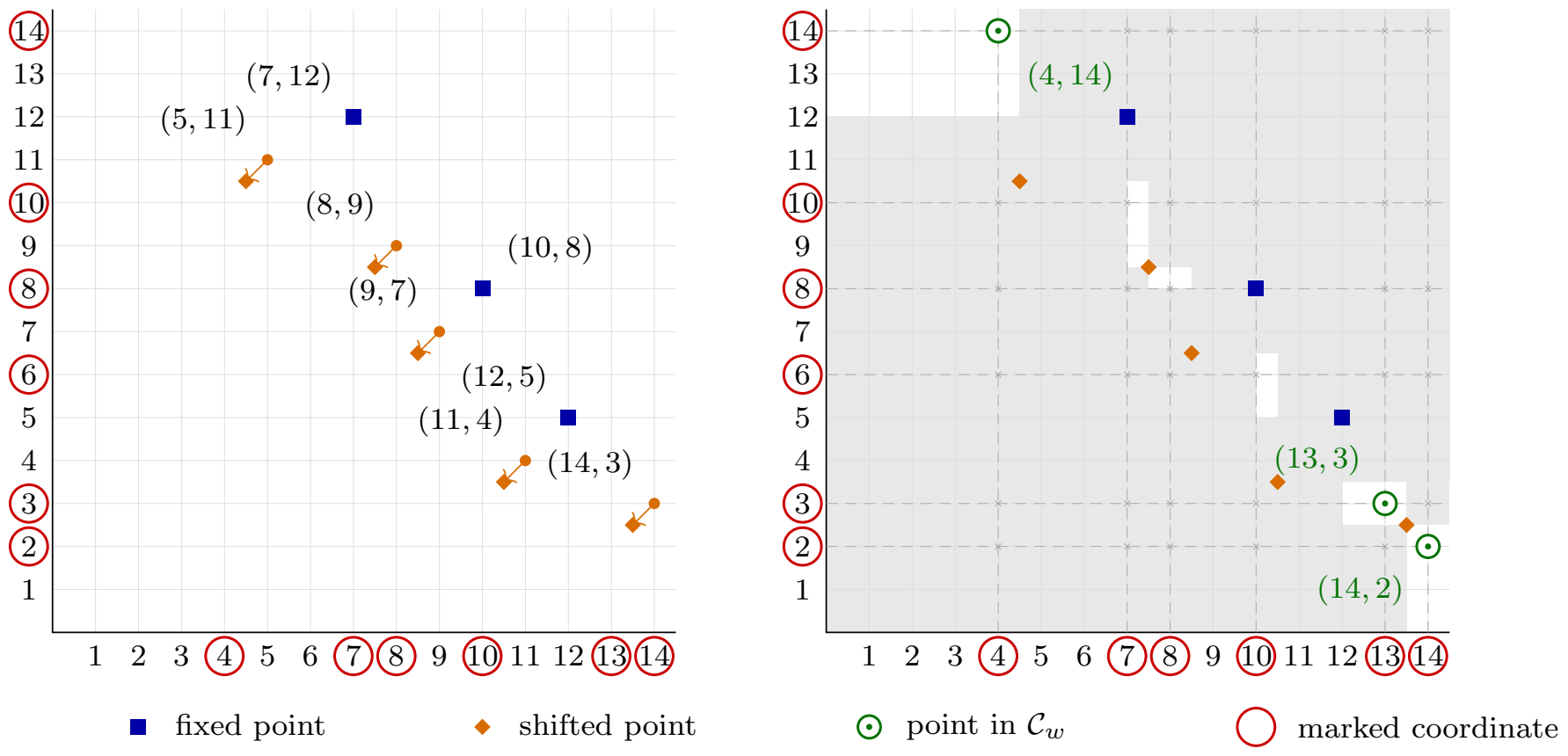
\begin{figure}[htbp]
	\centering
	\scalebox{0.75}{
		\resizebox{\linewidth}{!}{%
			\begin{tikzpicture}[
				x=.4cm,y=.4cm,
				every node/.style={font=\scriptsize},
				fixed/.style={
					rectangle,
					fill=blue!65!black,
					inner sep=0pt,
					minimum size=4pt
				},
				shifted/.style={
					diamond,
					fill=orange!85!black,
					inner sep=0pt,
					minimum size=4.5pt
				},
				recovered/.style={
					circle,
					draw=green!45!black,
					fill=white,
					line width=.7pt,
					inner sep=0pt,
					minimum size=6pt
				},
				marked/.style={
					circle,
					draw=red!80!black,
					line width=.7pt,
					inner sep=0pt,
					minimum size=3.5mm
				}
				]
				
				\begin{scope}
					\draw[step=1,gray!25,very thin]
					(0,0) grid (14.5,14.5);
					\draw[thin]
					(0,14.5)--(0,0)--(14.5,0);
					
					\foreach \k in {1,...,14}{
						\node at (\k,-.55) {$\k$};
						\node at (-.55,\k) {$\k$};
					}
					
					\foreach \k in {4,7,8,10,13,14}
					\node[marked] at (\k,-.55) {};
					
					\foreach \k in {2,3,6,8,10,14}
					\node[marked] at (-.55,\k) {};
					
					\foreach \x/\y in {5/11,8/9,9/7,11/4,14/3}{
						\draw[->,orange!85!black,thin]
						(\x,\y)--({\x-.5},{\y-.5});
						\fill[orange!85!black]
						(\x,\y) circle[radius=1.5pt];
						\node[shifted] at ({\x-.5},{\y-.5}) {};
						\node[anchor=south east]
						at ({\x-.12},{\y+.2}) {$(\x,\y)$};
					}
					
					\foreach \x/\y in {7/12,10/8,12/5}
					\node[fixed] at (\x,\y) {};
					
					\node[anchor=south east] at (6.88,12.2) {$(7,12)$};
					\node[anchor=south west] at (10.2,8.2) {$(10,8)$};
					\node[anchor=south east] at (11.88,5.2) {$(12,5)$};
				\end{scope}
				
				\begin{scope}[shift={(18,0)}]
					\foreach \r/\s in {
						4.5/10.5,
						7.5/8.5,
						8.5/6.5,
						10.5/3.5,
						13.5/2.5,
						7/12,
						10/8,
						12/5
					}{
						\fill[gray!18]
						(0,0) rectangle (\r,\s);
						\fill[gray!18]
						(\r,\s) rectangle (14.5,14.5);
					}
					
					\draw[step=1,gray!25,very thin]
					(0,0) grid (14.5,14.5);
					\draw[thin]
					(0,14.5)--(0,0)--(14.5,0);
					
					\foreach \k in {1,...,14}{
						\node at (\k,-.55) {$\k$};
						\node at (-.55,\k) {$\k$};
					}
					
					\foreach \k in {4,7,8,10,13,14}
					\node[marked] at (\k,-.55) {};
					
					\foreach \k in {2,3,6,8,10,14}
					\node[marked] at (-.55,\k) {};
					
					\foreach \x in {4,7,8,10,13,14}
					\draw[gray!65,densely dashed,very thin]
					(\x,0)--(\x,14.5);
					
					\foreach \y in {2,3,6,8,10,14}
					\draw[gray!65,densely dashed,very thin]
					(0,\y)--(14.5,\y);
					
					\foreach \x in {4,7,8,10,13,14}{
						\foreach \y in {2,3,6,8,10,14}{
							\draw[gray!70,very thin]
							({\x-.07},{\y-.07})--
							({\x+.07},{\y+.07})
							({\x-.07},{\y+.07})--
							({\x+.07},{\y-.07});
						}
					}
					
					\foreach \x/\y in {
						4.5/10.5,
						7.5/8.5,
						8.5/6.5,
						10.5/3.5,
						13.5/2.5
					}
					\node[shifted] at (\x,\y) {};
					
					\foreach \x/\y in {7/12,10/8,12/5}
					\node[fixed] at (\x,\y) {};
					
					\foreach \x/\y in {4/14,13/3,14/2}{
						\node[recovered] at (\x,\y) {};
						\fill[green!45!black]
						(\x,\y) circle[radius=.8pt];
					}
					
					\node[green!45!black,anchor=north west]
					at (4.3,13.65) {$(4,14)$};
					
					\node[green!45!black,anchor=south east]
					at (12.8,3.3) {$(13,3)$};
					
					\node[green!45!black,anchor=north east]
					at (13.8,1.7) {$(14,2)$};
				\end{scope}
				
				\node[fixed] at (2,-2.2) {};
				\node[anchor=west] at (2.5,-2.2) {fixed point};
				
				\node[shifted] at (10,-2.2) {};
				\node[anchor=west] at (10.5,-2.2) {shifted point};
				
				\node[recovered] at (19,-2.2) {};
				\fill[green!45!black]
				(19,-2.2) circle[radius=.8pt];
				\node[anchor=west] at (19.5,-2.2)
				{point in $\mathcal C_w$};
				
				\node[marked] at (28,-2.2) {};
				\node[anchor=west] at (28.6,-2.2) {marked coordinate};
				
			\end{tikzpicture}%
		}
	}
	\caption{The definition of $\mathcal C_w$ with
		$w=\jt^{-1}(p,q)$, where
		$p=11112132232312$ and
		$q=11223123212311$.}
	\label{fig:Cw}
\end{figure}

See Figure \ref{fig:Cw} as a description of the definition of $C_w$.

Define two words $p'=p'_1\cdots p'_n$ and $q'=q'_1\cdots q'_n$ by
letting
\begin{equation}\label{eq:path-pq}
\begin{aligned}
	\{i:p'_i=2\}&=\{x: (x,y) \in \E_w \cup \C_w\},
	\qquad\{i:p'_i=3\}=J_3(w),\\
	\{v:q'_v=2\}&=\{y: (x,y) \in \E_w \cup \C_w\},
	\qquad \{v:q'_v=3\}=T_3(w).
\end{aligned}
\end{equation}
Set the other entries of $p'$ and $q'$ to be $1$
and define $\sigma=\xi(w)=\pv(p',q')$. For $n=1$, 
set $\xi(H_r)=1$.

\begin{example}
	Consider the red-blue Motzkin path $w$ given in Figure \ref{fig:red-blue-path}, we have $(p,q)=\jt(w)$ as $p=123312123$ and $q=123112233.$
	Thus $J_2(w)=\{2,6,8\}$, $T_2(w)=\{7,6,2\}$,
	$J_3(w)=\{3,4,9\}$ and $T_3(w)=\{9,8,3\}$.
	Therefore, $\R_2(w)=\{(2,7),(6,6),(8,2)\}$
	and $\R_3(w)=\{(3,9),(4,8),(9,3)\}.	$
	Following the definition of $\xi$, we first split
	$\R_2(w)$ into	$\E_w=\{(2,7),(8,2)\}$ and  
	$\B_w=\{(6,6)\}.$
	Moreover, the marked $x$-coordinates are $1,5,7,9$ and the marked $y$-coordinates are $1,5,6,9$.
	Further, we have 
	$\C_w=\{(5,6),(7,5)\}.$
	Consequently, the words used to define $\xi(w)$ satisfy
	\begin{align*}
		&\{i:p'_i=2\}=\{2,5,7,8\},
		\qquad\{i:p'_i=3\}=\{3,4,9\},\\
		&\{v:q'_v=2\}=\{2,5,6,7\},
		\qquad\{v:q'_v=3\}=\{3,8,9\}.
	\end{align*}
	And hence, we obtain
	$p'=123321223,$ and $q'=123122233.$
	Finally, $\xi(w)=\pv(p',q')=479861523.$
\end{example}

\begin{lem}\label{lem:C}
	Let $\sigma\in\SI_n$ with $n \geq 2$ and let $w=f(\sigma)$. 
	Then, we have
	\[
	\C_w=\{(x,\sigma_x): \rr_x(\sigma)=2 \,\text{and}\,\sigma_x \,\text{is a right-to-left maximum of }\, \sigma\}.
	\]
\end{lem}

\pf	Firstly, distinct points of $\C_w$ have distinct vertical coordinates.
	Suppose $(t,v)$ and $(t',v')$ are two points in $\C_w$ with $t<t'$.
	Since $t<n$ is marked, we have $(t+1,y)\in\R_2(w)$ for some $y$. Since both $(t,v)$ and $(t',v')$ are not shaded, $v\geq y$ and $v'<y$. Thus $v>v'$, as desired.
	Following a similar argument, we may prove 
	two distinct points of $\C_w$ have distinct horizontal coordinates.
	
	Now let $(t,u)$ be a right-to-left maximum of $\sigma$ with rank $2$.  We prove that both $t$ and $u$ are marked.
	If $t=n$, it is immediate. Otherwise, let $(r,s)$ be the nearest right-to-left maximum to its right.
	Then, $(t+1,s+1)$ is the corresponding candidate corner. If $(t+1,s+1)$ is retained, then $t$ is marked.
	If it is rejected, there exists
	$(i,z)\in\E(\sigma)$ with
	$i\leq t+1$ and $z\leq s+1\leq u$.
	If $i\le t$, then $i<t$ and $z<u$,
	so $u$ has rank at least $3$, a contradiction.
	Consequently, $i=t+1$, and hence $t$ is marked.
	By a similar argument, we have $u$ is marked.
	
	The point $(t,u)$ is not shaded.
	It lies strictly northwest or southeast of each
	$3$-rank entry.	For a middle point $(i,z)$, if $i<t$, then $z>u$;	if $i>t$, then $z<u$.
	Thus $(t,u)$ is not in the northeast or southwest region of $(i-\frac12,z-\frac12)$.
	For each candidate corner
	$(t_j+1,u_{j+1}+1)$, either	$t\leq t_j$ and $u\geq u_{j+1}+1$, or $t\geq t_j+1$ and $u\leq u_{j+1}$.
	Hence $(t,u)$ also avoids the shaded regions
	of each shifted retained corner.
	Therefore $(t,u)\in\C_w$.
	
	Conversely, let $(t,u)\in\C_w$.
	We show that $t$ is the position of a right-to-left maximum of $\sigma$.
	Otherwise, since $\sigma_n$ is always a right-to-left maximum, we have $t<n$.	Since $t$ is marked, $(t+1,y)\in\Z(\sigma)$ for some $y$. Notice that 
	$(t+1,y)\notin\B(\sigma)$, since $t$ is not the position of a right-to-left maximum.
	Thus $(t+1,y)\in\E(\sigma)$. Suppose $(r,s)$ is the first right-to-left maximum to the	right of $t$.
	Then $r>t+1$ and $s>y$, so $s$ has rank $3$.
	Since $(t,u)$ is not shaded and $t<r$, we have $u>s$.
	This is impossible when $s=n$. When $s<n$,
	 let $(a,b)$ be the nearest right-to-left maximum to the left of $(r,s)$ in $\sigma$.
	We have $a<t<r$ and $(a+1,s+1)$ is a candidate corner.
	If $(a+1,s+1)$ is a retained corner, then
	$t\geq a+1$ and $u\geq s+1$. 
	Thus, $(t,u)$ is in the northeast region of $(a+\frac{1}{2},s+\frac{1}{2})$.
	If it is rejected, there exists $(i,z)\in\E(\sigma)$
	such that $i\leq a+1\leq t$ and $z\leq s+1\leq u$.
	Then $(t,u)$ lies in the northeast region
	of $(i-\frac12,z-\frac12)$.
	Both cases are impossible. Hence, $\sigma_t \in \Rlmaxv(\sigma)$. It can not have rank $3$, since each $3$-rank entry are shaded.
	Its rank cannot be $1$ either, otherwise $\sigma$ would have a nontrivial skew cut, a contradiction. Hence $\rr_t(\sigma)=2$.
	We have shown  that 
	\[
	\{(x,\sigma_x): \rr_x(\sigma)=2,\ \sigma_x
	\text{ is a right-to-left maximum}\}\subseteq \C_w,
	\]
hence $(t,\sigma_t)\in\mathcal C_w$.
	The uniqueness of the horizontal coordinate gives
	$u=\sigma_t$, and we complete the proof.
\qed

\begin{lem}\label{lem:path-pq-compatible}
	For $w \in \RM_{2n-1}$ with $n \geq 2$, its corresponding pair of words $(p',q')$ defined in (\ref{eq:path-pq}) is compatible.
\end{lem}

\pf For $w \in \RM_{2n-1}$,  write $\R_2(w)=\{(a_1,b_1),\ldots,(a_h,b_h)\}$ with
$a_1<\cdots<a_h$ and $b_1>\cdots>b_h$.
Set $a_{h+1}=n+1$ and $b_0=n+1$. From the geometrical definition of $\C_w$, we see  $(x,y) \in \C_w$ if and only if $(x,y)=(a_{j+1}-1,b_j-1)$ for some $0 \leq j \leq h$ and $(x,y)$ is not shaded by an entry from $\R_3(w)$.

Now write
$\Z'=\E_w\cup\C_w=\{(x_1,y_1),\ldots,(x_m,y_m)\}$
with $x_1<\cdots<x_m$ and $y_1>\cdots>y_m$.
Set $x_{m+1}=n+1$ and $y_0=n+1.$
Define
\[
\Delta(x,y)
=
\bigl|\{i<x:i\notin J_3(w)\}\bigr|
-
\bigl|\{v>y:v\notin T_3(w)\}\bigr|
\]
and
\[
\Gamma(x,y)
=
\bigl|\{i<x:i\in J_3(w)\}\bigr|
-
\bigl|\{v>y:v\in T_3(w)\}\bigr|.
\]
Let $(p,q)=\jt(w)$. By the definition of $p'$ and $q'$, we have $\{i:p'_i=3\}=J_3(w)=\{i:p_i=3\}$ and $\{v:q'_v=3\}=T_3(w)=\{v:q_v=3\}$. Upon this fact,
to prove this lemma is equal to show 
$\Delta(x_t,y_t)\geq 1$ for each $1\leq t\leq m$
and $\Gamma(x_{t+1},y_t)\leq 0$ for each $0\leq t\leq m$.

We now prove the first inequality. For $1\leq j\leq h$, let $k_j$ be
the position in $w$ of the $j$th blue horizontal step, reading from left to right. By the definitions of $p$ and $q$,
\begin{align*}
	\bigl|\{i<a_j:i\notin J_3(w)\}\bigr|
	&=1+\bigl|U_b(w_1\cdots w_{k_j-1})\bigr|+(j-1),\\
	\bigl|\{v>b_j:v\notin T_3(w)\}\bigr|
	&=\bigl|D_r(w_1\cdots w_{k_j-1})\bigr|+(j-1).
\end{align*}
Hence, by the red-blue rule,
\[
\Delta(a_j,b_j)
=1+\bigl|U_b(w_1\cdots w_{k_j-1})\bigr|
-\bigl|D_r(w_1\cdots w_{k_j-1})\bigr|
\geq 1.
\]
Since $\E_w\subseteq\R_2(w)$, this proves the first inequality for the case in $\E_w$.

Now consider $(x,y)\in\mathcal C_w$. As mentioned, $(x,y)=(a_{j+1}-1,b_j-1)$ for some $0\leq j\leq h$.
If $j \neq 0$ and $a_{j+1}-1 > a_j$, then
\[
\begin{aligned}
	&\bigl|\{i<a_{j+1}-1:i\notin J_3(w)\}\bigr|
	-\bigl|\{i<a_j:i\notin J_3(w)\}\bigr|\\
	&\qquad=1+\bigl|\{i:a_j<i<a_{j+1}-1,\ p_i=1\}\bigr|.
\end{aligned}
\]
Also, we have
\[
\bigl|\{v>b_j-1:v\notin T_3(w)\}\bigr|
-\bigl|\{v>b_j:v\notin T_3(w)\}\bigr|=1.
\]
Consequently,
\[
\Delta(a_{j+1}-1,b_j-1)
=
\Delta(a_j,b_j)
+\bigl|\{i:a_j<i<a_{j+1}-1,\ p_i=1\}\bigr|
\geq 1.
\]
If $j \neq 0$ and $a_{j+1}-1 = a_j$, then
$\Delta(a_{j+1}-1,b_j-1)=\Delta(a_j,b_j)-1$. 
It suffices to show $\Delta(a_j,b_j) \geq 2$.
Let $u=w_1\cdots w_{k_j-1}$ be the prefix preceding
the $j$th blue horizontal step.
Since $a_{j+1}-1=a_j$, we have
$G_j=(a_j,b_j-1)\in\C_w$.
For every $(r,s)\in\R_3(w)$ with $s>b_j$, we must
have $r<a_j$. Otherwise $G_j$ would be shaded, a contradiction.
Therefore, by the reading rules for $(p,q)$,
\[
|U_r(u)|
=
|\{(r,s)\in\R_3(w):s>b_j\}|
\leq
|\{(r,s)\in\R_3(w):r<a_j\}|
=
|D_b(u)|.
\]
By the red-blue rule, we further have
$|U_r(u)|=|D_b(u)|$.
Since $w\in\RM_{2n-1}$, the $j$th blue horizontal
step has positive height. Consequently,
\[
|U_b(u)|-|D_r(u)|
=
|U_b(u)|+|U_r(u)|-|D_b(u)|-|D_r(u)|
\geq1.
\]
It follows that
$\Delta(a_j,b_j)=1+|U_b(u)|-|D_r(u)|\geq2$, as desired.

If $j=0$, we have $(x,y)=(a_1-1,n)$. We claim that 
$p'_1=1$, namely, $a_1 \neq 2$. Otherwise, assume that $a_1=2$. Write $w=(B_0H_rR_0)H_b(B_1H_rR_1)H_b \cdots H_b (B_hH_rR_h)$.
Since $a_1=2$, then $B_0=\emptyset$.  Combining with the red-blue rule of $w$, there is no $D_r$ in $R_0$. 
Notice that each $H_b$ can not be at $x$-axis. Then, 
there is at least one $U_r$ in $R_0$. This means that there are at least one $3$ and no $1$ after the last $2$ in $q$.
Thus, $n \in T_3(w)$, and hence $(x,y)=(1,n)$ will be shaded by some $(r,n)$, a contradiction. The claim is verified.
Combining this claim with the fact that there are no values greater than
$n$,  we deduce $\Delta(x,y)\geq 1$. 

We next prove the second inequality. Suppose, to the contrary, that $\Gamma(x_{t+1},y_t)>0$
for some $0\leq t\leq m$.  This means that there is a point $(r,s)\in\R_3(w)$ such that $r<x_{t+1}$ and $s\leq y_t$.
Since $(p,q)$ is compatible, each point in $\R_3(w)$ has a point of $\R_2(w)$ strictly to the southwest of it. Therefore, there exists some $(a,b) \in \R_2(w)$ with
$a<r$ and $b<s$. So $(a,b)\in\E_w$. Hence $(a,b)\in\Z'$ with
$a<x_{t+1}$ and $b<y_t$. If $1\leq t<m$, this is impossible since $(x_t,y_t)$ and $(x_{t+1},y_{t+1})$ are consecutive in $\Z'$ and $b<y_t$.
If $t=0$, the same point $(a,b)\in\E_w\subseteq\Z'$ exists with $a<x_1$, a contradiction with the minimality of $x_1$.
If $t=m$, it would have $b<y_m$, and hence $a>x_m$, a contradiction with the maximality of $x_m$.
Therefore, $\Gamma(x_{t+1},y_t)\leq 0$ for $0\leq t\leq m$.
This completes the proof.
\qed

\begin{lem}\label{lem:norank1}
	Let $w\in\RM_{2n-1}$ with $n\geq 2$, and let
	$\sigma=\xi(w)$. Then $\sigma$ has no right-to-left maximum with rank $1$.
\end{lem}

\pf
Let $(p,q)=\jt(w)$ and $(p',q')$ be the words in (\ref{eq:path-pq}).
Write $\R_2(w)=\{A_1,\ldots,A_h\}$, where $A_i=(a_i,b_i)$ with $a_1<\cdots<a_h$ and $b_1>\cdots>b_h$. Set $a_{h+1}=n+1$, $b_0=n+1$ and
$G_j=(a_{j+1}-1,\,b_j-1)$ for $0\leq j\leq h$.
We shall use the following facts.
\begin{itemize}
	\item[(i)] By Lemmas \ref{lem:rb-rank}, \ref{lem:rb-bijection} and \ref{lem:path-pq-compatible},
	$\sigma\in\Sym_n(1234)$, the set of points of $\sigma$ with rank $2$ is $\E_w\cup\C_w$, and the set of points of
	$\sigma$ with rank $3$ is $\R_3(w)$. In particular, $\rr_i(\sigma)=p'_i$ and $\rr_{\sigma^{-1}(v)}(\sigma)=q'_v$.
	\item[(ii)] As shown in the proof of Lemma \ref{lem:path-pq-compatible}, $G_j\in\C_w$ if and only if $G_j$ lies
	strictly northwest or strictly southeast of every point of $\R_3(w)$.
	\item[(iii)] Since $(p,q)$ is compatible, every point of $\R_3(w)$ has a point of $\E_w$ strictly to its southwest. Moreover, $\gamma_0=0$ and $\alpha_t>\beta_t$ for $1\leq t\leq h$.
	\item[(iv)] As in the proof of Lemma \ref{lem:maptoRM}, the height of the $t$-th $H_b$-step of $w$ is
	$(\alpha_t-\beta_t-1)+(\delta_t-\gamma_{t-1})$.
\end{itemize}

Suppose, to the contrary, that $\sigma$ has a right-to-left maximum of rank $1$, and let $(m,v)$ be the leftmost one.
We consider two cases.

If $m=1$, then $\sigma_1=n$ and hence $q'_n=1$.
When $n\in T_3(w)$, by (i) we have $q'_n=3$, a contradiction. When $n\notin T_3(w)$, for $(x,y) \in \R_3(w)$, we have $y<n$.
Since $\gamma_0=0$, $x>a_1$.
Hence $G_0=(a_1-1,n)$ lies strictly northwest of every
point of $\R_3(w)$. By (ii), we have  $G_0\in\C_w$.  Thus
$q'_n=2$, a contradiction.

If $m\geq 2$,  set $c=m-1$ and $y=v+1$. Since $\sigma_m$ is both a left-to-right minimum and a
right-to-left maximum,  $\{\sigma_1,\ldots,\sigma_c\}=\{y,\ldots,n\}$ and every point of $\sigma$ other than $(m,v)$
lies in
\[
N=[1,c]\times[y,n]\qquad\text{or}\qquad S=[m+1,n]\times[1,v].
\]
By (i), we obtain:
\begin{itemize}
	\item[(1)] Every point of $\E_w\cup\C_w\cup\R_3(w)$ lies in $N\cup S$.
	\item[(2)] The value $y$ lies in $N$ and every entry to its left is greater than $y$, so $q'_y=1$.
	\item[(3)] Comparing the numbers of $3$-rank entries of $\sigma$ in $N$, we have $|\{i\leq c:p'_i=3\}|=|\{j\geq y:q'_j=3\}|$. By (i), $|\{i\leq c:p_i=3\}|=|\{j\geq y:q_j=3\}|$, and hence $|\{i\leq c:p_i\neq 3\}|=|\{j\geq y:q_j\neq 3\}|$.
\end{itemize}

\noindent\textbf{Claim 1.} There exists $1\leq t\leq h$ such that $a_t=m$, $A_t\notin\E_w$, and $s>b_t$ for every
$(r,s)\in\R_3(w)$ with $r<m$.

Clearly, $u=\sigma_c$ is a right-to-left maximum of $\sigma$, and hence $\rr_c(\sigma) \neq 1$. We consider two cases. If $\rr_c(\sigma) = 2$, then $(c,u)\notin\E_w$ and $(c,u)\in\C_w$. By (ii), 
$(c,u)=G_{t-1}$ for some $1\leq t\leq h+1$ with $a_t=m$ and $b_{t-1}=u+1$. If $A_t \in \E_w$, then $p'_{a_t}=2$, a contradiction. Hence, $A_t \notin \E_w$. Since $m\leq n$, we have $t\leq h$.
Because $G_{t-1}\in\C_w$, every $(r,s)\in\R_3(w)$ with $r<c$ satisfies $s>u=b_{t-1}-1\ge b_t.$
If $\rr_c(\sigma) = 3$,
there is an index $j'$ with $a_{j'}<c$ and $b_{j'}<u$ by (iii). Choose the largest $j$ among such $j'$,
then $A_j\in\E_w$ and $A_j$ is a point of $\sigma$ in $N$. Thus $b_j\geq y$. By the definition of $q'$, we have $q'_{b_j}=2$, while $q'_y=1$ by (2).
Therefore, $b_j>y$. By Proposition
\ref{prop:rankeq}, every $(r,s)\in\R_3(w)$ with $r<m$ satisfies $s\geq u>b_j$. By the choice of $j$, we have $a_{j+1} \geq m$. We claim that $a_{j+1} = m$. Otherwise, 
suppose $a_{j+1}>m$. If $(r,s)\in\R_3(w)$ with $r<m$, then 
$G_j$ lies strictly southeast of $(r,s)$ by $s>b_j-1$.
If $r>m$, then $r> a_{j+1}$. Otherwise, there is no $A_i$ strictly to the southwest of $(r,s)$, by (iii) we obtain a contradiction.
Thus, in this case $G_j$ lies strictly northwest of $(r,s)$.
By (ii), $G_j\in\C_w$, but
$G_j\in[m,n]\times[y,n]$, contrary to (1). Therefore $a_{j+1}=m$ and the claim is verified. Take $t=j+1$. Then every $(r,s)\in\R_3(w)$ with $r<m$ satisfies $s>b_j>b_t$.
Since $p'_m=1$, $A_t \notin \E_w$.

\noindent\textbf{Claim 2.} $b_t=y$.

Suppose first that $b_t>y$. Since $a_{t+1}-1\geq m$ and $b_t-1\geq y$, we have $G_t \in [m,n] \times [y,n]$.
Put $(r,s)\in\R_3(w)$. If $r<m$, then by Claim 1 we have $s>b_t>b_t-1$ , so $G_t$ lies strictly to the southeast  of $(r,s)$. The case $r=m$ is impossible by $p'_m=1$.
If $r>m$, then $r>a_{t+1}$. If not, there is no $A_i$ lying strictly to the southwest of $(r,s)$, a contradiction with (iii). Since $s<y\leq b_t-1$,  $G_t$ lies strictly northwest of $(r,s)$. By (ii),
$G_t\in\C_w\cap([m,n]\times[y,n])$, a contradiction with  (1).

Suppose next that $b_t<y$. Put  $(r,s)\in\R_3(w)$. Then, 
$s < b_t$ or $s \geq y$. Otherwise, by (1) $(r,s)$ would lie in $S$, so $r>m=a_t$. Thus $A_t\in\E_w$, contrary to
Claim 1. Since $a_t=m$, we have
\[
|\{i\leq c:p_i\neq 3\}|=\alpha_t+(t-1),
\]
and
\[
|\{v\geq y:q_v\neq 3\}|
=|\{v>b_t:q_v\neq 3\}|-(y-1-b_t)
=\beta_t+(t-1)-(y-1-b_t).
\]
By (3), $\alpha_t=\beta_t-(y-1-b_t)\leq\beta_t$, which contradicts (iii). This proves Claim 2.

Combining Claims 1 and 2, we have $A_t=(m,y)=(c+1,n-c+1)$. Since  $q_{b_t}=2$, by (3) we have
\[
\gamma_{t-1}=|\{i<a_t:p_i=3\}|=|\{v>b_t-1:q_v=3\}|=\delta_t.
\]
 Moreover,
\[
a_t-1=\alpha_t+(t-1)+\gamma_{t-1},
\qquad
n-b_t+1=\beta_t+t+\delta_t .
\]
Since $a_t-1=c=n-b_t+1$, we have $\alpha_t-\beta_t-1=0$ and $\delta_t-\gamma_{t-1}=0$. By (iv), the $t$-th $H_b$-step of $w$ lies on the $x$-axis,
a contradiction with the fact $w\in\RM_{2n-1}$. This completes the proof.
\qed

\begin{lem}\label{lem:path2per-BCE}
	For $w \in \RM_{2n-1}$ with $n\geq 2$ and $\sigma=\xi(w)$, we have
		$\E(\sigma)=\E_w$,
		 \begin{equation}\label{eq:Cw}
		 \{(x,\sigma_x): \rr_x(\sigma)=2 \,\text{and}\,\sigma_x \,\text{is a right-to-left maximum of }\, \sigma\}=\C_w,
		 \end{equation}
 and $\B(\sigma)=\B_w$.
\end{lem}	

\pf Recall that in the construction of $\xi$, we have 
$\sigma=\xi(w)=\pv(p',q')$. By (\ref{eq:path-pq}) and the definition of $\E_w$ and $\C_w$, we have $\E_w =\E(\sigma)$.

Let $\R_2(\sigma)=\{(x,\sigma_x):\rr_x(\sigma)=2\}.$
By the definition of $\xi$, we have $\R_2(\sigma)=\E_w\cup\C_w.$
Clearly, $\E_w\cap\C_w=\varnothing$.
Since $\E(\sigma)=\E_w$, it follows that
$\C_w=\R_2(\sigma)\setminus\E(\sigma).$
A $2$-rank entry of $\sigma$ is either a midpoint or a right-to-left maximum of $\sigma$. Thus,
\[
\C_w
=
\{(x,\sigma_x):r_x(\sigma)=2
\text{ and }\sigma_x\text{ is a right-to-left maximum of }\sigma\}.
\]

Now, we want to show $\B(\sigma)=\B_w$. 
As in the proof of Lemma \ref{lem:norank1}, write $\R_2(w)=\{A_1,\ldots,A_h\}$, where
$A_i=(a_i,b_i)$ with $a_1<\cdots<a_h$ and $b_1>\cdots>b_h$. Set $a_{h+1}=n+1$, $b_0=n+1$ and 
$G_j=(a_{j+1}-1,\,b_j-1)$ for $0\leq j\leq h$.


Let $M$ be the set of points $(t,\sigma_t)$ such that $t \in \Rlmaxp(\sigma)$. By Proposition \ref{prop:rankeq}, Lemma \ref{lem:norank1}, (\ref{eq:Cw}) and (i), we have $M=\C_w\cup\R_3(w)$.
By (ii), every point of $\C_w$  arises from $G_j$. We need the following observation.

\noindent\textbf{(O)} Let $L=(\ell,u)$ and $R=(r,v)$ be consecutive points of $M$
with $\ell<r$. Then no $G_j$ lies in the open rectangle
$D=(\ell,r)\times(v,u).$

Suppose to the contrary that  there is some $G_j$ lying in  $D$. Since $\R_3(w)\subseteq  M$ and every point of $M$ either has position at most $\ell$ and value at least
$u$, or has position at least $r$ and value at most $v$,
$G_j$ lies
strictly northwest or strictly southeast of every point of $\R_3(w)$. By fact (ii), 
$G_j\in\C_w\subseteq M$. Therefore, this point has position strictly between
$\ell$ and $r$, contradicting to the fact that $L$ and $R$ are consecutive.

%
%

We next show that $\B(\sigma) \subseteq \B_w$.
Suppose $(x,y)\in\B(\sigma)$. Then there are consecutive points $L=(\ell,u)$ and $R=(r,v)$ in
$M$ such that $\ell<r$, $u>v$, $x=\ell+1$ and $y=v+1$. 
Since $\E(\sigma)=\E_w$, we have
\begin{equation}\label{eq:retained-w}
	\E_w\cap\bigl([1,x]\times[1,y]\bigr)=\varnothing .
\end{equation}
Since $L$ and $R$ belong to $\C_w$ or $\R_3(w)$, there are four cases to consider. For each case, we show that $(x,y)=A_t \in \B_w$ for some $1\leq t\leq h$.

If $L\in\C_w$ and $R\in\C_w$, 
write $L=G_i$ and $R=G_k$ with $i<k$. If $k\geq i+2$, then  $G_{i+1}$ lies in $D$, contrary to (O). Therefore $k=i+1\leq h$ and
$(x,y)=(\ell+1,v+1)=(a_{i+1},b_{i+1})=A_{i+1}.$
Notice that every point in $\R_3(w)$ lies strictly to the southeast of $R$ or strictly  to the northwest of $L$. 
This gives $A_{t} \in \B_w$ with $t=i+1.$

If $L\in\C_w$ and $R\in\R_3(w)$,
write $L=G_{t-1}$. Then $a_t=\ell+1=x$ and $b_{t-1}=u+1$. 
Since $a_t=\ell+1\leq r\leq n$,
we have $t\leq h$.  It suffices to prove $y=b_t$.
If $A_t \in \E_w$, (\ref{eq:retained-w}) gives $b_t>y$.
Since every point in $\R_3(w)$ either lies to the  
southeast of $R$ or strictly to the northwest of $L$, there is no point of $\R_3(w)$ to the northeast of $A_t$, against the definition of $\E_w$. Hence, $A_t \in B_w$ with $b_t>v=y-1$. 
Following (iii), there is at least one point of $\E_w$ lying strictly to the southwest of $R$. Choose the leftmost one as $(m,n)$. In view of (\ref{eq:retained-w}), $m>a_t$. 
It follows that the first coordinate of $G_t$, which is $a_{t+1}-1$, is not larger than $m-1$ and its second coordinate is $b_t-1$. By (O), we have $b_t-1 \leq y-1$ or $b_t-1 \geq u$.
Combining with
the fact that $y-1<b_t \leq b_{t-1}-1=u$, we obtain $b_t=y$, as desired.

If $L\in\R_3(w)$ and $R\in\C_w$, write $R=G_t$. 
If $t=0$, then $L \in \R_3(w)$ lies to the left of $(a_1,b_1)$, which is against (iii). Thus, $1 \leq t \leq h$. Moreover, $a_{t+1}=r+1$ and $b_t=v+1=y$. 
We claim that $A_t \in \B_w$. Otherwise, by (\ref{eq:retained-w}) we have $a_t>x$. Since each point of $\R_3(w)$ either lies to the northwest of $L$ or strictly to the southeast of $R$, there is no point in $\R_3(w)$ lying to the northeast of $A_t$, contrary to the definition of $\E_w$. The claim is verified. It suffices to show $a_t=x$. By (iii), there is at least one point in $\E_w$ to the southwest of $L$.  Choose the rightmost one as $(m,n)$ and  we have $t \geq 2$. In view of (\ref{eq:retained-w}),
$m<x$ and $n>y$. The first coordinate of $G_{t-1}$ is $a_t-1$ and its second coordinate is not larger than $n-1$.
By (O), $a_t-1$ is either not smaller than $r$ or not larger than $\ell=x-1$. By the definition of $\B_w$, $a_t \geq \ell$. Since the first coordinates of $\R_2(w)$ and $\R_3(w)$ do not intersect, $a_t > \ell$.
Combining with $a_t \leq a_{t+1}-1 =r$, we deduce that $a_t=x$, as desired.

Consider the remaining case  $L\in\R_3(w)$ and $R\in\R_3(w)$. By (iii), there is at least one point in $\E_w$ strictly to the southwest of $L$, as well as $R$.
Choose the rightmost one $(m_1,n_1)$ of $L$ and the leftmost one $(m_2,n_2)$ of $R$. By $(x,y)$ is retained and (\ref{eq:retained-w}),  $m_1<x<m_2$ and $n_1>y>n_2$. If $(m_1,n_1)$ and $(m_2,n_2)$ is consecutive in $\R_2(w)$,
then there is some $G_j=(m_2-1,n_1-1)$. Since $r-1>m_2-1 \geq x$ and $u-1 > n_1-1 \geq y$, then $G_j$ lies in $D$, a contrary to (O). Combining with the choices of $(m_1,n_1)$ and $(m_2,n_2)$ in $\E_w$, there is no point in $\E_w$ and at least one point in $\B_w$ between  $(m_1,n_1)$ and $(m_2,n_2)$.
By the fact that the points in $\R_2(w)$ increase in the first coordinate and decrease in the second coordinate, and the definition of $\B_w$, all such points lie in the open rectangle $(x-1,m_2) \times (y-1,n_1)$.
 Indeed, there is only one such point in $\B_w$.
Otherwise, there will be some $G_j$ lying in the open rectangle $(x-1,m_2) \times (y-1,n_1)$, and hence in
$D$, against (O). We write it as $A_t$ with $2 \leq t \leq h-1$. Again by (O), $G_t=(m_2-1,b_t-1)$ and $G_{t-1}=(a_{t}-1,n_1-1)$ lies out of $D$. Hence, $b_t-1 \leq v=y-1$ and $a_t-1 \leq l=x-1$. Since $a_t>x-1$ and $b_t>y-1$, we have $a_t=x$ and $b_t=y$, as desired.

We now proceed to show  $\B_w\subseteq\B(\sigma)$.
Assume that $A_t\in\B_w$. Since $p_1=q_1=1$ and
$\B_w\subseteq\R_2(w)$, we have $a_t \geq 2 $ and $b_t\geq2$. We prove three claims.

\noindent\textbf{(C1)} No point of $M$ lies in $Q=[a_t,n]\times[b_t,n]$.

Since $A_t\notin\E_w$ and the first (resp. second) coordinates of $\R_2(w)$ and $\R_3(w)$ do not intersect, each point $(r,s)\in\R_3(w)$ satisfies that $r<a_t$ or $s<b_t$. For $G_j \in C_w$ and $j \leq t-1$, we have $a_{j+1}-1 \leq a_t-1<a_t$. Hence, we need only take care of the case $j \geq t$. In this case, $a_{j+1}-1 \geq a_{t+1}-1 \geq a_t$ and $b_j -1\leq b_t-1<b_t$. This confirms the claim.

\noindent\textbf{(C2)} There exists some point $L=(\ell,u)$ in $M$ satisfies that $\ell=a_t-1$ and $u\geq b_t$.

If $a_t-1 \in J_3(w)$, 
put $u=\sigma_{a_t-1}$. By (iii), $u>b_{t-1}>b_t$ and
the claim  holds.
Otherwise, when $t=1$, it is obvious that $G_0=(a_1-1,n) \in C_w$. When $t>1$, consider $G_{t-1}=(a_{t}-1,b_{t-1}-1)$. Since $(a_t,b_t) \in \B_w$, $G_{t-1}$ is to the northwest of each point in $\R_3(w)$ lying to the right of $A_t$, if any. Assume the rightmost point in  $\R_3(w)$ to the left of $A_t$ is $(r,v)$, if any. If $a_{t-1}>r$, then $b_{t-1}<v$. Hence, $G_{t-1}$ lies to the southeast of all the points in $\R_3(w)$ lying to the left of $A_t$. Hence, by (ii) $G_{t-1} \in \C_w$. If $a_{t-1}<r$, we also have $G_{t-1}$ lies to the southeast of all the points in $\R_3(w)$ lying to the left of $A_t$ through
$a_t-1 \notin J_3(w)$ and $b_{t-1}<v$. Then, (ii) gives $G_{t-1} \in \C_w$. In view of $b_{t-1}-1 \geq b_t$, the claim is verified.

\noindent\textbf{(C3)} There exists some  point $R=(r,v)\in M$ has $v=b_t-1$ and $r\geq a_t$.

The proof is analogous to that of (C2), and hence it is omitted.

Let $L$ and $R$ be the points in (C2) and (C3). Then $\ell<a_t\leq r$ and $u\geq b_t>v$. By (C1) and the fact that the first coordinate of $M$ increases and the second coordinate of $M$ decreases,  $L$ and $R$ are consecutive points in $M$. Hence, $A_t=(a_t,b_t)=(\ell+1,v+1)$ is a candidate corner. Since $\R_2(w)=\B_w \cup \E_w$, $\E_w=\E(\sigma)$ and the fact that the first coordinate of $\R_2(w)$ increases and its  second coordinate decreases, we have there is no point in $\E(\sigma)$ lying in $[1,a_t] \times [1,b_t]$. Thus, $A_t$ is retained, and hence $A_t \in \B(\sigma)$. 

In conclusion, we have $\B(\sigma)=\B_w$. This completes the proof.
\qed

\begin{lem}\label{lem:maptoSI}
	For $w \in \RM_{2n-1}$, we have $\xi(w) \in \SI_n$.
\end{lem}

\pf
	Let $\sigma=\xi(w)$. By Lemmas \ref{lem:rb-rank} and \ref{lem:path-pq-compatible},
	$\sigma$ avoids $1234$. Suppose to the contrary that 
    $\sigma$ has a skew cut $c$ with
	$1\leq c<n$, namely,  $\{\sigma_1,\ldots,\sigma_c\}=\{n-c+1,\ldots,n\}$.
	Set $\tau^1=\sigma_1\cdots\sigma_c$ and $\tau^2=\sigma_{c+1}\cdots\sigma_n$. 
	By the definition of skew sum, we have $\rr_{\tau^1}(x)=\rr_{\sigma}(x)$ and $\rr_{\tau^2}(y)=\rr_{\sigma}(y)$ for $x \in \{\sigma_1,\ldots,\sigma_c\}$ and $y \in \{\sigma_{c+1},\ldots,\sigma_n\}$. By reduction, $\tau^1$ corresponds to a $1234$-avoiding permutation with length $c$. It implies that for each $k\in\{1,2,3\}$,
	\[
	\bigl|\{i\leq c:p'_i=k\}\bigr|
	=
	\bigl|\{v>n-c:q'_v=k\}\bigr|.
	\]
	
	We now wish to show $(c+1,n-c+1) \in \B_w$. In view of Lemma \ref{lem:path2per-BCE}, it suffices to show  $(c+1,n-c+1) \in \B(\sigma)$.
	
By the structure of $\sigma$,
$c \in \Rlmaxp(\sigma)$
and $n-c \in \Rlmaxv(\sigma)$.
Thus $(c,\sigma_c)$ and $(\sigma^{-1}(n-c),n-c)$ are
consecutive right-to-left maxima, and $(c+1,n-c+1)$ is
a candidate corner. The only points in
	$[1,c+1]\times[1,n-c+1]$ are
	$(\sigma^{-1}(n-c+1),n-c+1)$ and $(c+1,\sigma_{c+1})$.
	Both correspond to left-to-right minima.
    Thus, neither point is therefore a middle point, so
	$(c+1,n-c+1)$ is retained, and
	$(c+1,n-c+1)\in\B(\sigma)=\B_w$.
	
	Write $(p,q)=\jt(w)$ and
	$\R_2(w)=\{A_1,\ldots,A_h\}$ with
	$A_i=(a_i,b_i)$ with
	$a_1<\cdots<a_h$ and $b_1>\cdots>b_h$.
	There exists an index $1\leq t\leq h$ with
    $a_t=c+1$ and $b_t=n-c+1$.
	By the definition of $\xi$, $p_i=3$ if and only if $p'_i=3$ and $q_v=3$ if and only if $q'_v=3$.
	Since $q_{b_t}=2$, 
	\[
	\begin{aligned}
		\gamma_{t-1}
		&=|\{i\leq c:p_i=3\}|
		=|\{i\leq c:p'_i=3\}|\\
		&=|\{v\geq b_t:q'_v=3\}|
		=|\{v>b_t:q_v=3\}|
		=\delta_t.
	\end{aligned}
	\]
	On the other hand, counting the letters in the
	prefix $p_1\cdots p_{a_t-1}$ and the suffix
	$q_{b_t}\cdots q_n$, respectively, gives
    $a_t-1=\alpha_t+(t-1)+\gamma_{t-1}
	$ and $n-b_t+1=\beta_t+t+\delta_t.$
	Since $a_t-1=n-b_t+1=c$ and
	$\gamma_{t-1}=\delta_t$, we obtain
	$\alpha_t-\beta_t-1=0$.
    By (iv), the height of the $t$-th $H_b$-step in $w$ is
	$(\alpha_t-\beta_t-1)+(\delta_t-\gamma_{t-1})=0.$
	This contradicts $w\in\RM_{2n-1}$.
	Hence, $\sigma\in\SI_n$.
	\qed

\noindent { \textbf{Proof of Theorem \ref{thm:odd-indecomposable}}}

\pf The case $n=1$ is immediate. For $n\geq2$,
Lemmas~\ref{lem:maptoRM} and~\ref{lem:maptoSI}
show that $f$ and $\xi$ are well-defined.
By Lemmas~\ref{lem:BEZ} and~\ref{lem:C},
the construction of $\xi(f(\sigma))$ recovers
the rank words of $\sigma$, so
$\xi(f(\sigma))=\sigma$ for every $\sigma\in\SI_n$.
Conversely, Lemma~\ref{lem:path2per-BCE},
together with the preservation of the rank-$3$
positions and values, shows that the construction
of $f(\xi(w))$ recovers $\jt(w)$.
Hence $f(\xi(w))=w$ for every $w\in\RM_{2n-1}$.
Thus $f$ and $\xi$ are inverse to each other.
\qed

Given $\pi \in \Sym_n(1234)$, we now write 
\begin{equation*}
	\pi=\tau^{(1)}\ominus\tau^{(2)}\ominus\cdots\ominus\tau^{(k)}
\end{equation*}
into skew-indecomposable permutations with 
$\Cut(\pi)=\{c_1,c_2,\ldots,c_{k-1}\}.$
Clearly, $|\tau^{(i)}|=c_{i}-c_{i-1}$, where $1 \leq i \leq k$, and we set $c_0=0$ and $c_k=n$.
By the decomposition of $\pi$, we have 
$\pi$ avoids $1234$-pattern exactly when each $\tau^{(i)}$ does.

Define
\begin{equation*} 
	w=F(\pi)
	=f(\tau^{(1)})\,H\,
	f(\tau^{(2)})\,H\cdots H\,
	f(\tau^{(k)}).
\end{equation*}
The length of $w$ is $\sum_{i=1}^k(2(c_i-c_{i-1})-1)+(k-1)=2n-1.$
Since $f(\tau^{i})$ has length $2(c_i-c_{i-1})-1$ and the number of up steps equals the number of down steps, then
each separating $H$ in $w$ is blue and lies at height zero.

Given $w \in \M_{2n-1}$, split $w$ at every blue horizontal step on the $x$-axis.
This gives the unique decomposition
$w=w^{(1)}H_b w^{(2)}H_b\cdots H_b w^{(k)}$,
where $w^{(i)}\in\mathcal{RM}_{2n_i-1}$.
  Let the \emph{ground-level separator set} of $w$ be \[\Sep(w)=\{n_1,n_1+n_2,\ldots,n_1+n_2+\cdots+n_{k-1}\}.\]
Define 
\[\pi=\Xi(w)=\xi(w^{(1)})\ominus\xi(w^{(2)})\ominus\cdots \ominus\xi(w^{(k)}).\] Clearly, the length of $\pi$ is $\sum_{i=1}^k n_i=n$ and $\pi$ avoids $1234$ by the fact each $\xi(w^{(i)})$ avoids $1234$.

\noindent { \textbf{Proof of Theorem \ref{thm:bijodd-alternative}}}

\pf
Since $f$ and $\xi$ are inverses to each other, it is routine to check that $F$
and $\Xi$ are inverses to each other and the statistics are kept.
\qed

\section{Proofs of Theorems \ref{thm:perdecom} and \ref{thm:percore-baxter}}
\label{sec:1234-core-baxter}
In this section, we give a direct decomposition of $\T_{n+1}$ into an
explicit permutation core and a subset of $[n-1]$, and further establish a bijection 
between these permutation cores and non-intersecting triples of Lattice paths.

\begin{lem}\label{lem:rank-preservation}
Given $\pi=\pi_1\pi_2 \cdots \pi_n \in \Sym_n(1234)$ and $D=\{d_1,\ldots,d_s\}\subseteq V_3(\pi)$, 
let $\rho$ be the permutation obtained by deleting the element in $D$ and then applying $\std$ to each remaining
entry, where $\std:[n]\setminus D\to[n-|D|]$ is the
order-preserving map.  Then
$i\in V_j(\pi)$ if and only if 
$\std(i)\in V_j(\rho)$
with $j=1,2,3$.
\end{lem}
\pf
For a value $x$ of a permutation $\sigma$, let
$l_{\sigma}(x)$ denote the length of a longest increasing
subsequence of $\sigma$ ending at $x$. It is immediate that the rank of $x$ in $\sigma$
is $l_{\sigma}(x)$.

Now fix $i\in[n]\setminus D$, and let
$a_1<a_2<\cdots<a_{l-1}<i$
be the values of a longest increasing subsequence of $\pi$ ending
at $i$, where $l=l_{\pi}(i)$. None of
$a_1,\ldots,a_{l-1}$ belongs to $D$. Indeed, if
$a_t\in D$, then $a_t\in V_3(\pi)$. Then $a_1, a _2, \ldots, a_t,i$ will contain an increasing subsequence of $\pi$ ending at $i$ with length $4$. This contradicts that $\pi$ is $1234$-avoiding.
Since $\std$ preserves relative order, we obtain the increasing
subsequence
$\std(a_1)<\std(a_2)<\cdots<\std(a_{r-1})<\std(i)$
in $\rho$. Therefore, $l_{\rho}(\std(i))\geq l_{\pi}(i).$

Conversely, every increasing subsequence of $\rho$ ending at
$\std(i)$ pulls back under $\std^{-1}$ to an increasing
subsequence of $\pi$ ending at $i$. Consequently,
$l_{\rho}(\std(i))\leq l_{\pi}(i).$
It follows that
$l_{\rho}(\std(i))=l_{\pi}(i)$, as desired.
\qed

\begin{lem}\label{lem:third-rank-pattern}
For $\pi\in\T_m$, two $3$-rank entries are adjacent if and only if $\pi$
contains $12\underline{43}$. Moreover, the last entry of $\pi$ is not a
$3$-rank entry if and only if it equals $2$.
\end{lem}

\pf Let $x,y$ be two adjacent $3$-rank entries.  By
Proposition \ref{prop:rankeq}, we have $x>y$.  Choose an increasing
triple $u<v<y$ ending at $y$.  Since $x>y$ and $x,y$ are adjacent, both
$u$ and $v$ occur before $x$.  Hence $u,v,x,y$ form an occurrence of
$12\underline{43}$.

Conversely, suppose that $u,v,x,y$ form an occurrence of
$12\underline{43}$.  Then $u<v<y<x$ with $x,y$ adjacent. So the rank of 
$x$ and $y$ are both at least $3$ in $\pi$. In fact, they are
exactly $3$ because $\pi$ avoids $1234$.

Finally, since $1$ occurs to the left of $2$ in $\pi$, every entry to the right of $2$ exceeds
 $2$.
 Hence, if the last entry is not $2$, it must be a $3$-rank entry; if it is $2$, it is a $2$-rank entry.
 \qed

By Lemma~\ref{lem:third-rank-pattern}, $\K_k$ consists precisely of those permutations in $\T_{k+2}$ whose $3$-rank entries are pairwise nonadjacent and whose last entry is not a $3$-rank entry.

To prove Theorem \ref{thm:perdecom},
we first define a compression operation on $\T_{n+1}$.  Given $\pi \in\T_{n+1}$, consider each
maximal consecutive block of $3$-rank entries. From every block not containing $\pi_{n+1}$,
retain only the last entry; the block containing $\pi_{n+1}$, if it exists, is deleted 
entirely.
 Standardize the retained word and denote the result
by $\kappa$.  Write the set of the retained values as
$\{1,2\}\cup\{s+2:s\in S\}$ with $S\subseteq[n-1]$
and set $k=|S|$.

The entries $1$ and $2$ have ranks $1$ and $2$, respectively, and are therefore always retained.  Lemma~\ref{lem:rank-preservation} shows that the ranks of
the retained entries are unchanged.  Since all $1$-rank and $2$-rank
entries survive, no two $3$-rank entries are adjacent.  Moreover, every entry to the right of $2$ is a $3$-rank entry, so deleting the terminal $3$-rank block makes $2$ the last entry of $\kappa$.  It follows from Lemma 
\ref{lem:third-rank-pattern} that $\kappa\in\K_k$.  

We next describe a recovery operation, which will be proved to be the  inverse of the compression operation.
Fix $S\subseteq[n-1]$ and $\kappa\in\KK_{|S|}$, and let
$ A=\{1,2\}\cup\{s+2:s\in S\}$ and 
$ X=[n+1]\setminus A$.

Write $A=\{a_1<a_2<\cdots<a_{k+2}\}$ and replace each value $i$ of
$\kappa$ by $a_i$.  This gives an unstandardized retained permutation $q$.
Let its $3$-rank entries, in order of occurrence, be
$t_1>t_2>\cdots>t_h.$ Set
$t_0=n+2$ and $t_{h+1}=2$.
For $1\leq j\leq h+1$, let $D_j$ be the decreasing word formed by the
elements of
\begin{equation}\label{eq:3rankseq}
 \{x\in X:t_j<x<t_{j-1}\}.
 \end{equation}
For $1\leq j\leq h$, insert $D_j$ immediately before $t_j$, and append
$D_{h+1}$ to the end of $q$. Denote the resulting permutation by $p$.

Theorem \ref{thm:perdecom} now follows directly from the next result.

\begin{thm}
\label{thm:direct-1234-core}
For every $n\geq1$, the compression above is a bijection from $\T_{n+1}$ to \[\bigsqcup_{k=0}^{n-1}\binom{[n-1]}{k}\times\K_k.\]
\end{thm}

\pf The compression has already been shown to produce a pair
$(S,\kappa)$ in the stated disjoint union.  It remains to show that $p \in \T_{n+1}$ and the two operations are mutually inverse.

We first claim that every newly inserted value of $q$ is a right-to-left maximum of $p$. Indeed, if
$x\in D_j$ with $j\leq h$, then $x>t_j$, while every entry of $q$ to
the right of $t_j$ is smaller than $t_j$.  The entries following $x$ within the same $D_j$ are smaller than $x$, as are all elements of $D_l$ with $l>j$.  If $x\in D_{h+1}$, only smaller entries from the same
terminal decreasing word occur to its right. The claim is verified. Consequently, 
 no inserted value  can
contribute to an increasing subsequence ending at an entry of $q$, so
the ranks of all entries of $q$ remain unchanged.

We next show that $p$ avoids $1234$. Suppose, to the contrary, that the recovery operation creates an occurrence $abyx$ of $1234$ in $p$. Since every inserted value is a right-to-left maximum of $p$, the last entry $x$ must belong to some $D_k$ with $1\leq k\leq h+1$, while $a$, $b$ and $y$ are entries of $q$. Then $y$ has rank $3$ in $p$, and hence also in $q$, since the ranks of the entries of $q$ are preserved. 
But $y<x$ and $y$ precedes $x$, which contradicts the definition of the recovery operation: every $3$-rank entry of $q$ preceding $x$ in $p$ is larger than $x$.
Thus, $p$ avoids $1234$-pattern. The relative order of $1$ and 
$2$ is unchanged, so $p\in\T_{n+1}$.

It remains to verify that the two operations are inverse to each other.  
First, for any given pair $(S,\kappa)$ where $S$ is a $k$-subset of $[n-1]$ and $\kappa \in 
\K_k$, 
we aim to prove that applying the recovery operation followed by the compression operation
returns exactly $(S,\kappa)$. To this end,
we claim that  every element in $D_l$ 
with $1 \leq l \leq h+1$ has rank $3$ in $p$.
Suppose that $x\in D_j$ with $j\leq h$. Since the rank of $t_j$ in $q$ is $3$, there is an
increasing triple $u<v<t_j$ of $q$, and then $uvx$ is an increasing subsequence of $p$. It follows 
that the rank of $x$ in $p$ is at least $3$.
If $x\in D_{h+1}$, then $12x$ is an increasing subsequence of $p$, so again $x$ has rank at least $3$. Since $p$ avoids $1234$, its rank is exactly $3$.  The claim is verified.
It follows that,  the maximal consecutive block of $3$-rank entries contains no $p_{n+1}$
is exactly $D_jt_j$ with $1\leq j\leq h$, while that contains $p_{n+1}$ is $D_{h+1}$ whenever
$p_{n+1} \neq 2$.
After the deletion of all the $D_j$ with $1\leq j\leq h+1$, 
no two $3$-rank entries in $q$ are adjacent and its last entry is $2$. Standardization then recovers the pair $(S,\kappa)$.

Conversely, starting from $p\in\T_{n+1}$, we show that compression followed by recovery reproduces $p$. Let $q$ be the word obtained from $p$ by deleting each consecutive $3$-rank block as prescribed.
If $t_j$
is the retained last entry of a nonterminal $3$-rank block, then every
deleted value in that block is larger than $t_j$ and, for $j>1$, smaller 
than $t_{j-1}$, since
the $3$-rank entries of $p$ decrease from left to right.  Every value in the
terminal $3$-rank block is smaller than the final retained $3$-rank entry.
If no $3$-rank entry is retained, all deleted values belong to the terminal
block.  Hence \eqref{eq:3rankseq} uniquely determines the block containing each deleted value, and the decreasing order within each block uniquely determines the position of that value. Therefore the recovery operation reproduces $p$.\qed

\begin{example}
Let
$p = 13872654$, then $(\rr_1,\dots,\rr_8) = (1,2,3,3,2,3,3,3).$
The nonterminal $3$-rank block $87$ is replaced by its last entry $7$, and the terminal block $654$ is deleted. Hence
\[
13872654 \longmapsto 1372 \longmapsto 1342 \in \mathcal{K}_2,\quad S = \{1,5\} \subseteq [6].
\]
Conversely, $S = \{1,5\}$ and $\kappa = 1342$ give
\[
A = \{1,2,3,7\},\quad q = 1372,\quad X = \{4,5,6,8\}.
\]
The unique $3$-rank entry of $q$ is $7$. Hence, $D_1=8$ and $D_2=654$.
$D_1$ is inserted immediately before $7$, while $D_2$ are appended to the end:
\[1372 \longmapsto 13872654.\]
\end{example}

We next identify $\K_n$ with $\PP_{n}$, the family of non-intersecting triples of lattice paths of length $n$. Such path triples constitute a classical Baxter family and have been
related bijectively to Baxter permutations and several other Baxter
objects \cite{Dilks2012,DulucqGuibert1994,DulucqGuibert1998,LinLiu2023}. 
Our construction starts from the permutation class
$\Sym_n(1234,12\underline{43})$
and uses a type of ordered forest as a bridge: three statistics of
the forest define three sets, which are then encoded the path triple.

\begin{lem}
For $n \geq 1$, there is a bijection between $\K_n$ and $\Sym_{n+1}(1234,12\underline{43})$.
\end{lem}
\pf
Given $\pi=\pi_1\cdots \pi_{n+2} \in \K_n$, define a map $\delta$ such that $\delta(\pi)=\std(\pi_1\cdots \pi_{n+1})$. Since the appended $\pi_{n+2}=2$ cannot participate
 in a $1234$ occurrence or in a $12\underline{43}$ occurrence, we have $\delta(\pi) \in \Sym_{n+1}(1234,12\underline{43})$.  
The inverse of $\delta$ may be defined by sending
$q\in\Sym_{n+1}(1234,12\underline{43})$ to the permutation obtained by increasing every value
at least $2$ by $1$ and appending $2$. Then, it is routine to check that $\delta$ is a 
bijection from $\K_n$ to $\Sym_{n+1}(1234,12\underline{43})$.
\qed

It remains to construct a bijection $\Psi:\Sym_{n+1}(1234,12\underline{43})\rightarrow\PP_n.$
For a permutation $\pi=\pi_1 \pi_2 \cdots \pi_n \in\Sym_n$ and a subset
$S\subseteq[n]$, let $\pi|_S$ denote the restriction of $\pi$
to $S$, namely, the subword obtained by deleting from $\pi$ all
entries whose values do not belong to $S$.
For a permutation $\sigma$ and a value $x$ occurring in
$\sigma$, let $\pos_{\sigma}(x)$ denote the position of $x$ in the one-line notation of $\sigma$.

For ordered trees $T_1,\ldots,T_s$ with pairwise disjoint vertex
sets and a vertex $r$ not contained in any of them, we write
$[r;T_1,\ldots,T_s]$
for the ordered tree obtained by adjoining $r$ as a new root and
making the roots of $T_1,\ldots,T_s$, from left to right, the children
of $r$. 
 We write $[r]$ for the one-node
tree rooted at $r$.

We now recursively define a family $\F_n$ of admissible ordered
forests. Let
$\F_1=\{([1])\}.$
Suppose that $F_m=(T_1,\ldots,T_{M_m})\in\F_m,$
and let $c_m$ be the unique index such that the root of $T_{c_m}$
is the current maximum $m$.
For $c_m\leq p\leq M_m+1$, define the leaf extension
\[
L_p(F_m)=(T_1,\ldots,T_{p-1},[m+1],T_p,\ldots,T_{M_m}).
\]
For $1\leq s<c_m$, define the grafting operation
\[
G_s(F_m)=\bigl(
T_1,\ldots,T_{s-1},
[ m+1;T_s,T_{s+1},\ldots,T_{c_m}],
T_{c_m+1},\ldots,T_{M_m}
\bigr).
\]
We then set
\[
\F_{m+1}
=
\bigcup_{F_m\in\F_m}
\left(
\{L_p(F_m):c_m\leq p\leq M_m+1\}
\cup
\{G_s(F_m):1\leq s<c_m\}
\right).
\]
The elements of $\F_n$ are called \emph{admissible ordered
forests} on $[n]$.

Now assume that $\pi \in \Sym_{n}(1234,12\underline{43})$, we construct a forest $\varphi(\pi)$ on $[n]$ inductively.  Start with the one-node forest
$F_1=([1])$.  Suppose that the restriction $\sigma_m=\pi|_{[m]}$ has produced the
ordered forest $ F_{m}=(T_1,\ldots,T_{M_m}).$
Set $ g=\pos_{\pi|_{[m+1]}}(m+1)-1$. Let 
 $\alpha_m=\idr(\sigma_m)+1$.  Construct
$F_{m+1}$ according to the following three rules:
\begin{enumerate}
\item If $g=0$, insert the one-node tree $m+1$ immediately before $T_{c_m}$.
\item If $1\leq g\leq\alpha_m$, insert the one-node tree $m+1$ so that exactly
      $g-1$ old trees remain to its right.
\item If $g>\alpha_m$, replace
      $T_{g-\alpha_m},T_{g-\alpha_m+1},\ldots,T_{c_m}$
      by the single tree $[m+1;T_{g-\alpha_m},T_{g-\alpha_m+1},\ldots,$ $T_{c_m}]$.
\end{enumerate}

\begin{example}\label{ex:common-phi}
		Let
		$\pi=584693721		\in\Sym_9(1234,12\underline{43}).$
		The construction of $\varphi$ gives
	\[
	\renewcommand{\arraystretch}{1.25}
	\begin{array}{c|c|c|c|l}
		m&\sigma_m&\alpha_m&g&F_m\\ \hline
		1&1&1&0&
		([1])\\
		2&21&2&0&
		([2],[1])\\
		3&321&3&0&
		([3],[2],[1])\\
		4&4321&4&0&
		([4],[3],[2],[1])\\
		5&54321&5&2&
		([5],[4],[3],[2],[1])\\
		6&546321&2&4&
		([5],[4],[3],[2],[6],[1])\\
		7&5463721&2&1&
		([5],[7;[4],[3],[2],[6]],[1])\\
		8&58463721&1&4&
		([5],[7;[4],[3],[2],[6]],[1],[8])\\
		9&584693721&1&-&
		([5],[7;[4],[3],[2],[6]],[9;[1],[8]])
	\end{array}
	\]
		Here the insertions of $2,3,4,5$ use Rule~1,
		those of $6,8$ use Rule~2, and those of $7,9$
		use Rule~3.
		Thus $F=\varphi(\pi)$ is shown below.
		
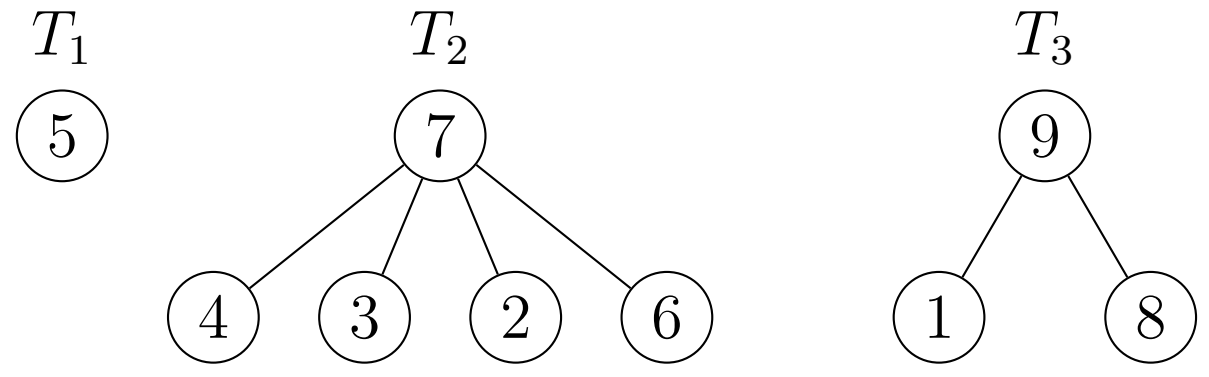
\begin{center}
			\begin{tikzpicture}[
				x=1cm,y=1cm,
				vertex/.style={
					circle,draw,fill=white,
					inner sep=2pt,minimum size=6mm
				}
				]
				\node[vertex] (v5) at (-4,0) {$5$};
				\node[vertex] (v7) at (-1.5,0) {$7$};
				\node[vertex] (v9) at (2.5,0) {$9$};
				
				\node[vertex] (v4) at (-3,-1.2) {$4$};
				\node[vertex] (v3) at (-2,-1.2) {$3$};
				\node[vertex] (v2) at (-1,-1.2) {$2$};
				\node[vertex] (v6) at (0,-1.2) {$6$};
				\node[vertex] (v1) at (1.8,-1.2) {$1$};
				\node[vertex] (v8) at (3.2,-1.2) {$8$};
				
				\foreach \v in {4,3,2,6}{
					\draw (v7)--(v\v);
				}
				\foreach \v in {1,8}{
					\draw (v9)--(v\v);
				}
				
				\node at (-4,.65) {$T_1$};
				\node at (-1.5,.65) {$T_2$};
				\node at (2.5,.65) {$T_3$};
			\end{tikzpicture}
		\end{center}
\end{example}

For convenience, let 
\[
\mu_m=\mu(\sigma_m)=
\begin{cases}
m, & \text{if $\sigma_m$ has no $3$-rank entry},\\
q_m-2, & \text{if $q_m$ is the position of the first $3$-rank
entry of $\sigma_m$}.
\end{cases}
\]
The following lemma shows that $\varphi$ is
well defined and records two identities maintained at every stage.

\begin{lem}
\label{lem:forest-well-defined}
For $1 \leq m \leq n$, we have $F_m \in \F_m$, $\mu_m=M_m$ and  $\alpha_m=M_m-c_m+1$. 
\end{lem}
\pf
We prove all three assertions simultaneously by induction on $m$. The case $m=1$ is clear.
Suppose the assertions hold at stage $m$,
we first prove $g\leq M_m$. If $\sigma_m$ has no $3$-rank
entry, then $g\leq m=\mu_m$. Otherwise, let $q_m$ be the
position of its first $3$-rank entry $x$ and there exists an increasing
pair $a<b<x$ of $\sigma_m$. If $g\geq q_m$, then
$abx(m+1)$
forms a $1234$-pattern of $\pi$. If $g=q_m-1$, then  $ab(m+1)x$
forms an occurrence of $12\underline{43}$ of $\pi$. Both are impossible,
so $g\leq q_m-2=\mu_m=M_m.$

We now distinguish the three cases in the construction. 

If $g=0$, the one-node tree is inserted immediately before
$T_{c_m}$. Hence,
$F_{m+1}=L_{c_m}(F_m)\in\F_{m+1}.$
In this case, $M_{m+1}=M_m+1$, $c_{m+1}=c_m$, 
$\alpha_{m+1}=\alpha_m+1$ and $\mu_{m+1}=\mu_m+1$.
Hence $M_{m+1}=\mu_{m+1}$ and 
$M_{m+1}-c_{m+1}+1=\alpha_m+1=\alpha_{m+1}$.

Suppose next that $1\leq g\leq\alpha_m$. The one-node tree is
inserted in position
$p=M_m-g+2.$  By the induction hypothesis, $\alpha_m=M_m-c_m+1,$
and hence $c_m+1\leq p\leq M_m+1.$
Thus, $F_{m+1}=L_p(F_m)\in\F_{m+1}.$
Moreover,
$M_{m+1}=M_m+1$, $c_{m+1}=p$, $\alpha_{m+1}=g$ and $\mu_{m+1}=\mu_m+1$.
It follows that $M_{m+1}=\mu_{m+1}$ and 
$M_{m+1}-c_{m+1}+1=M_m+1-(M_m-g+2)+1=g=\alpha_{m+1}.$

Finally, suppose that $g>\alpha_m$ and set
$s=g-\alpha_m.$
Since $g\leq M_m$ and $\alpha_m=M_m-c_m+1$, we have
$1\leq s=g-\alpha_m \leq M_m-\alpha_m=c_m-1<c_m.$
Hence, $F_{m+1}=G_s(F_m)\in\F_{m+1}.$
In this case, $\mu_{m+1}=g-1$ and $\alpha_{m+1}=\alpha_m$, and
the grafting operation gives
$M_{m+1}=M_m-c_m+s=g-1$ and $c_{m+1}=s=g-\alpha_m$.
Consequently,
$M_{m+1}=\mu_{m+1}$ and $M_{m+1}-c_{m+1}+1=(g-1)-(g-\alpha_m)+1
=\alpha_{m+1}.$ This completes the proof.
\qed

For a permutation $\sigma=\sigma_1\cdots\sigma_{m-1}\in \Sym_{m-1}(1234,12\underline{43})$
and $0\leq g\leq m-1$, let 
\[I_g(\sigma)=\sigma_1\cdots\sigma_g\,m\,\sigma_{g+1}\cdots\sigma_{m-1}.\]
We recursively define a map $\psi:\F_n\rightarrow\Sym_{n}(1234,12\underline{43}).$
Set $\psi([1])=1.$
Let $m\geq2$ and suppose that
$F_m=(T_1,\ldots,T_{K_m})\in\F_m.$
Let $T_{b_m}$ be the unique component whose root is $m$.
If $T_{b_m}=[m]$, let $F^-=(T_1,\ldots,T_{b_m-1},T_{b_m+1},\ldots,T_{K_m})$
   $\in\F_{m-1}.$
Let $b_{m-1}$ be the position in $F^-$ of the component rooted at $m-1$.
By the recursive definition of $F_m$, we have $b_m\geq b_{m-1}$ in this case.
Define
\[
g=
\begin{cases}
0, & b_m=b_{m-1},\\
K_m-b_{m}+1, & b_m>b_{m-1}.
\end{cases}
\]
If $T_{b_m}=[m;U_1,\ldots,U_r]$ with $r\geq2$,
let \[
F^-=(T_1,\ldots,T_{{b_m}-1}, U_1,\ldots,U_r, T_{{b_m}+1},\ldots,T_{K_m})
\in\F_{m-1}
\]
and define $g=K_m+1.$
In either case, set $\psi(F_m)=I_g\bigl(\psi(F^-)\bigr).$

\begin{lem}
\label{lem:forest-inverse}
The map $\psi$ is well defined, and $\mu(\psi(F_m))=K_m$
and $\idr\bigl(\psi(F_m)\bigr)=K_m-b_m$ for $1 \leq m \leq n$.
Moreover, $\varphi\circ\psi=\operatorname{id}_{\F_m}$
and $\psi\circ\varphi=\operatorname{id}_{\Sym_{m}(1234,12\underline{43})}.$
\end{lem}

\pf
We first prove, by induction on $m$, that $\psi(F_m)\in\Sym_{m}(1234,12\underline{43})$, $\mu(\psi(F_m))=K_m$ and 
$\idr(\psi(F_m))=K_m-b_m$;
we then establish the two inverse identities, again by induction on $m$.
All assertions are immediate for $m=1$.
Assume they hold for $m-1$, we show that they hold for $m$.

Suppose first that $T_{b_m}=[m]$. The forest $F^-$ has $K_m-1$
components. If $b_m=b_{m-1}$, then $g=0$ and $m$ is inserted at the beginning of 
$\psi(F^-)$, creating
neither an occurrence of $1234$-pattern nor one of
$12\underline{43}$. Thus, $\psi(F_m)\in\Sym_{m}(1234,12\underline{43})$.
Moreover,
$\mu(\psi(F_m))=\mu(\psi(F^-))+1$ and $\idr(\psi(F_m))=\idr(\psi(F^-))+1$. 
So the induction hypothesis gives $\mu(\psi(F_m))=(K_m-1)+1=K_m$ and 
$\idr(\psi(F_m))=(K_m-1)-b_{m-1}+1=K_m-b_m.$

If $b_m>b_{m-1}$,
the induction
hypothesis gives
$\idr(\psi(F^-))=K_m-b_{m-1}-1$,
and hence $1\le g=K_m-b_m+1
\le K_m-b_{m-1}.$
Thus $m$ is inserted within, or immediately after,
the initial decreasing run, which creates no occurrence of $1234$ or $12\underline{43}$; again $\psi(F_m)\in\Sym_{m}(1234,12\underline{43})$.
In this case  $\mu(\psi(F_m))=\mu(\psi(F^-))+1$ and $\idr(\psi(F_m))=g-1$,
so the induction hypothesis and the definition of $g$ yield $\mu(\psi(F_m))=K_m$ and $\idr(\psi(F_m))=K_m-b_m$.

Suppose next that
$T_{b_m}=[m;U_1,\ldots,U_r].$
Then $F^-$ has $K_m+r-1$ components, and its component rooted at
$m-1$ is $U_r$, which occupies position $b_m+r-1$. By
induction,
$\idr(\psi(F^-))=(K_m+r-1)-(b_m+r-1)=K_m-b_m.$
Since $g=K_m+1$, we have
$g>\idr(\psi(F^-))+1$
and $g\leq K_m+r-1=\mu(\psi(F^-))$.
If $\psi(F^-)$ has no $3$-rank entry, the insertion of $m$  creates 
no occurrence of $1234$ or $12\underline{43}$.
Otherwise, $m$ is inserted before, but not immediately before, the leftmost $3$-rank entry of 
$\psi(F^-)$, and again no new occurrence of $1234$ or $12\underline{43}$ arises. Hence $\psi(F_m)\in\Sym_{m}(1234,12\underline{43})$.
Moreover, $\mu(\psi(F_m))=g-1$ and $\idr(\psi(F_m))=\idr(\psi(F^-))$, from which 
 $\mu(\psi(F_m))=K_m$ and  $\idr(\psi(F_m))=K_m-b_m$, as desired.

We next prove  that $\varphi\circ\psi=\operatorname{id}_{\F_m}.$
Let $F_m\in\F_m$, and let $F^-\in\F_{m-1}$
be the forest obtained in the recursive definition of $\psi$.
By induction, $\varphi(\psi(F^-))=F^-.$
Therefore, it remains only to check that the final insertion step
recovers $F_m$.
 If $T_{b_m}=[m]$ and $b_m=b_{m-1}$, then $g=0$, which recovers the
first rule of $\varphi$. If $T_{b_m}=[m]$ and $b_m>b_{m-1}$, then  
$K_m-b_m=g-1$ components lie to the right of $T_{b_m}$, recovering the second rule.
Finally, if $T_{b_m}=[m;U_1,\ldots,U_r],$
then
$\alpha_{m-1}=\idr(\psi(F^-))+1=K_m-b_m+1$, $g=K_m+1$, $g-\alpha_{m-1}=b_m$ and $b_{m-1}=b_m+r-1$.
Hence the third rule of $\varphi$ grafts precisely
$U_1,\ldots,U_r$ under the root $m$, recovering $T_{b_m}$. 
Thus, in each case, the final step of $\varphi$ reconstructs
$F_m$.  Together with the induction hypothesis, this proves
$\varphi(\psi(F_m))=F_m,$ as desired.

It remains to prove $\psi\circ\varphi=\operatorname{id}_{\Sym_{m}(1234,12\underline{43})}.$  Let
$\pi\in\Sym_{m}(1234,12\underline{43})$, and write $\tau=\pi|_{[m-1]}$ and 
$g=\operatorname{pos}_{\pi}(m)-1$.
Let $\alpha=\idr(\tau)+1$, $F^-=\varphi(\tau)$ and $F_m=\varphi(\pi).$
By induction, $\psi(F^-)=\tau.$
If $g=0$, the first rule of $\varphi$ creates a singleton
component $[m]$ whose position equals that of the component rooted
at $m-1$ after $[m]$ is deleted.  Hence the inverse rule assigns
$g=0$.
If $1\leq g\leq\alpha$, the second rule creates a singleton
component $[m]$ with exactly $g-1$ components to its right.
Therefore, $b_m=K_m-g+1,$ i.e., $g=K_m-b_m+1.$
Since $b_m>b_{m-1}$, the inverse construction recognizes
this as the second case and recovers the same value of $g$.
Finally, if $g>\alpha$, the third rule creates a non-singleton
component rooted at $m$. A direct computation gives
$K_m=g-\alpha+K_{m-1}-b_{m-1}$, and by the induction hypothesis 
$\alpha=K_{m-1}-b_{m-1}+1$. It follows that $g=K_m+1$.
The inverse construction therefore recognizes the non-singleton
case and insert $m$ immediately after the first $K_m+1$ element of $\tau$.
 In all three cases, $\psi(F_m)=I_g(\tau)=\pi$, as desired.
 This completes the proof.
\qed

We now construct a map $\theta$ from $\F_n$ to $\PP_{n-1}$. 
Given an ordered forest $F=(T_1,\ldots,T_k),$
its preorder word $\pre(F)$ is obtained by traversing
the trees $T_1,\ldots,T_k$ in this order; within each tree, the root
is read first, followed recursively by its subtrees from left
to right. Write $\pre(F)=v_1 v_2\cdots v_n$.
For a vertex $v_i$, let
$\ell^F_<(v_i)$ and $r^F_<(v_i)$
be the nearest vertices smaller than $v_i$ to its left and right,
respectively.  If either of them does not exist, we
assign it the value $0$.  Similarly, let
$\ell^F_>(v_i)$ and $r^F_>(v_i)$
be the nearest vertices larger than $v_i$ to its left and right,
respectively, with  missing entries assigned the value $n+1$.
Define three subsets of $[n-1]$ by
\begin{align*}
S_1(F)&=\{i\in[n-1]: \ell^F_<(v_i)<r^F_<(v_i)\},\\[1pt]
S_2(F)&=\{i\in[n-1]:v_i>v_{i+1}\},\\[1pt]
S_3(F)&=\{i\in[n-1]:\ell^F_>(v_i)<r^F_>(v_i)\text{ or }v_i=n\}.
\end{align*}
For each $i\in[n-1]$, set $b_i=1$ if $i\in S_1(F)$, $m_i=1$ if $i\in S_2(F)$, and $t_i=1$ if
$i\in S_3(F)$, with all other entries equal to $0$.
This yields a triple
$\theta(F)=(P_{\mathrm{b}},P_{\mathrm{m}},P_{\mathrm{t}})$
of lattice paths of length $n-1$. 

\begin{example}\label{ex:three-components-paths}
	Take the forest $F$ from
	Example~\ref{ex:common-phi}. We have $\pre(F)=574326918$
	and the three subsets associated with $\pre(F)$ are
	\[
	\begin{aligned}
		S_1(F)&=\{1,3,4,5\},\\
		S_2(F)&=\{2,3,4,7\},\\
		S_3(F)&=\{4,5,6,7\}.
	\end{aligned}
	\]
	Equivalently,
	\[
	\renewcommand{\arraystretch}{1.15}
	\begin{array}{c|ccccccccc}
		i&1&2&3&4&5&6&7&8&9\\ \hline
		v_i&5&7&4&3&2&6&9&1&8\\ \hline
		b_i&1&0&1&1&1&0&0&0&-\\
		m_i&0&1&1&1&0&0&1&0&-\\
		t_i&0&0&0&1&1&1&1&0&-
	\end{array}
	\]
	Its corresponding non-intersecting triple of path is given in Figure \ref{fig:forest-path-example}.
	\begin{figure}[htbp]
		\centering
		\begin{tikzpicture}[x=0.8cm,y=0.8cm]
			\draw[step=1,gray!20,very thin] (0,0) grid (6,6);
			\draw[->,gray] (0,0)--(6.7,0) node[right] {$x$};
			\draw[->,gray] (0,0)--(0,6.7) node[above] {$y$};
			\node[below left,font=\scriptsize] at (0,0) {$0$};
			\foreach \i in {1,...,6}{
				\node[below,font=\scriptsize] at (\i,0) {$\i$};
				\node[left,font=\scriptsize] at (0,\i) {$\i$};
			}
			
			\draw[very thick,blue!75!black]
			plot[mark=*,mark size=1.5pt] coordinates {
				(2,0) (3,0) (3,1) (4,1) (5,1)
				(6,1) (6,2) (6,3) (6,4)
			};
			
			\draw[very thick]
			plot[mark=*,mark size=1.5pt] coordinates {
				(1,1) (1,2) (2,2) (3,2) (4,2)
				(4,3) (4,4) (5,4) (5,5)
			};
			
			\draw[very thick,red!75!black]
			plot[mark=*,mark size=1.5pt] coordinates {
				(0,2) (0,3) (0,4) (0,5) (1,5)
				(2,5) (3,5) (4,5) (4,6)
			};
			
			\node[right,text=blue!75!black] at (6,4) {$P_b$};
			\node[right] at (5,5) {$P_m$};
			\node[right,text=red!75!black] at (4,6) {$P_t$};
		\end{tikzpicture}
		\caption{The triple $T=\theta(F)\in\mathcal P_8$ encoded by
			$b=10111000$, $m=01110010$, and $t=00011110$.}
		\label{fig:forest-path-example}
	\end{figure}
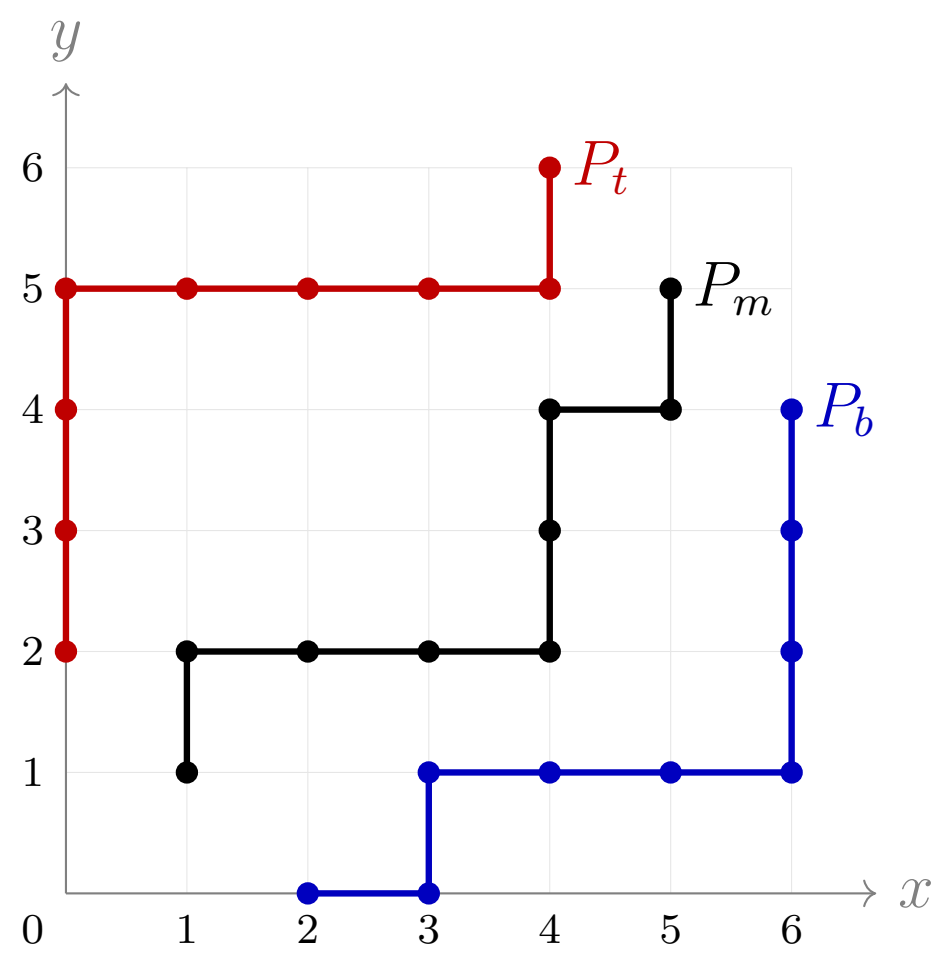
	
\end{example}

To show $\theta(F) \in \PP_{n-1}$, we need the following  properties of admissible ordered forests.

\begin{prop}\label{prop:proforest}
\begin{itemize}
  \item[(1).] Labels decrease along every edge of the forest.
 \item[(2).] The relative preorder of existing vertices is unchanged by either a leaf extension $L_p$ or a grafting operation $G_s$. 
\end{itemize}
\end{prop}

\begin{prop}\label{prop:boundarydescent}
Let $F=(T_1,\ldots,T_{M_n})\in\F_n$
and let $T_{b_n}$ be the component whose root is $n$.  If $j>b_n$, then
\[
\ro(T_{j-1})\geq\last(T_{j-1})>\ro(T_j),
\]
where $\last(T)$ denotes the last vertex of $T$ in preorder.
\end{prop}

\begin{prop}\label{prop:boundary-separation}
	Let $F=(T_1,\ldots,T_{M_n}) \in \F_n$. For any  
	$2\le p\le M_n$, put $u_p=\last(T_{p-1})$ and $w_p=\ro(T_p)$.
	Suppose that  $a$ is a vertex of one of $T_1,T_2,\ldots,T_{p-1}$
	and $c$ is a vertex of one of $T_{p},\ldots,T_{M_n}$. 
	We have
	\begin{enumerate}
	\item if $a>u_p$ and $c>u_p$, then $a<c$.
	\item if $u_p>w_p$ and $c>w_p$, then $u_p<c$.
	\end{enumerate}
\end{prop}

\pf
	We proceed by induction on the number of vertices.
	The assertion is immediate for the one-vertex admissible ordered forest. Suppose that it holds for
    $F$. We show that it holds for each forest obtained from $F$ by a leaf extension or a grafting operation.
	
	Let $T_{b_n}$ be the component rooted at $n$. 
	If $H=(A_1,\ldots,A_{M_n+1})=L_q(F)$ with $b_n \leq q \leq M_n+1$, we consider three cases. 	If $2 \leq p \leq q-1$, then for $c \neq n+1$ both conclusions follow from the induction hypothesis, while for $c=n+1,$ they are immediate.
    If $p=q$, conclusion $1$ follows from the
    induction hypothesis when $c \neq n+1$ and is immediate
    when $c=n+1$. For conclusion $2$, the hypothesis can not hold since $w_p=n+1$.
    If $q+1 \leq p \leq M_n+1$, we claim that the hypotheses of both conclusions fail. Indeed, when $p=q+1$, we have $u_p=n+1$, so no such $a$, $c$ exist for conclusion $1$. 
	For conclusion $1$ in the case $p=q+1$ and for 
	$q+1 < p \leq M_n+1$, Proposition \ref{prop:boundarydescent} gives
	$\last(A_{r-1})>\ro(A_r)$ for $r>q$.
	Combining with Proposition \ref{prop:proforest} $(1)$,  we deduce that 
	$\last(A_{p-1})>\ro(A_p)\geq y$ for every vertex $y$ in one of the trees $A_p,\ldots,A_{M_n+1}$. The claim is verified.

	Next suppose that $H=G_s(F)$ for $1\leq s<b_n$, so that
	\[H=(A_1,\ldots,A_{s+M_n-b_n})=(T_1,\ldots,T_{s-1},U,T_{b_n+1},\ldots,T_{M_n})\]
	where $U=[n+1;T_s,T_{s+1},\ldots,T_{b_n}]$.
	We consider four cases. If $s=1$, then $U$ is the first component of $H$ and contains the largest vertex $n+1$. 
	For $2 \leq p \leq s+M_n-b_n$, by Proposition \ref{prop:proforest} $(1)$
	and Proposition \ref{prop:boundarydescent}, $\last(A_{p-1})>\ro(A_p)\geq y$ for each
	vertex $y$ in the trees $A_p,\ldots,A_{s+M_n-b_n}$.
	So the hypotheses of the proposition cannot hold.
	If $s \geq 2$ and $2 \leq p \leq s-1$, it is obvious by the induction hypothesis by distinguishing whether or not $c=n+1$ as above. If $s \geq 2$ and $p=s$, conclusion $1$ follows from the induction hypothesis,
	while the hypothesis of conclusion $2$ fails since $w_p=n+1$.
	If $s \geq 2$ and $p \geq s+1$,  the same argument based on Proposition \ref{prop:boundarydescent} and Proposition \ref{prop:proforest} $(1)$ gives $\last(A_{p-1})>\ro(A_{p})\geq y$ for each vertex $y$ in one of the trees $A_p,\ldots,A_{s+M_n-b_n}$. 
	Hence no vertex
	$c$ can satisfy the hypothesis of either conclusion.
	
	We have now verified the assertion for both a leaf extension and a	grafting of $F$, which completes the induction step.
	\qed
	
 \begin{lem}\label{lem:forest-path-welldefined}
For every $F\in\mathcal F_n$, we have 
$|S_1(F)|=|S_2(F)|=|S_3(F)|$
and
\[| S_1(F)\cap[j]|\geq|S_2(F)\cap[j]|\geq| S_3(F)\cap[j]|\]
for each $0\leq j\leq n-1$. 
\end{lem}

\pf
We proceed by induction on $n$. The assertion is immediate for
$F=([1])$.
Suppose that it holds for every forest in $\F_{n-1}$; we show that it holds for every forest in $\F_n$.

Let $F=(T_1,\ldots,T_{M_n})\in\F_n$ and 
let $F^-$ be the forest on $[n-1]$ from which $F$ was generated.
Write $\pre(F^-)=v_1v_2\cdots v_{n-1}$
 and let $x$ be the position occupied by $n$ in
$\pre(F)$. By Proposition \ref{prop:proforest} $(2)$, $\pre(F)=v_1 \ldots v_{x-1} n v_x \ldots v_{n-1}$,
 where $x=n$ means that $n$ is appended at the end. Put $v_{d}=n-1$.
 For a set $X\subseteq[n-1]$, define
\[
\delta_z(X)=\{i\in X:i<z\}\cup\{i+1:i\in X,\ i\geq z\}.
\]
For brevity, write
\begin{align*}
	A^-&=S_1(F^-),\quad B^-=S_2(F^-),\quad C^-=S_3(F^-),\\
	A&=S_1(F),\quad B=S_2(F),\quad C=S_3(F).
\end{align*}
We first record how the insertion of $n$ into the preorder word affects these sets according to the following tables, where the verification for $A$ and $B$ is routine.
\begin{table}[h]
	\centering
	\vspace{-4mm}
	\[
	\begin{array}{c|c|c}
		\text{condition} & A & B\\
		\hline
		x=1\ \text{or}\ 
		\bigl(1<x\le n-1,\ v_{x-1}<v_x\bigr)
		&
		\delta_x(A^-)\cup\{x\}
		&
		\delta_x(B^-)\cup\{x\}
		\\[1mm]
		1<x\le n-1,\ v_{x-1}>v_x
		&
		\delta_x(A^-)
		&
		\bigl(\delta_x(B^-)\setminus\{x-1\}\bigr)\cup\{x\}
		\\[1mm]
		x=n
		&
		\delta_x(A^-)
		&
		\delta_x(B^-)
	\end{array}
	\] 	\vspace{-4mm}
	\caption{The recursive definitions of $A$ and $B$}
	\label{tab:A-and-B}
\end{table}
\begin{table}[h]
	\centering
	\vspace{-4mm}
\[
\begin{array}{c|c}
	\text{condition} & C\\
	\hline
	x=1\ \text{or}\ 
	\bigl(1<x\le n-1,\ v_{x-1}<v_x,\ d\ge x\bigr)
	&
	\delta_x(C^-)\cup\{x\}
	\\[1mm]
	1<x\le n-1,\ v_{x-1}<v_x,\ d<x
	&
	\text{impossible}
	\\[1mm]
	1<x\le n-1,\ v_{x-1}>v_x,\ d<x
	&
	\delta_x(C^-\setminus\{d\})\cup\{x\}
	\\[1mm]
	1<x\le n-1,\ v_{x-1}>v_x,\ d\ge x
	&
	\delta_x(C^-\setminus\{x\})\cup\{x\}
	\\[1mm]
	x=n
	&
	(C^-\cup\{n-1\})\setminus\{d\}.
\end{array}
\] 		
\vspace{-4mm}
\caption{The recursive definitions of $C$}
	\label{tab:C}
\end{table}

To Table \ref{tab:C}, we establish the following two claims. First, we claim that for any index $j$ with $j<x<n$ and $j \neq d$, we have $j \in S_3(F)$ if and only if $j \in S_3(F^-)$. Indeed, if $j \in S_3(F^-)$, since $\ell^F_>(v_j)=\ell^{F^-}_>(v_j)$ and $r^F_>(v_j) \geq r^{F^-}_>(v_j)$, so certainly $j \in S_3(F)$. 
Suppose now that $j \notin S_3(F^-)$ in which case  $r^{F^-}_>(v_j) \neq n+1$. Put $a=\ell^{F^-}_>(v_j)$.
If $a=n$,  we have $\ell^F_>(v_j)=n+1>r^F_>(v_j)$.
 If $a \neq n$ and $r^{F}_>(v_j)=r^{F^-}_>(v_j)$, again we have  $\ell^F_>(v_j)> r^{F}_>(v_j)$. Finally, suppose that  $a \neq n$ and $r^{F}_>(v_j) \neq r^{F^-}_>(v_j)$, then $\max(v_{j+1},\ldots,v_{x-1})<v_j$.
Set $c=r^{F^-}_>(v_j)$, which lies to the right of $v_{x-1}$.
Since $n$ is at position $x$ in $\pre(F)$, thus 
$v_{x-1}=\last(T)$ for some tree $T$ in $F^{-}$. By conclusion $1$ in Proposition \ref{prop:boundary-separation}, we have $a<c$, a contradiction. Consequently, $j \notin S_3(F)$ in every case, which proves the claim.

Next, we claim that for each $j$ with $x<j\le n-2$ and $j\neq d$, we have
$j\in S_3(F^-)$ if and only if $j+1\in S_3(F)$.
Since $n$ is inserted to the left of $v_j$, then $r^F_{>}(v_j)\geq r^{F^-}_{>}(v_j)$. 
Suppose first that $j\in S_3(F^-)$. Since $j\neq d$, we have
$\ell^{F^-}_{>}(v_j)<r^{F^-}_{>}(v_j)\leq n$. 
We aim to show $\ell^F_{>}(v_j)=\ell^{F^-}_{>}(v_j)$, it suffices to
show that some vertex among
$v_x,\ldots,v_{j-1}$ is larger than $v_j$.
Let $T$ be the component of $F^-$ that contains $v_j$.
If $v_j$ is not the root of $T$, then the root of $T$, which is larger than $v_j$ by $(1)$ in
Proposition~\ref{prop:proforest}, lies strictly between $n$ and $v_j$ in $\pre(F)$. 
If $v_j$ is the root of $T$, then $u=v_{j-1}$ is the last vertex of the
component preceding $T$, and $u$ lies between $n$ and $v_j$ because $j-1\geq x$.
If $u<v_j$, then $\ell^{F^-}_{>}(v_j)>u$ and $v_j>u$, so conclusion $1$ in
Proposition~\ref{prop:boundary-separation} implies  $\ell^{F^-}_{>}(v_j)<v_j$, a contradiction. Hence $u>v_j$.
In either case,
$\ell^F_{>}(v_j)=\ell^{F^-}_{>}(v_j)<r^{F^-}_{>}(v_j)\leq r^F_{>}(v_j)$, so $j+1\in S_3(F)$.
Conversely, suppose that $j\notin S_3(F^-)$. We have 
$\ell^{F^-}_{>}(v_j) \neq r^{F^-}_{>}(v_j)$ since $j \neq d$. Hence $\ell^{F^-}_{>}(v_j)>r^{F^-}_{>}(v_j)$ and
$r^F_{>}(v_j)=r^{F^-}_{>}(v_j)$. Since
$\ell^F_{>}(v_j)\geq\ell^{F^-}_{>}(v_j)$, we obtain
$\ell^F_{>}(v_j)>r^F_{>}(v_j)$ and $j+1\notin S_3(F)$. This proves the claim.

By the two claims above, to determine the relation between $C$ and $C^-$,
it remains only to consider the elements $n-1$ and $v_x$. 
We omit most of the details, which are routine, and record only the facts needed to justify the equality of the cardinality increments obtained from the two tables.
First, if $d<x<n$, then $d\in C^-.$
Indeed, $v_d=n-1$ is the maximum vertex of $F^-$. Since
$d<x\leq n-1,$ the vertex $v_d$ is not the terminal entry of
$\pre(F^-).$ Thus, $d\in C^-.$
Second, suppose that $1<x\leq n-1$, $v_{x-1}>v_x$ and $d\geq x$, then $x\in C^-.$ To see this, note first that $v_x \neq n-1$, so $d>x$, and the nearest larger vertex to the right of $v_x$ exists. Let $r_{>}^{F^-}(v_x)=v_r.$ 
Since $v_{x-1}>v_x$ and $v_{x-1}$ lies immediately to the left of $v_x$, we have $\ell_{>}^{F^-}(v_x)=v_{x-1}.$ 
Proposition \ref{prop:boundary-separation} $2$ gives
$v_{x-1}<v_r$. Hence $x \in C^-$.
Finally, whenever $1<x\leq n-1$ and $v_{x-1}>v_x,$ we have
$x-1\in B^-.$
Combining these three facts and the two tables for $A,B$ and $C$, we have $|A|-|A^-|=|B|-|B^-|=|C|-|C^-|$.
By induction, we obtain $|A|=|B|=|C|$.

One case in the table above still requires justification.
If $1 <x \leq n-1$, $v_{x-1}<v_x$ and $d<x$, then $v_{x-1}=\last(T_{j-1})$ and $v_x=\ro(T_j)$ for some $T_j$ lying to the right of the tree rooted at $n-1$, contradicting
Proposition \ref{prop:boundarydescent}. Hence this case cannot occur.

We now prove the two prefix inequalities.  Define
\begin{align*}
D_j=|A\cap[j]|-|B\cap[j]| \quad \text{and} \quad
D_j^-=|A^-\cap[j]|-|B^-\cap[j]|. 
\end{align*}
If $x=1$, then
$D_j=D^-_{j-1}$ for $1\leq j\leq n-1$.
If $1<x\leq n-1$, then
\[
D_j=
\begin{cases}
	D^-_j,&0\leq j<x-1,\\
	D^-_{x-1},&j=x-1 \,\,\text{and}\,\,  v_{x-1}<v_{x},\\
	D^-_{x-1}+1,&j=x-1 \,\,\text{and}\,\,  v_{x-1}>v_{x},\\
	D^-_{j-1},&x\leq j\leq n-1.
\end{cases}
\]
If $x=n$, then $D_j=D^-_j$ for $0\leq j\leq n-1$.
By the induction hypothesis, all these quantities are nonnegative, and hence
$|A\cap[j]|\geq |B\cap[j]|$ for $0\leq j\leq n-1$.

It remains to compare the middle and top sets.  Put
\begin{align*}
E_j=|B\cap[j]|-|C\cap[j]| \quad \text{and} \quad E_j^-=|B^-\cap[j]|-|C^-\cap[j]|.
\end{align*}
If $x=1$, then
$E_j=E_{j-1}^-$ for $1\leq j\leq n-1$.
Suppose that $1<x\le n-1$, $v_{x-1}<v_x$ and $d\ge x$;
then
\begin{equation}\label{eq:Ecase1}
E_j=
\begin{cases}
	E_j^-,   & 0\le j<x,\\
	E_{j-1}^-, & x\le j\le n-1.
\end{cases}
\end{equation}
Suppose that $1<x\leq n-1$, $v_{x-1}>v_x$ and $d<x$;
then
\begin{equation}\label{eq:Ecase2}
E_j=
\begin{cases}
	E_j^-,       & 0\le j<d,\\
	E_j^-+1,     & d\le j<x-1,\\
	E_{x-1}^-,   & j=x-1,\\
	E_{j-1}^-,   & x\le j\le n-1.
\end{cases}
\end{equation}
Suppose that $1<x\leq n-1$, $v_{x-1}>v_x$ and $d\geq x$;
then
\begin{equation}\label{eq:Ecase3}
E_j=
\begin{cases}
	E_j^-,             & 0\leq j<x-1,\\
	E_{x-1}^- -1,      & x-1\leq j\leq x,\\
	E_{j-1}^-,         & x<j\leq n-1.
\end{cases}
\end{equation}
Finally, if $x=n$, then
\begin{equation}\label{eq:Ecase4}
E_j=
\begin{cases}
	E_j^-,       & 0\leq j<d,\\
	E_j^-+1,     & d\leq j<n-1,\\
	E_{n-1}^-,   & j=n-1.
\end{cases}
\end{equation}

By the induction hypothesis, all the terms on the right-hand sides are
nonnegative, except possibly the term
$E_{x-1}^- -1$ in the case $v_{x-1}>v_x$ and $d\geq x$.  In this case,
we claim that every element of $C^-\cap[x-1]$ can be matched with a distinct element of $B^-\cap[x-2]$. 
Set $i=x-1$. We construct an injection from $C^-\cap [i]$ to $B^-\cap [i-1].$
Since $d\geq x=i+1,$ the exceptional element corresponding to the maximum vertex is not contained in $C^-\cap [i].$
Thus every element $k\in C^-\cap [i]$ satisfies
$\ell_{>}^{F^-}(v_k)<r_{>}^{F^-}(v_k)$,
and, in particular, has a larger vertex to its left.
For such $k$, let $\lambda(k)=\max\{t<k:v_t>v_k\}.$
By definition, every vertex between
$v_{\lambda(k)}$ and $v_k$ 	is smaller than $v_k$.
Consequently, $v_{\lambda(k)}>v_{\lambda(k)+1},$
and hence $	\lambda(k)\in B^-.$
Moreover, $\lambda(k)<k\leq i,$	so $\lambda(k)\in B^-\cap[i-1].$
	
It remains to prove that $\lambda$ is injective. Suppose, to the contrary, that there exist $k,l \in C^-\cap [i]$ such that $k<l$ and $\lambda(k)=\lambda(l)=q$. Then 
$v_q>v_l>v_k$. Let $r$ be the position of the nearest larger vertex to the right of $v_k$. Since $v_l>v_k$ and $l$ lies to the right of $k$, we have $r \leq l$.
Because $k\in C^-,$ we have $v_q<v_r$, and hence $v_q<v_l$, a contradiction. Thus $\lambda$ is injective.

Consequently,
$|C^-\cap[x-1]|\leq|B^-\cap[x-2]|.$
On the other hand, $v_{x-1}>v_x$ implies that $x-1\in B^-$,  so $|B^-\cap[x-1]|=|B^-\cap[x-2]|+1.$
Combining these two relations gives
$$
\begin{aligned}
	E_{x-1}^-
	&=
	|B^-\cap[x-1]|-|C^-\cap[x-1]|\\
	&=
	|B^-\cap[x-2]|+1-|C^-\cap[x-1]|\\
	&\geq 1.
\end{aligned}
$$
It follows that
$|B\cap [j]|\ge |C\cap [j]|$ for $0\leq j\leq n-1.$
This completes the proof. \qed

We now describe the inverse construction $\alpha:\PP_{n-1} \rightarrow \F_n$ by induction.
Let 
$P=(P_{\mathrm b},P_{\mathrm m},P_{\mathrm t})\in\PP_{n-1},$
and let $A,B,C\subseteq[n-1]$
be the position sets of the $H$-steps of the three paths. We construct
a forest $F=\alpha(P)$. For $n=1$, sent the unique element of $\PP_0$
to $([1])$.
Assume now that $n\geq 2$.
Define the critical position
\[
x=\min\Bigl(\{\,i\in C:E_{i-1}=0\,\}\cup\{n\}\Bigr).
\]
For a set $X\subseteq[n]$, define 
\[\kappa_z(X)=\{i\in X:i<z\}\cup\{i-1:i\in X,\ i>z\}.\]
We define three sets $A^-$, $B^-$, and $C^-$ by the following three cases.
\begin{enumerate}
	\item If $x\in A$, set
	$A^-=\kappa_x(A)$, $B^-=\kappa_x(B)$ and 	$C^-=\kappa_x(C)$.
	
	\item Suppose that $x\notin A$ and that $x=n$, $x=n-1$ or $x+1\in C$. Let
	\[d=\max\{i\in[x-1]:E_{i-1}=0\}.\]
	If $x<n$, define 
	\[A^-=\kappa_x(A),\quad B^-=\kappa_x(B\cup\{x-1\}), \quad C^-=\kappa_x(C\cup\{d\}).\]
	If $x=n$, define 
	\[A^-=A,\quad B^-=B, \quad C^-=(C\cup\{d\})\setminus\{n-1\}.\]
	
	\item Suppose that $x\notin A$, $x<n-1$ and $x+1\notin C$.	Define
	\[A^-=\kappa_x(A),\quad B^-=\kappa_x(B\cup\{x-1\}), \quad C^-=\kappa_x(C\cup\{x+1\}).\]
\end{enumerate}
For each $i\in[n-2]$, set $b^-_i=1$ if $i\in A^-$, $m^-_i=1$ if $i\in B^-$, and $t^-_i=1$ if
$i\in C^-$, with all other entries equal to $0$.
This yields a triple
$P^-=(P_{\mathrm{b}^-},P_{\mathrm{m}^-},P_{\mathrm{t}^-})$
of lattice paths of length $n-2$. 
Lemma \ref{lem:path-delete} below shows that $P^- \in \PP_{n-2}$.
Set $F^-=\alpha(P^-)=(T_1,\ldots,T_K)$
and let $T_r$ be the unique component whose root is $n-1$.
If $x=n$, set $q=K+1.$
If $x<n$, consider the $x$-th vertex in 
$\pre(F^-)$.  Lemma \ref{lem:root-gap} below implies that this vertex is the root of a
unique component, say $T_q$.
We then define
\[
F=\alpha(P)
=
\begin{cases}
	L_q(F^-),& r\leq q\leq K+1,\\[1mm]
	G_q(F^-),& 1\leq q<r.
\end{cases}
\]

 \begin{lem}\label{lem:path-delete}
Given $P\in\mathcal P_{n-1}$ with $n \geq 2$, we have $P^-\in\mathcal P_{n-2}.$
\end{lem}

\pf
	Put $N=n-1$.  Since $P\in\mathcal P_N$, we have
	$D_j\geq 0$, $E_j\geq 0$ for $0\leq j\leq N$ with
	$D_N=E_N=0$. For $X\subseteq[N]$ and $x\leq N$, a direct count gives
\begin{align}\label{eq:card}
	\bigl|\kappa_x(X)\cap[j]\bigr|
	=
	\begin{cases}
		|X\cap[j]|,&0\leq j<x,\\[1mm]
		|X\cap[j+1]|-\mathbf 1_{\{x\in X\}},
		&x\leq j\leq N-1.
	\end{cases}
\end{align}
	Before treating the cases, we record three facts about the critical position $x$. First, if $x<n$, then $x\in C$ and $E_{x-1}=0$.  Since
	$E_x=E_{x-1}+\mathbf 1_{\{x\in B\}}
	-\mathbf 1_{\{x\in C\}},$ we necessarily have
    $x\in B.$ Second, if $x\notin A$, then $x\geq 2$.  Indeed, if $x=1$, then $1\in C$, and the prefix inequalities imply $1\in B$ and $1\in A$, a contradiction. 
    Finally, $x-1\notin B$ for every $x$.
    This is immediate when
    $x=1$.  If $x>1$, then $E_{x-1}=0$: this follows from the definition
    of $x$ when $x<n$, and from $|B|=|C|$ when $x=n$.
    If $x-1\in B$, then $0=E_{x-1}=E_{x-2}+\mathbf 1_{\{x-1\in B\}}-\mathbf 1_{\{x-1\in C\}}$
	forces $x-1\in C$ and $E_{x-2}=0$, contradicting the definition of $x$.
	
	We now verify the three cases in the construction of $\alpha$.
	
	\noindent\textbf{Case 1: $x\in A$.}
	
	Here $x<n$, so $x\in B$ and $x\in C$.  Thus $x$ belongs to all three sets.  Formula (\ref{eq:card}) gives
	\[
	D_j^-=
	\begin{cases}
		D_j,&0\le j<x,\\
		D_{j+1},&x\le j\le N-1,
	\end{cases}
	\qquad
	E_j^-=
	\begin{cases}
		E_j,&0\le j<x,\\
		E_{j+1},&x\le j\le N-1.
	\end{cases}
	\]
	Hence $D_j^-\geq0$ and $E_j^-\geq0$ for every $0\leq j\leq N-1$.  Moreover,
	$|A^-|=|A|-1$, $|B^-|=|B|-1$ and $|C^-|=|C|-1$, so the three cardinalities remain equal.

	\noindent\textbf{Case 2a: $x\notin A$, $x<n$, and either $x=n-1$ or $x+1\in C$.}
	
	In this case $x\notin A$, $x\in B\cap C$ and
	$x-1\notin B$.  Since $x \geq 2$ and $E_0=0$, we have $d=\max\{i\in[x-1]:E_{i-1}=0\}$
	is well defined.  Since $d<x$, the minimality of $x$ implies $d\notin C.$ Recall that
	\[
	A^-=\kappa_x(A),\qquad
	B^-=\kappa_x(B\cup\{x-1\}),\qquad
	C^-=\kappa_x(C\cup\{d\}).
	\]
	Using (\ref{eq:card}), together with
	$x\in B\cap C$, $x\notin A$, $x-1\notin B$ and
	$d\notin C$, we obtain
	\[
	D_j^-=
	\begin{cases}
		D_j,&0\le j<x-1,\\
		D_{x-1}-1,&j=x-1,\\
		D_{j+1},&x\le j\le N-1,
	\end{cases}
    \qquad
	E_j^-=
	\begin{cases}
		E_j,&0\le j<d,\\
		E_j-1,&d\le j\le x-2,\\
		0,&j=x-1,\\
		E_{j+1},&x\le j\le N-1.
	\end{cases}
	\]
	Since $x\in B$ and $x\notin A$,	we have $D_x=D_{x-1}-1\ge0,$ so $D_{x-1}\ge1$.  By the maximality of $d$, we have
	$E_j\geq 1$ for $d\leq j\leq x-2$. It follows that
	$D_j^-\geq 0$ and $E_j^-\geq 0$ for $0\le j\le N-1$.

	The cardinalities are $|A^-|=|A|,$	$|B^-|=|B|+1-1=|B|$
	and $|C^-|=|C|+1-1=|C|$.
	Thus the three sets again have the same cardinality.
	
	\noindent\textbf{Case 2b: $x=n$.}
	
	Here $d=\max\{i\in[N]:E_{i-1}=0\}.$
	We first show that $N\notin A\cup B.$
	If $N\in B \setminus C$, then
	$E_N=E_{N-1}+1>0,$ which is impossible.  If $N\in B\cap C$, then $E_N=E_{N-1}=0,$
	contradicting the definition of $x=n$.  Hence $N\notin B$. If $N\in A$, then $D_N=D_{N-1}+1>0$ since $N\notin B$, again impossible. 
	Not also that  $d\notin C$, because $d<n$ and $E_{d-1}=0$.
	Recall that 
	\[
	A^-=A,\qquad
	B^-=B,\qquad
    C^-=(C\cup\{d\})\setminus\{N\}.
	\]
	If $d=N$, then $N\notin C$, so $C^-=C.$
	If $d<N$, then $N\in C$.  Otherwise,
	since $N \notin B$, the equality $0=E_{N}=E_{N-1}+\mathbf 1_{\{N\in B\}}-\mathbf 1_{\{N\in C\}}$ would give 
	$E_N=E_{N-1}=0$ and hence $d=N$.  Hence, when
	$d<N$, the set $C^-$ is obtained from $C$ by replacing $N$	with $d$.  In either case, we have 
	$|A^-|=|B^-|=|C^-|.$
	
	Clearly, $D_j^-=D_j\geq0$ for $0\leq j\leq N-1$.
	If $d=N$, then also $E_j^-=E_j\geq0$.  If $d<N$,
	then
	\[
	E_j^-=
	\begin{cases}
		E_j,&0\le j<d,\\
		E_j-1,&d\le j\le N-1.
	\end{cases}
	\]
	By the maximality of $d$, we have $E_j\geq1$ for $d\leq j\leq N-1$. Therefore, $E_j^-\geq 0$ for $0\leq j\leq N-1$.
	
	\noindent\textbf{Case 3: $x \notin A$, $x<n-1$ and $x+1\notin C$.}
	
	As before, we have $x\in B\cap C$, $x\notin A$ and $x-1\notin B$.
	Recall that
	\[
	A^-=\kappa_x(A),\qquad
	B^-=\kappa_x(B\cup\{x-1\}),\qquad
	C^-=\kappa_x(C\cup\{x+1\}).
	\]
	Applying (\ref{eq:card}) again, we have
	\[
	D_j^-=
	\begin{cases}
		D_j,&0\le j<x-1,\\
		D_{x-1}-1,&j=x-1,\\
		D_{j+1},&x\le j\le N-1,
	\end{cases}
   \qquad
	E_j^-=
	\begin{cases}
		E_j,&0\le j<x-1,\\
		E_{x-1}+1=1,&j=x-1,\\
		E_{j+1},&x\le j\le N-1.
	\end{cases}
	\]
	As in Case 2a, $D_{x-1}\geq 1$, so  
	$D_j^-\geq 0$ and $E_j^-\geq 0$ for $0\le j\leq N-1$.
	Finally, $	|A^-|=|A|$, $|B^-|=|B|+1-1=|B|$ and $|C^-|=|C|+1-1=|C|$.

	In all cases, the three sets \(A^-,B^-,C^-\) are subsets of $[n-2]$ of equal cardinality satisfying
	\[
	|A^-\cap[j]|
	\ge
	|B^-\cap[j]|
	\ge
	|C^-\cap[j]|
	\]
	for $0\leq j\leq n-2.$
	Consequently, $P^-\in\PP_{n-2}.$ This completes the proof.
\qed
 
\begin{prop}\label{prop:admissible-insertion-positions}
	Let $F=(T_1,\ldots,T_{M_n})\in\F_n$, and let $T_{b_n}$ be	the component rooted at $n$.  In the recursive construction of	$\F_{n+1}$, the new vertex $n+1$ is inserted either immediately before the root of some component or in the final gap of
	$\pre(F)$.  More precisely, for the two operations 
	\[
	\begin{array}{ll}
		G_s(F),&1\leq s<b_n,\\
		L_p(F),&b_n\leq p\leq M_n+1,
	\end{array}
	\]
	the first one inserts $n+1$ immediately before the root of $T_s$, while the second one inserts $n+1$ immediately before the root of $T_p$ when $p\leq M_n$, or in the final gap when $p=M_n+1$.
\end{prop}
This follows directly from the recursive definition of admissible ordered forests together with
Proposition \ref{prop:proforest} $(2)$.

\begin{lem}\label{lem:root-gap}
Given any $P\in\PP_{n-1}$ with $n \geq 2$, let $F^-=\alpha(P^-)$.  If $x<n$, then the $x$-th vertex of $\pre(F^-)$ is the root of some component in $F^-$.
\end{lem}
\pf
For $P\in\PP_{n-1}$, write $F=\alpha(P)$. Let $A,B,C,A^-,B^-,C^-,x,E_j,E^-_j$ be defined for $P$ and $P^-$ as above. Set $H_0=H_{n+1}=0$ and, for $1\leq i\leq n$,
\[
H_i=E_{i-1}+\mathbf 1_{\{i\geq x\}}
-\mathbf 1_{\{i\in C\cup\{n\}\}}.
\]
Define
\[
\RootPos(F)
=
\{\pos_{\pre(F)}(\ro(T)):
T\text{ is a component of }F\}.
\]
We prove simultaneously, by induction on $n$,  that
\begin{equation}\label{eq:root-zero}
r\in\RootPos(F) \Longleftrightarrow H_{r-1}=H_r=0 \,\,\text{for each $1\leq r\leq n$}
\end{equation} 
 and that $x=\pos_{\pre(F)}(n).$

The case $n=1$ is immediate. Assume that the assertions hold for $n-1$.
By Lemma~\ref{lem:path-delete}, $P^-\in\PP_{n-2}$. Put
\[a=\min\Bigl(\{\,i\in C^-:E^-_{i-1}=0\,\}\cup\{n-1\}\Bigr) \,\,\text{and}\,\,
H^-_i=E^-_{i-1}+\mathbf 1_{\{i\geq a\}}
-\mathbf 1_{\{i\in C^-\cup\{n-1\}\}}.
\]
To relate $H_i$ and $H^-_i$, we first identify $a$ with $d$ introduced in case $2$ of the construction of $\alpha$.
Suppose first that $x \notin A$, $x<n$ and either $x=n-1$ or $x+1 \in C$. Since $d<x$ and $E_{d-1}=0$, the minimality of $x$ implies $d\notin C$. By the definition of $C^-$ and the relations between $E_i$ and $E^-_i$ established in case 2a of Lemma \ref{lem:path-delete} , we have $d\in C^-$ and $E_{d-1}^-=E_{d-1}=0.$
Moreover, if $i<d$ and $i\in C^-$, then
$i\in C$ and $E^-_{i-1}=E_{i-1}$. The equality
$E^-_{i-1}=0$ would contradict the definition of $x$.
Therefore, $a=d$ in this case.
If $x=n$, then $C^-=(C\cup\{d\})\setminus\{n-1\},$ and
the same argument shows that  $E^-_{i-1}\neq0$ whenever $i<d$ and $i\in C^-$.
When $d<n-1$, then $d\in C^-$ and
$E^-_{d-1}=0$; when $d=n-1$, the definition of $a$ gives $a=n-1$. In both subcases $a=d$, and
by the induction hypothesis, $a=d$ is the position of $n-1$ in $\pre(F^-)$.

Using formula (\ref{eq:card}), a direct substitution gives the following two possibilities:
\begin{equation}\label{eq:H1}
	H_i=
	\begin{cases}
		H_i^-,&0\leq i<x,\\
		H_{i-1}^-,&x\leq i\leq n+1,
	\end{cases}
	\qquad\text{if }a\leq x,
\end{equation}
and
\begin{equation}\label{eq:H2}
H_i=
\begin{cases}
	H_i^-,&0\leq i\leq x,\\
	H_{i-1}^-+1,&x<i\leq a,\\
	H_{i-1}^-,&a<i\leq n+1,
\end{cases}
\qquad\text{if }x<a.
\end{equation}
The identities \eqref{eq:H1} and \eqref{eq:H2}
follow by direct substitution in the definitions of $H_i$ and $H^-_i$, the relations between $E_i$ and $E^-_i$ established in Lemma \ref{lem:path-delete} and the definition of $C^-$ in the construction of $\alpha$.
Since the
verification is purely algebraic, we omit the case-by-case calculation.

We next prove that
\begin{equation}\label{eq:Hx0}
H_{x-1}=H_x=0.
\end{equation}
Suppose first that $x<n$, so that $x\in C$ and
$E_{x-1}=0$. If $x=1$, the assertion follows directly from the definition of $H_0$ and $H_1$. Assume $x>1$. We claim that $x-1\notin B$. Indeed, if $x-1\in B \setminus C$, then
$E_{x-1}\geq 1$. If $x-1\in B\cap C$, then
$E_{x-2}=0$, contradicting the minimality of $x$. Hence $E_{x-2}=\mathbf 1_{\{x-1\in C\}}$, and consequently
\[H_{x-1}=E_{x-2}-\mathbf 1_{\{x-1\in C\}}=0, \quad H_x=E_{x-1}+1-\mathbf 1_{\{x\in C\}}=0.\]
If $x=n$, then there is no $i \in C$ with $E_{i-1}=0$.
Since $E_{n-1}=E_{n-2}+\mathbf 1_{\{n-1\in B\}}-\mathbf 1_{\{n-1\in C\}}$, the possibility $n-1\in B\cap C$ would give $E_{n-2}=0$, contradicting the
definition of $x=n$, while $n-1\in B \setminus C$ would give $E_{n-1} \geq 1$. Hence $n-1 \notin B$, so
$E_{n-2}=\mathbf 1_{\{n-1\in C\}},$  and therefore
\[H_{n-1}=E_{n-2}-\mathbf 1_{\{n-1\in C\}}=0, \quad H_n=E_{n-1}+1-\mathbf 1_{\{n\in\{n\}\}}=0.\]
This proves (\ref{eq:Hx0}).

We now prove that
\begin{equation}\label{eq:H-x0}
H^-_{x-1}=H^-_x=0.
\end{equation}
If $x<a$, then \eqref{eq:H2} together with
\eqref{eq:Hx0}  immediately gives
$H^-_{x-1}=H_{x-1}=0$ and $H^-_x=H_x=0.$
It remains to consider the case $a\leq x$. 
By (\ref{eq:H1}), $H^-_{x-1}=H_x=0$, so we need only show
that $H^-_x=0$. We distinguish three cases.
If $x\in A$, assume to the contrary that $a<x$. Then for some $i<x$, we have $i\in C^-$ and $E^-_{i-1}=0$. Since
$C^-=\kappa_x(C)$ and $E^-_{i-1}=E_{i-1}$, 
this would give $i\in C$ and $E_{i-1}=0$, contradicting the
definition of $x$. Hence $a=x$.
Since $x\in A$, $x<n$, and the relations between $E_i$ and $E^-_i$ established in case $1$ of Lemma \ref{lem:path-delete} gives $E^-_{x-1}=E_{x-1}=0.$
If $x<n-1$, then $a=x$ and the definition of $a$ imply that
$x\in C^-$. If $x=n-1$, then trivially
$x\in\{n-1\}$. In either case,
$x\in C^-\cup\{n-1\}$, and consequently
$H^-_x=E^-_{x-1}+\mathbf 1_{\{x\geq a\}}-\mathbf 1_{\{x\in C^-\cup\{n-1\}\}}=0.$

Suppose next that $x\notin A$ and either $x=n-1$ or $x+1\in C$. Through the proof of case 2a  in Lemma \ref{lem:path-delete}, we have $E^-_{x-1}=0$ and $x\in C^-\cup\{n-1\}$, so
$H^-_x=E^-_{x-1}+1-\mathbf 1_{\{x\in C^-\cup\{n-1\}\}}=0.$
If $x=n$, then $H^-_n=0$ holds by definition.

Suppose finally that $x \notin A$, $x<n-1$ and $x+1\notin C$. By definition, $B^-=\kappa_x(B\cup\{x-1\})$ and 
$C^-=\kappa_x(C\cup\{x+1\})$.
For every $i<x$, $i\in C^-$ if and only if $i\in C$, and 
$E^-_{i-1}=E_{i-1}$. Thus, $a \geq x$.
Moreover, $x\in C^-$ while $E^-_{x-1}=1,$
so $a \neq x$, and hence $a>x.$
But then this case falls under the situation $x<a$ treated above, for which \eqref{eq:H-x0} has already been established.
Thus \eqref{eq:H-x0} holds in all cases.

If $x<n$, then (\ref{eq:H-x0}) and the induction hypothesis show that the $x$-th vertex of $\pre(F^-)$ is a component root. If $x=n$, the final gap is allowed by
Proposition~\ref{prop:admissible-insertion-positions}. Hence the recursive construction of $F$ is well-defined.

We now determine the position of the vertex $n$.
Suppose first that $x<n$.  By the
definition of $\alpha$ and Proposition~\ref{prop:proforest} $(2)$, both operations $L_q(F^-)$ and $G_q(F^-)$
insert $n$ immediately before the root of $T_q$ in $\pre(F^-)$.
If $x=n$, then $n$ is appended at
the end of $\pre(F^-)$.
In either case, $x=\pos_{\pre(F)}(n).$

Finally suppose that $a\leq x$. The vertex $n$ is inserted by a leaf extension, so
\[
\RootPos(F)
=
\bigl(\RootPos(F^-)\cap[x-1]\bigr)
\cup\{x\}
\cup
\{r+1:r\in\RootPos(F^-),\ r\geq x\}.
\]
By \eqref{eq:root-zero} applied to $F^-$, together with (\ref{eq:H1}) and (\ref{eq:Hx0}), the zero
pairs of the sequence $(H_i)$ have exactly the same transformation.
Therefore \eqref{eq:root-zero} holds for $F$.

Suppose next that $x<a$.  The vertex $n$ is inserted by a
grafting operation, which merges the components whose roots occur at position $x,x+1,\ldots,a-1$ or $a$; consequently,
\[
\RootPos(F)
=
\bigl(\RootPos(F^-)\cap[x-1]\bigr)
\cup\{x\}
\cup
\{r+1:r\in\RootPos(F^-),\ r\geq a+1\}.
\]
Formula (\ref{eq:H2}) shows that the zero pairs at positions
$x+1,\ldots,a$ disappear, while all zero pairs beginning at a position at least $a+1$ shift one position to the right.
Together with \eqref{eq:Hx0}, this establishes \eqref{eq:root-zero} in this case as well.
This completes the induction.
\qed

The proof of the lemma above immediately yields the following corollary.
\begin{cor}\label{cor:x-n-position}
We have $x=\pos_{\pre(F)}(n)$ and $d=\pos_{\pre(F^-)}(n-1)$.  
\end{cor}

\begin{lem}\label{lem:forest-path-bijection}
The map $\theta:\F_n\rightarrow\PP_{n-1}$
is a bijection, with inverse 
$\alpha: \PP_{n-1}\rightarrow\F_n.$
\end{lem}

\pf
It suffices to show that $\alpha \circ \theta=\operatorname{id}_{\F_n}$ and $\theta \circ \alpha=\operatorname{id}_{\PP_{n-1}}$.

We first prove that $\alpha \circ \theta=\operatorname{id}_{\F_n}$ by induction on $n$.
The case $n=1$ is immediate.
Given $F \in \F_n$, let $F^-$ be the unique forest from which $F$ is generated. Set $A=S_1(F)$, $B=S_2(F)$, $C=S_3(F)$ and $P=\theta(F)$. Let $y$ be the position of $n$ in $\pre(F)$. We will show that $y=x$ and $P^-=\theta(F^-)$, where $P^-$ is defined in the construction of $\alpha$.

We first prove that $y=x$.  Write $\pre(F^-)=v_1v_2\cdots v_{n-1}$, and let $d=\pos_{\pre(F^-)}(n-1)$, and set
$\widetilde A=S_1(F^-)$, $\widetilde B=S_2(F^-)$ and  
$\widetilde C=S_3(F^-)$.
For $0\leq j\leq n-2$, put
$\widetilde E_j=|\widetilde B\cap[j]|-|\widetilde C\cap[j]|.$
Define $\widetilde H_i=\widetilde E_{i-1}+\mathbf 1_{\{i\geq d\}}
-\mathbf 1_{\{i\in\widetilde C\cup\{n-1\}\}}.$
Since $\alpha(\theta(F^-))=F^-$ by induction, through the proof of Lemma \ref{lem:root-gap}, we have 
\begin{equation}\label{eq:d-critical}
	d=
	\min\Bigl(
	\{i\in\widetilde C:\widetilde E_{i-1}=0\}
	\cup\{n-1\}
	\Bigr).
\end{equation}
and
\[
r\in\RootPos(F^-)
\quad\Longleftrightarrow\quad
\widetilde H_{r-1}=\widetilde H_r=0.
\]
If $y<n$, Proposition~\ref{prop:admissible-insertion-positions}
shows that $v_y$ is the root of a component of $F^-$, and hence
\begin{equation}\label{eq:Hy-root}
	\widetilde H_{y-1}=\widetilde H_y=0.
\end{equation}
We now apply Tables~\ref{tab:A-and-B} and~\ref{tab:C}, with the insertion position denoted by $y$.

If $y=1$, then $1\in C$ and $E_0=0$, so $x=1=y$.

Suppose next that $1<y<n$ and $v_{y-1}<v_y$.  The impossible
case in Table~\ref{tab:C} shows that $d\geq y$.  If $d=y$, then $\widetilde E_{y-1}=0$ by \eqref{eq:d-critical}; this remains true  when $d=y=n-1$, since
$|\widetilde B|=|\widetilde C|$.  If $d>y$, then
\eqref{eq:Hy-root} implies $0=\widetilde H_{y-1}=\widetilde E_{y-2}-\mathbf1_{\{y-1\in\widetilde C\}}.$
Since $v_{y-1}<v_y$, we have $y-1\notin\widetilde B$, and therefore
$\widetilde E_{y-1}=\widetilde E_{y-2}-\mathbf1_{\{y-1\in\widetilde C\}}=0.$
In either subcase, Tables \ref{tab:A-and-B} and \ref{tab:C} give $y\in C$ and $E_{y-1}=\widetilde E_{y-1}=0$.
Moreover, for every $i<y$, we have
$i\in C$ if and only if $i\in\widetilde C$ and 
$E_{i-1}=\widetilde E_{i-1}.$
Since $i<y\leq d$, \eqref{eq:d-critical} implies that no such $i$ satisfies both $i\in C$ and $E_{i-1}=0$.  Hence $x=y$.

Now suppose that $1<y<n$, $v_{y-1}>v_y$ and $d<y$.  From (\ref{eq:Hy-root}) we have
$0=\widetilde H_y=\widetilde E_{y-1}+1-\mathbf 1_{\{y\in\widetilde C\cup\{n-1\}\}}.$
Since $\widetilde E_{y-1}\geq0$, it follows that
$\widetilde E_{y-1}=0$ and $y\in\widetilde C\cup\{n-1\}$.  
Tables \ref{tab:A-and-B} and \ref{tab:C}, and (\ref{eq:Ecase2}) then yield
$y\in C$ and  $E_{y-1}=\widetilde E_{y-1}=0$.
For $i<d$, since $i\in C$ if and only if $i\in\widetilde C$ and $E_{i-1}=\widetilde E_{i-1}$,  condition  \eqref{eq:d-critical} again implies that no such $i$ satisfies $i \in C$ and $E_{i-1}=0$. 
For $d<i<y$, Tables \ref{tab:A-and-B} and \ref{tab:C}, and (\ref{eq:Ecase2})  give
$E_{i-1}=\widetilde E_{i-1}+1>0.$ Since $d \notin C$, we conclude that $x=y$.

Suppose next that $1<y<n$, $v_{y-1}>v_y$ and $d\geq y$.  In fact $d>y$, since $d=y$ would give $v_y=n-1$, contradicting $v_{y-1}>v_y$.  As $y-1<d$, \eqref{eq:Hy-root} gives
$\widetilde E_{y-2}=\mathbf 1_{\{y-1\in\widetilde C\}}.$
Moreover, $v_{y-1}>v_y$ implies $y-1\in\widetilde B$, and   hence $\widetilde E_{y-1}=\widetilde E_{y-2}+1-\mathbf 1_{\{y-1\in\widetilde C\}}=1.$
Tables \ref{tab:A-and-B} and \ref{tab:C}, and (\ref{eq:Ecase3}) then give
$y\in C$ and $E_{y-1}=\widetilde E_{y-1}-1=0$.
For every $i<y$, we have $i\in C$ if and only if $i\in\widetilde C$ and $E_{i-1}=\widetilde E_{i-1}$.
Since $i<y<d$, no such $i$ satisfies both $i \in C$ and $E_{i-1}=0$. Thus $x=y$.

Finally, suppose that $y=n$.  By the tables and (\ref{eq:Ecase4}), for $i<d$ we have $i\in C$ if and only if $i\in\widetilde C$ and $E_{i-1}=\widetilde E_{i-1}$.
In view of \eqref{eq:d-critical}, none of these $i$ satisfies $i\in C$ and $E{i-1}=0$.
Moreover, $d \notin C$, and if $d<i<n$, then from (\ref{eq:Ecase4}) $E_{i-1}=\widetilde E_{i-1}+1>0.$
Consequently, there is no $i<n$ with
$i\in C$ and $E_{i-1}=0$.  By definition, $x=n=y$.
Therefore $x=y$ in all cases.

We next prove that $P^-=\theta(F^-)$, that is, $A^-=S_1(F^-)$, $B^-=S_2(F^-)$ and $C^-=S_3(F^-)$, where $x,A^-,B^-$ and $C^-$ are defined in the construction of $\alpha$. Since $x=y$, we may apply Tables~\ref{tab:A-and-B}
and~\ref{tab:C} with the insertion position denoted by $x$.
For every $X\subseteq[n-2]$, we have
\begin{equation}\label{eq:kappa-delta-inverse}
	\kappa_x(\delta_x(X))
	=
	\kappa_x(\delta_x(X)\cup\{x\})
	=
	X.
\end{equation}
We compare the three cases in the definition of $\alpha$ with the corresponding rows of Tables \ref{tab:A-and-B} and \ref{tab:C}.

\noindent\textbf{Case 1: $x\in A$.}
By Table~\ref{tab:A-and-B}, this occurs precisely when $x=1$, or when $1<x<n$ and $v_{x-1}<v_x$.  The impossible case in Table~\ref{tab:C} shows that, in the latter situation $d\geq x$.
Consequently, the two tables give
\[
A=\delta_x(\widetilde A)\cup\{x\},\qquad
B=\delta_x(\widetilde B)\cup\{x\},\qquad
C=\delta_x(\widetilde C)\cup\{x\}.
\]
The first case of $\alpha$, together with
\eqref{eq:kappa-delta-inverse}, therefore gives
\[
A^-=\kappa_x(A)=\widetilde A,\qquad
B^-=\kappa_x(B)=\widetilde B,\qquad
C^-=\kappa_x(C)=\widetilde C.
\]

\noindent\textbf{Case 2: $x\notin A$ and either $x=n$, $x=n-1$, or $x+1\in C$.}
Consider first the subcase $x<n$.  This subcase is characterized by $v_{x-1}>v_x$ and $d<x$.
Indeed, the preceding proof of $x=y$ shows that
$x\in\widetilde C\cup\{n-1\}$.  Hence if $x<n-1$, then
$x\in\widetilde C$, and Table~\ref{tab:C} gives $x+1\in C$;
if $x=n-1$, the second case of $\alpha$ applies directly.
Thus,  
\begin{equation*}
	A=\delta_x(\widetilde A), \quad
	B=\bigl(\delta_x(\widetilde B)\setminus\{x-1\}\bigr)
	\cup\{x\},\quad
	C=\delta_x(\widetilde C\setminus\{d\})
	\cup\{x\}.
\end{equation*}
Let $d_\alpha=\max\{i\in[x-1]:E_{i-1}=0\}$.  Relation (\ref{eq:Ecase2}) gives
$E_{d-1}=\widetilde{E}_{d-1}=0$ and 
$E_{i-1}=\widetilde E_{i-1}+1>0$ for $d<i<x$.
It follows that $d_\alpha=d$.  Therefore,
\begin{equation*}
	A^-	=\kappa_x(A)=\widetilde A, \quad
	B^-	=\kappa_x(B\cup\{x-1\})	=\widetilde B, \quad
	C^-=\kappa_x(C\cup\{d_\alpha\})	=\kappa_x(C\cup\{d\})
	=\widetilde C.
\end{equation*}
We now consider the subcase $x=n$.  Tables \ref{tab:A-and-B} and \ref{tab:C} give 
\[
A=\widetilde A,\qquad
B=\widetilde B,\qquad
C=(\widetilde C\cup\{n-1\})\setminus\{d\}.
\]
As above, relation (\ref{eq:Ecase4}) gives
$E_{d-1}=\widetilde{E}_{d-1}=0$ and $E_{i-1}>0$ for $d<i<n$,
so $d_\alpha=d$. The second case of $\alpha$ gives
\[
A^-=A=\widetilde A,\quad
B^-=B=\widetilde B, \quad
C^-=(C\cup\{d\})\setminus\{n-1\}
=\widetilde C.
\]

\noindent\textbf{Case 3: $x\notin A$, $x<n-1$, and
	$x+1\notin C$.}
This case is characterized by $v_{x-1}>v_x$ and $d>x$.
Indeed, the proof of $x=y$ gives
$x\in\widetilde C$, so Table~\ref{tab:C} shows that $x+1\notin C$.  Since $d>x$, we also have $x<n-1$.
Tables \ref{tab:A-and-B} and \ref{tab:C} give
\begin{equation*}
	A=\delta_x(\widetilde A), \quad
	B=\bigl(\delta_x(\widetilde B)\setminus\{x-1\}\bigr)
	\cup\{x\}, \quad
	C=\delta_x(\widetilde C\setminus\{x\})\cup\{x\}.
\end{equation*}
Hence, $B\cup\{x-1\}=\delta_x(\widetilde B)\cup\{x\}$
and $C\cup\{x+1\}=\delta_x(\widetilde C)\cup\{x\}.$
The third case of $\alpha$ and \eqref{eq:kappa-delta-inverse} therefore give
\begin{equation*}
A^-=\kappa_x(A)=\widetilde A, \quad
B^-=\kappa_x(B\cup\{x-1\})=\widetilde B, \quad
C^-=\kappa_x(C\cup\{x+1\})=\widetilde C.
\end{equation*}

In all three cases, $A^-=S_1(F^-)$,
$B^-=S_2(F^-)$ and $C^-=S_3(F^-)$.
Consequently, $P^-=\theta(F^-).$

Finally, the induction hypothesis then gives $\alpha(P^-)=F^-.$ write $F^-=(T_1,\ldots,T_K)$, and let $T_r$ be the component rooted at $n-1$. If $y<n$, let $T_q$ be the component whose root occupies position $y$ in $\pre(F^-)$; if $y=n$, set $q=K+1$.  Since $x=y$ and $\alpha(P^-)=F^-$, the construction of $\alpha(P)$ uses $L_q$ when $r\leq q$ and $G_q$ when $q<r$. These are exactly the two operations
in the construction of $F$ from $F^-$. Hence
$\alpha(P)=F$, as desired.

The identity $\theta\circ\alpha=\operatorname{id}_{\PP_{n-1}}$ follows  by a similar induction: Corollary \ref{cor:x-n-position} identifies the parameters $x$ and $d$ used by $\alpha$ with the corresponding positions in the resulting forests, and Tables~\ref{tab:A-and-B} and~\ref{tab:C} then recover the original sets $A$,  $B$, and $C$. We omit the analogous details.
\qed

\section{Proofs of Theorems \ref{thm:setpardecom} and \ref{thm:setcore-baxter}}
\label{sec:partition-decomposition-baxter}

In this section, we first introduce the blue-chain contraction to decompose $\N_n$ into a $2$-colored set partition core and a subset of $[n-1]$. This decomposition yields a proof of Theorems \ref{thm:setpardecom}. We then construct a bijection 
between these $2$-colored set partition cores and non-intersecting triples of Lattice paths, which proves Theorem \ref{thm:setcore-baxter}.

It is worth noting that the constructions developed in this section were inspired by the reduction
algorithm proposed by Chen, Deng and Du \cite{ChenDengDu2005}. The algorithm replaces each arc $(i,j)$ by $(i,j-1)$ and then deletes the last vertex.  
When adapted to our $2$-colored framework, this idea motivates us to first contract maximal
chains of blue adjacent arcs, yielding a colored core and a cut set. We then apply the same endpoint shift to the core during the right-end normalization. 

Given $\eta\in\N_n$, let $G_{\mathrm b}(\eta)$ denote the graph on $[n]$ whose edges consist of the blue adjacent arcs of $\eta$. Let $I_1<I_2<\cdots<I_{k+1}$ stand for the vertex sets of its connected components, ordered from left to right.  Since every edge of $G_{\mathrm b}(\eta)$ takes the form $(i,i+1)$, each $I_j$ forms a consecutive interval, and these intervals partition  $[n]$.  We now define the \emph{blue-chain contraction}: contract each interval $I_j$ to a single vertex, while preserving
every arc whose endpoints lie in different intervals together with its original color. Write $\core(\eta)$
for the resulting colored partition of $[k+1]$, and define its cut set via
\begin{equation*}
 S(\eta)=\{\max I_1,\ldots,\max I_k\}\subseteq[n-1].
\end{equation*}

Conversely, we introduce the \emph{expansion} operation. Fix $S=\{s_1<\cdots<s_k\}\subseteq[n-1]$ and
$c\in\DD_{k}$.  The set $S$ determines the intervals
\[
I_1=[1,s_1],\quad I_2=[s_1+1,s_2],\quad\ldots,\quad
I_{k+1}=[s_k+1,n].
\]
In the special case $k=0$, we have $I_1=[1,n]$. Replace each vertex $j$ of
$c$ by the interval $I_j$, connect consecutive elements within $I_j$ using blue adjacent arcs, and replace every arc $(r,s)$ of $c$ by an arc from
$\max I_r$ to $\min I_s$ with its original color.

\begin{lem}\label{lem:partition-contraction-core}
The contracted partition $\core(\eta)$ belongs to $\DD_{k}$.
\end{lem}

\pf
     We first verify that $\core(\eta)$ is noncrossing.
	Suppose, to the contrary, that $\core(\eta)$ contains 
	a pair of crossing arcs sharing the same color, denoted 
	$(i,k)$ and $(j,l)$ with $i<j<k<l$.
	The arc $(i,k)$ of $\core(\eta)$ arises from an arc of $\eta$ connecting $\max I_i$ to $\min I_k$.
	Similarly, $(j,l)$ corresponds to an arc of $\eta$ from $\max I_j$ to $\min I_l$.
	Since these intervals are ordered consecutively and
	$i<j<k<l$, we have	$\max I_i<\max I_j<\min I_k<\min I_l.$
	It follows that  $(\max I_i,\min I_k)$ and $(\max I_j,\min I_l)$ cross in $\eta$ with the same color, contradicting the assumption that $\eta$ is $2$-colored	noncrossing.  Hence $\core(\eta)$ is $2$-colored noncrossing.

  We next prove that $\core(\eta)$ has no blue adjacent arc. Assume for contradiction that $(i,i+1)$ is such an arc. Then $(\max I_i, \min I_{i+1})$ would be a blue  adjacent arc in $\eta$. In that case, the two intervals $I_i$ and $I_{i+1}$ would  belong to the same maximal blue chain, a contradiction. This completes the proof. \qed

\begin{lem}
\label{lem:partition-binomial-decomposition}
For every $n\geq1$, the blue-chain contraction defines a bijection from $\N_n$ to
\[\bigsqcup_{k=0}^{n-1}\binom{[n-1]}{k}\times\DD_{k}.\]
\end{lem}

\pf 
	It suffices to conform that blue-chain contraction and the	expansion operation defined above are mutual inverses. The endpoint order in
	the expansion preserves the noncrossing property for each color, and	the construction clearly produces a $2$-colored set partition.

	Under expansion, each core arc $(r,s)$ is replaced by an arc from
	$\max I_r$ to $\min I_s$. Therefore, any arc between two distinct intervals always takes the right endpoint of the left interval and the left endpoint of the right interval; contraction therefore recovers both the core and the cut set.  Conversely, for the original partition $\eta$, any arc leaving a
	blue-chain interval must start at final vertex, and any arc entering such an interval must terminate at its first vertex.
	The condition $\core(\eta)\in\mathcal D_{k}$ ensures
	that the recovered intervals are exactly the maximal blue chains. So applying the expansion operation to $\bigl(S(\eta),\operatorname{core}(\eta)\bigr)$ reconstructs $\eta$.
\qed

\begin{example}
A $2$-colored noncrossing set partition  $\eta$ and its blue-chain contraction are given in the following figure.
	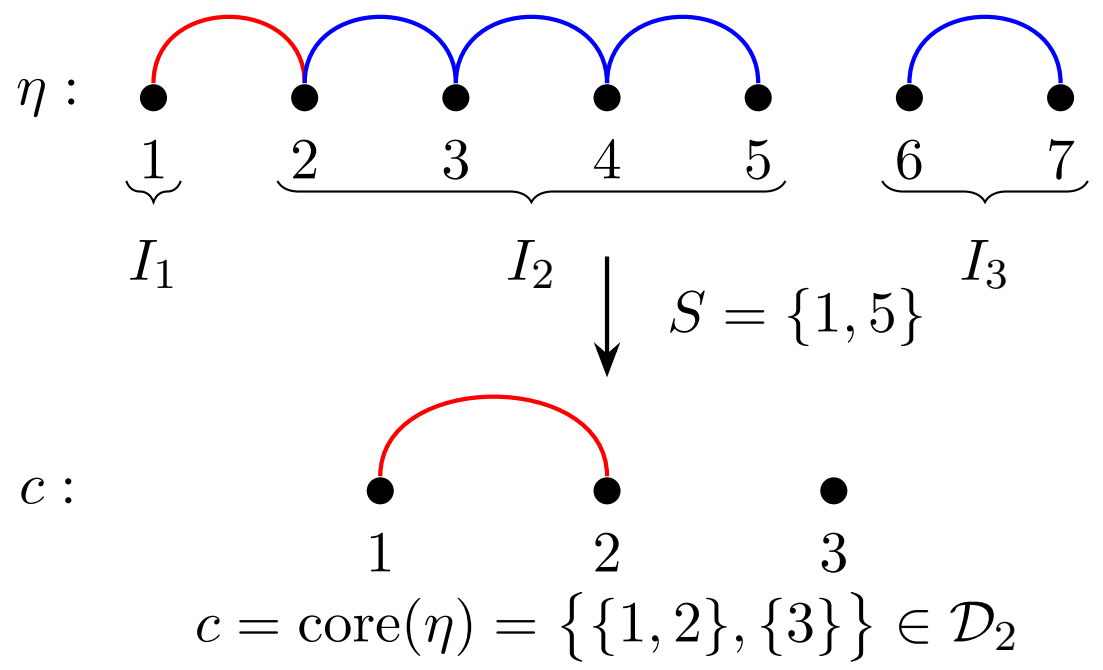
\begin{figure}[htbp]
		\centering
		\begin{tikzpicture}[
			x=1cm,y=1cm,
			vertex/.style={circle,fill=black,inner sep=1.8pt},
			redarc/.style={red,thick},
			bluearc/.style={blue,thick},
			every node/.style={font=\small},
			>=Stealth
			]
			\node at (-0.7,0) {$\eta:$};
			
			\foreach \i in {1,...,7}{
				\node[vertex] (v\i) at ({\i-1},0) {};
				\node[below=4pt] at (v\i) {$\i$};
			}
			
			\draw[redarc]
			(v1) to[out=90,in=90,looseness=1.5] (v2);
			
			\foreach \i/\j in {2/3,3/4,4/5,6/7}{
				\draw[bluearc]
				(v\i) to[out=90,in=90,looseness=1.5] (v\j);
			}
			
			\draw[decorate,decoration={brace,mirror,amplitude=4pt}]
			(-0.18,-0.55) -- (0.18,-0.55)
			node[midway,below=7pt] {$I_1$};
			
			\draw[decorate,decoration={brace,mirror,amplitude=4pt}]
			(0.82,-0.55) -- (4.18,-0.55)
			node[midway,below=7pt] {$I_2$};
			
			\draw[decorate,decoration={brace,mirror,amplitude=4pt}]
			(4.82,-0.55) -- (6.18,-0.55)
			node[midway,below=7pt] {$I_3$};
			
			\draw[->,thick] (3,-1.05) -- (3,-1.85)
			node[midway,right=7pt] {$S=\{1,5\}$};
			
			\node at (-0.7,-2.6) {$c:$};
			
			\foreach \i/\x in {1/1.5,2/3,3/4.5}{
				\node[vertex] (u\i) at (\x,-2.6) {};
				\node[below=4pt] at (u\i) {$\i$};
			}
			
			\draw[redarc]
			(u1) to[out=90,in=90,looseness=1.2] (u2);
			
			\node at (3,-3.5)
			{$c=\operatorname{core}(\eta)
				=\bigl\{\{1,2\},\{3\}\bigr\}\in\mathcal{D}_2$};
			
		\end{tikzpicture}
		\caption{A blue-chain contraction of $\eta$.}
		\label{fig:blue-chain-example}
	\end{figure}
	Conversely, $S=\{1,5\}$ determines the intervals
	$I_1=[1,1]$, $I_2=[2,5]$ and $ I_3=[6,7].$
	Expand each vertex $j$ of $c$ into $I_j$ and join
	consecutive vertices within each interval by blue arcs.
	The red core arc $(1,2)$ becomes the red arc
	$(\max I_1,\min I_2)=(1,2).$
	This recovers the original colored partition $\eta$.
\end{example}

We now construct an explicit bijection from $\DD_n$ to
$\PP_{n}$. Take any $c\in\DD_n$. We substitute each colored arc $(i,j)$ of $c$ with a directed edge from $i$ to $j-1$ with the same color, and then delete the vertex $n+1$.
Call the resulting graph on $[n]$ the \emph{right-end normalization} $G(c)$.  A red adjacent arc turns into a red loop, while the definition
of $\DD_n$ forbids blue loops.
At any vertex of $G(c)$, there exists at most one incoming and at most one outgoing nonloop edge.  If both occur, 
they must carry different colors:
a pair of edges of the same color would arise from crossing arcs in $c$.
Consequently, each vertex of $G(c)$ has exactly one of the  eight states below.  We encode the state at vertex $i$ via the column vector $(b_i,m_i,t_i)$, which is given in Table \ref{tab:setpar2triplepath}.
\begin{table}[ht]
	\centering
\[
\begin{array}{c|c|c@{\quad\quad}c|c|c} 
status\, number &\text{state}	& (b_i,m_i,t_i) & status \,number & \text{state}& (b_i,m_i,t_i)\\
	\hline
	1&
	\stateisolated
	&000&
	2
	&
	\stateredloop
	&111
	\\[2mm]
	3&
	\stateredout
	&100&
	4
	&
	\stateredin
	&011
	\\[2mm]
	5&
	\stateblueout
	&110&
	6
	&
	\statebluein
	&001
	\\[2mm]
	7&
	\stateredinblueout
	&010&
	8
	&
	\stateblueinredout
	&101
\end{array}
\]
\vspace{-4mm}
\caption{The eight states of the right-end normalization and corresponding binary words}
\label{tab:setpar2triplepath}
\end{table}
Traversing vertices in left-to-right order yields three binary words
$b=b_1\cdots b_n$, $m=m_1\cdots m_n$ and $t=t_1\cdots t_n$.

\begin{lem}\label{lem:colored-core-three-words}
For each integer $0\leq j\leq n$, we have
\begin{equation}\label{eq:8stategeq}
 \sum_{i=1}^{j}b_i\geq
 \sum_{i=1}^{j}m_i\geq
 \sum_{i=1}^{j}t_i,
\end{equation}
and 
\begin{equation*}
	\sum_{i=1}^{n}b_i=
	\sum_{i=1}^{n}m_i=
	\sum_{i=1}^{n}t_i,
\end{equation*}
\end{lem}

\pf 
For $1\leq r\leq8$, let $N_r(j)$ denote the number of vertices with status $r$ among the first $j$ vertices.  By the table above,
\[
\sum_{i=1}^{j} b_i
=
N_2(j)+N_3(j)+N_5(j)+N_8(j),
\]
while
\[
\sum_{i=1}^{j} m_i
=
N_2(j)+N_4(j)+N_5(j)+N_7(j).
\]
Therefore,
\[
\sum_{i=1}^{j}b_i-\sum_{i=1}^{j}m_i
=
N_3(j)+N_8(j)-N_4(j)-N_7(j).
\]
The first difference in (\ref{eq:8stategeq}) equals the
number of red nonloop edges that have started but not yet terminated among the first $j$ vertices of $G(c)$.
By a parallel analysis, the second difference equals the
corresponding number of open blue edges.  Both quantities are nonnegative.  At the
end every  edge has been closed, so both differences vanish and the equality follows. \qed

Set $1$ to be an $H$-step and $0$ an $N$-step. Let $\gamma(c)$ be the triple of lattice path, where $b,m,t$
are the bottom, middle, and top rows.  Lemma \ref{lem:colored-core-three-words} shows that $\gamma(c) \in \PP_{n}$.

\begin{example}\label{ex:three-components-partition}
	The left diagram in Figure \ref{fig:partition-normalization} shows $c\in\D_8$, and the right diagram shows its right-end normalization $G(c)$.
\begin{figure}[htbp]
	\centering
	\begin{tikzpicture}[
		x=0.65cm,y=0.75cm,
		vertex/.style={circle,fill=black,inner sep=1.5pt},
		every node/.style={font=\small}
		]
		\begin{scope}
			\node[left] at (-0.45,0) {$c:$};
			
			\foreach \i in {1,...,9}{
				\node[vertex] (c\i) at ({\i-1},0) {};
				\node[below=4pt] at (c\i) {$\i$};
			}
			
			\draw[red,thick]
			(c1) to[out=90,in=90] (c3);
			\draw[red,thick]
			(c4) to[out=90,in=90] (c5);
			\draw[red,thick]
			(c5) to[out=90,in=90] (c8);
			
			\draw[blue,thick]
			(c2) to[out=90,in=90] (c7);
			\draw[blue,thick]
			(c3) to[out=90,in=90] (c6);
		\end{scope}
		
		\begin{scope}[shift={(11,0)}]
			\node[left] at (-0.45,0) {$G(c):$};
			
			\foreach \i in {1,...,8}{
				\node[vertex] (g\i) at ({\i-1},0) {};
				\node[below=4pt] at (g\i) {$\i$};
			}
			
			\draw[red,thick]
			(g1) to[out=90,in=90] (g2);
			\draw[red,thick]
			(g4) edge[loop above,min distance=8mm] (g4);
			\draw[red,thick]
			(g5) to[out=90,in=90] (g7);
			
			\draw[blue,thick]
			(g2) to[out=90,in=90] (g6);
			\draw[blue,thick]
			(g3) to[out=90,in=90] (g5);
		\end{scope}
	\end{tikzpicture}
	\caption{The arc diagram of $c$ and its right-end
		normalization $G(c)$.}
	\label{fig:partition-normalization}
\end{figure}
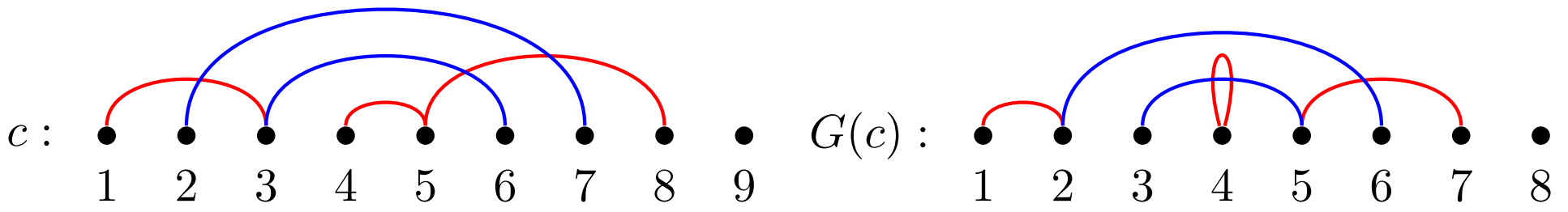
 The state table is
	\[
	\renewcommand{\arraystretch}{1.2}
	\begin{array}{c|cccccccc}
		i&1&2&3&4&5&6&7&8\\ \hline
		\text{status number}&3&7&5&2&8&6&4&1\\
		(b_i,m_i,t_i)&100&010&110&111&101&001&011&000
	\end{array}
	\]
	Reading the three rows gives
    $b=10111000$, $m=01110010$ and $t=00011110$.
	$\gamma(c)$ is exactly the path given in Figure \ref{fig:forest-path-example}. 
\end{example}

\begin{lem}
\label{lem:colored-core-path-bijection}
$\gamma$ is a bijection from $\DD_n$ to $\PP_{n}$.
\end{lem}

\pf We first construct the inverse of $\gamma$. Given a triple of non-intersecting lattice path $P \in \PP_{n}$, we extract the step sequences of its bottom, middle and top paths, which are written as $b=b_1\cdots b_n$, $m=m_1\cdots m_n$ and $t=t_1\cdots t_n$,
respectively. Let $c=\beta(P)$ be given through the following three steps.

 First, for each $1\leq i\leq n$, assign vertex $i$ a state specified by $(b_i,m_i,t_i)$ in Table \ref{tab:setpar2triplepath}.

Second, connect the half-edges with the same color as follows.  Scan the vertices from left to right. Whenever a red incoming half-edge is
encountered, connect it to the nearest unmatched red outgoing
half-edge to its left. The inequality
$\sum_{i=1}^{j}b_i\geq\sum_{i=1}^{j}m_i$
for every $1 \leq j \leq n$
guarantees that such an outgoing half-edge always exists.  Similarly, each blue incoming half-edge is matched to the nearest unmatched blue outgoing half-edge to its left; this is ensured by
$\sum_{i=1}^{j}m_i\geq\sum_{i=1}^{j}t_i$
for every $1 \leq j \leq n$.
The equality of the total row sums ensures that all non-loop incoming
and outgoing half-edges are paired.  One readily verifies that the resulting arcs of each color are noncrossing

Finally, replace each edge $(i,j)$ with $(i,j+1)$ and adjoin the vertex $n+1$.
Since every half-edge is paired and each state of the vertex has at most one incoming and at most one outgoing half-edge, the resulting graph is
the linear representation of a $2$-colored set partition.
The right-end shift introduces no crossing: 
by Table \ref{tab:setpar2triplepath},
no vertex simultaneously carries an incoming and an outgoing half-edge with the same color, except for a loop.
Combined with the fact that the resulting graph in the second step is noncrossing with each color, the arcs of each color of the $2$-colored set partition generated  in this step remain noncrossing. 
The absence of blue loops gives no blue
adjacent arc.  Hence the resulting partition belongs to $\DD_n$.  

Both the right-end normalization and
Table \ref{tab:setpar2triplepath} can be reversed uniquely, so $\gamma$ and $\beta$ are inverses to each other. \qed

\section*{Acknowledgement}
This work was supported by the National Natural Science Foundation of China (No. 12271301) .

\section*{Declaration of AI use}
AI tools were used to assist in exploring proof strategies,
performing calculations, and polishing the language.
All AI-assisted arguments and suggestions were independently
verified and revised as necessary by the author, who takes
full responsibility for the content of the article.

\end{document}